\documentclass{amsart}
\usepackage{amssymb, amsmath, amscd, latexsym, mathrsfs, appendix}
\usepackage{amsmath}
\usepackage[mathscr]{eucal} 
\usepackage{amsfonts}
\usepackage[all]{xy}

\numberwithin{equation}{section}
\newtheorem{thm}[equation]{Theorem} 
\newtheorem{cor}[equation]{Corollary}
\newtheorem{lem}[equation]{Lemma}
\newtheorem{conv}[equation]{Conventions}
\newtheorem{notn}[equation]{Notation}
\newtheorem{prop}[equation]{Proposition}
\newtheorem{defn}[equation]{Definition}
\newtheorem{con}[equation]{Construction}
\newtheorem{rem}[equation]{Remark}
\newtheorem{expl}[equation]{Example}

\newcommand{\set}[1]{\left\{#1\right\}}

\newcommand{\lra}{\longrightarrow}
\newcommand{\co}{\colon\!}
\newcommand{\blank}{\textup{---}}
\newcommand{\ld}{^{\textup{ld}}} 
\newcommand{\lf}{^{\textup{lf}}} 

\newcommand{\smin}{\smallsetminus}
\newcommand{\pt}{\star\,}
\newcommand{\ubul}{^{\bullet}}
\newcommand{\lbul}{_{\bullet}}

\newcommand{\nin}{\noindent}
\newcommand{\cp}{{}^c}
\newcommand{\upr}{{}^{\%}}
\newcommand{\lpr}{{}_{\%}}
\newcommand{\thf}{{}^{th\ZZ/2}}
\newcommand{\hf}{{}^{h\ZZ/2}}
\newcommand{\ho}{{}_{h\ZZ/2}}
\newcommand{\cc}{\,;\,c}

\newcommand{\red}{\textup{rd}}
\newcommand{\Str}{\textup{Str}}

\newcommand{\cone}{\textup{cone}}
\newcommand{\id}{\textup{id}}

\newcommand{\holim}{\textup{holim}}
\newcommand{\hocolim}{\textup{hocolim}}
\newcommand{\colim}{\textup{colim}}
\newcommand{\fiber}{\textup{fiber}}
\newcommand{\hofiber}{\textup{hofiber}}
\newcommand{\mor}{\textup{mor}}

\newcommand{\map}{\textup{map}}
\newcommand{\homeo}{\mathscr{H}\textit{om}}
\newcommand{\Homeo}{\textup{\textsf{Hom}}}
\newcommand{\emb}{\mathscr{E}\textit{mb}}
\newcommand{\Emb}{\textsf{Emb}}
\newcommand{\rel}{\textup{rel }}
\newcommand{\st}{\textup{st}}
\newcommand{\Fr}{\textup{Fr}}

\newcommand{\Wh}{\textup{Wh}}

\newcommand{\TOP}{\textup{TOP}}
\newcommand{\G}{\textup{G}}
\newcommand{\bLA}{\mathbf L\mathbf A}
\newcommand{\bLK}{\mathbf L\mathbf K}
\newcommand{\bVL}{\mathbf V\mathbf L} 
\newcommand{\bVLA}{\mathbf V\mathbf L\mathbf A}
\newcommand{\Th}{\mathbf T\mathbf h}
\newcommand{\Bo}{\mathbf B\mathbf m}

\newcommand{\bA}{\mathbf A}
\newcommand{\bC}{\mathbf C}
\newcommand{\bE}{\mathbf E}
\newcommand{\bF}{\mathbf F}
\newcommand{\bH}{\mathbf H}
\newcommand{\bK}{\mathbf K}
\newcommand{\bL}{\mathbf L}
\newcommand{\bP}{\mathbf P}
\newcommand{\bS}{\mathbf S}

\newcommand{\bV}{\mathbf V}

\newcommand{\sA}{\mathcal A}
\newcommand{\sB}{\mathcal B}
\newcommand{\sC}{\mathcal C}
\newcommand{\sD}{\mathcal D}
\newcommand{\sE}{\mathcal E}
\newcommand{\sG}{\mathcal G}
\newcommand{\sH}{\mathcal H}
\newcommand{\sI}{\mathcal I}
\newcommand{\sK}{\mathcal K}
\newcommand{\sM}{\mathcal M}
\newcommand{\sP}{\mathcal P}
\newcommand{\sQ}{\mathcal Q}
\newcommand{\sR}{\mathcal R}
\newcommand{\sS}{\mathcal S}
\newcommand{\sW}{\mathcal W}

\newcommand{\NN}{\mathbb N}
\newcommand{\RR}{\mathbb R}
\newcommand{\ZZ}{\mathbb Z}
\newcommand{\JJ}{\mathbb J}

\newcommand{\ee}{\'{e} }

\newcommand{\colimsub}[1]{\begin{array}[t]{cc} \textup{colim} \\
[-1.7mm] \scriptstyle{#1} \end{array}}

\newcommand{\hocolimsub}[1]{\begin{array}[t]{cc} \textup{hocolim} \\
[-1.7mm] \scriptstyle{#1} \end{array}}

\date{23 Feb 2012}
\begin{document}
\title[WWIII]{Automorphisms of manifolds and algebraic $K$--theory: Part III }%
\author{Michael S. Weiss and E. Bruce Williams}%
\address{Math. Institut, Universit\"at M\"unster, Einsteinstr. 62, 48149 M\"unster, Germany} \email{m.weiss@uni-muenster.de}%
\address{Dept. of Math., University of Notre Dame, 255 Hurley, Notre Dame, IN 46556, USA} \email{williams.4@nd.edu}
\thanks{}%

\begin{abstract}
The structure space $\sS(M)$ of a closed topological $m$-manifold
$M$ classifies bundles whose fibers are closed $m$-manifolds
equipped with a homotopy equivalence to $M$. We construct a highly
connected map from $\sS(M)$ to a concoction of algebraic $L$-theory
and algebraic $K$-theory spaces associated with $M$. The
construction refines the well-known surgery theoretic analysis of
the block structure space of $M$ in terms of $L$-theory.
\end{abstract}
\maketitle

\tableofcontents
\section{Introduction} \label{sec-intro}
The structure space $\sS(M)$ of a closed topological $m$-manifold
$M$ is the classifying space for bundles $E\to X$
with an arbitrary $CW$-space $X$ as base, closed topological manifolds as fibers
and with a fiber homotopy trivialization
\[ E\to M\times X \]
(a homotopy equivalence and a map over $X$).
The points of $\sS(M)$ can loosely be imagined as pairs
$(N,f)$ where $N$ is a closed $m$-manifold and $f\co N\to M$ is a
homotopy equivalence.
To explain the relationship between $\mathcal S(M)$ and automorphisms of $M$, we invoke
$\homeo(M)$, the topological group of homeomorphisms from $M$ to $M$, and $\G(M)$, the grouplike topological
monoid of homotopy equivalences from $M$ to $M$. In practice we work with simplicial models of $\homeo(M)$ and
$\G(M)$. The homotopy fiber of the inclusion $B\homeo(M)\to B\G(M)$
is homotopy equivalent to a union of connected components of $\sS(M)$.

\medskip
The main result of this paper is a calculation of the homotopy type
of $\sS(M)$ in the so-called concordance stable range, in terms of
$L$- and algebraic $K$-theory. With $m$ fixed as above, we construct a
homotopy invariant functor
\[ (Y,\xi) \mapsto \bLA\lbul\lpr(Y,\xi,m) \]
from spaces $Y$ with spherical fibrations $\xi$ to spectra. The spectrum
$\bLA\lbul\lpr(Y,\xi,m)$ is a concoction of the $L$-theory and the
algebraic $K$-theory of spaces \cite{WaldRutgers} associated with $Y$,
compounded with an assembly construction \cite{RanickiTopMan}. (The subscript
$\lpr$ is for homotopy fibers of assembly maps.)
In the case where $Y=M$ (nonempty and connected for simplicity)
and $\xi$ is $\nu$, the normal fibration of $M$,
there is a ``local degree'' map
\[ \Omega^{\infty+m}\bLA\lbul\lpr(M,\nu,m)\lra  8\ZZ\subset\ZZ. \]
There is then a highly connected map
\begin{equation} \label{eqn-mainmap}
\sS(M) \lra
\fiber\big[\CD\,\Omega^{\infty+m}\bLA\lbul\lpr(M,\nu,m) @>\textup{local deg.\,\,}>\rule{0mm}{1mm}> 8\ZZ\endCD\,\big]
\end{equation}
where $\fiber$ in this case means the fiber over $0\in 8\ZZ$, an infinite loop space.
The connectivity estimate is
given by the concordance stable range. In practice that translates into $m/3$
approximately, but in theory it is more convoluted and the reader is referred to
definition~\ref{defn-concordstable}. \newline
The result has a generalization to the case in which $M$ is compact with nonempty
boundary. It looks formally the same. Points of $\sS(M)$ can be imagined as pairs
$(N,f)$ where $N$ is a compact manifold with boundary and
$f\co(N,\partial N)\to (M,\partial M)$ is a homotopy equivalence of pairs restricting to
a homeomorphism of $\partial N$ with $\partial M$.

\medskip
We now give a slightly more detailed, although still sketchy,
definition of the spectrum $\Omega^m\bLA_{\bullet\%}(Y,\xi,m)$. (Details can be found in
section~\ref{sec-hocharsig}.) It is the total homotopy fiber of
a commutative square
\begin{equation} \label{eqn-defLA}
\CD
\Omega^m\bL_\bullet\upr(Y,\xi) @>>> S^1\wedge \bA^\%(Y,\xi,m)_{h\ZZ/2} \\
@VVV   @VVV \\
\Omega^m\bL_\bullet(Y,\xi) @>>> S^1\wedge \bA(Y,\xi,m)_{h\ZZ/2}~.
\endCD
\end{equation}
The left-hand column is the quadratic $L$-theory assembly map. (This has a variety of equivalent descriptions and in
the simplest of these it depends only on $Y$ and the orientation double covering $w_\xi$ of $Y$ determined by $\xi$.
We do not insist on such a simple description because that would make the horizontal maps in the diagram more
obscure; therefore $\bL\lbul(Y,\xi)$ rather than $\bL\lbul(Y,w_\xi)$ is the notation which we prefer.) The right-hand column
is the Waldhausen $A$-theory assembly map with $S^1\wedge$ and homotopy orbit construction inflicted. Both columns use a category
of (finitely dominated) retractive spaces or spectra over $Y$, subject to finiteness conditions
and equipped with a notion of Spanier-Whitehead duality which depends on $\xi$ and $m$.
The horizontal maps are variants of a natural transformation $\Xi$ which was defined in \cite{WW2}; the precise
relationship will be clarified in a moment. By opting for finitely dominated retractive spaces in both columns
we have implicitly selected the decoration $p$. (The infinite loop space
$\Omega^{\infty+m}\bLA_{\bullet\%}(Y,\xi,m)$ is decoration independent, i.e.,
any consistent choice of decoration from the list $h,p,\dots,\langle -i\rangle,...,\langle -\infty\rangle$ gives
the same result up to a homotopy equivalence.)
\newline
This is a definition which relates $\Omega^m\bLA_{\bullet\%}(Y,\xi,m)$ to known and trusted concepts
in algebraic L- and K-theory. For our constructions we prefer another definition of
$\Omega^m\bLA_{\bullet\%}(Y,\xi,m)$
as the homotopy fiber of the map between homotopy pullbacks of the rows in the commutative diagram
\begin{equation} \label{eqn-pracdefLA}
\CD
\Omega^m\bVL^{\bullet\%}(Y,\xi) @>\Xi>> \bA^\%(Y,\xi,m)^{th\ZZ/2} @<\textup{incl.}<< \bA^\%(Y,\xi,m)^{h\ZZ/2} \\
@VVV   @VVV  @VVV \\
\Omega^m\bVL^\bullet(Y,\xi) @>\Xi >> \bA(Y,\xi,m)^{th\ZZ/2} @<\textup{incl.}<<  \bA(Y,\xi,m)^{h\ZZ/2}
\endCD
\end{equation}
where the left-hand column is the assembly map in a form of \emph{visible symmetric} $L$-theory.
The visible symmetric $L$-theory will be reviewed in section~\ref{sec-visrevis}.
There are forgetful natural transformations
\[ \bL_\bullet(Y,\xi)\lra \bVL^{\bullet}(Y,\xi) \]
which fit into a homotopy cartesian square
\[
\CD
\bL_\bullet^\%(Y,\xi)@>>> \bVL^{\bullet\%}(Y,\xi) \\
@VVV   @VVV  \\
\bL_\bullet(Y,\xi)@>>> \bVL^{\bullet}(Y,\xi).
\endCD
\]
This will also be reviewed in section~\ref{sec-visrevis}.
Together with the norm fibration sequence
\[
\CD
S^1\wedge \bA(Y,\xi,m)_{h\ZZ/2} @<<<    \bA(Y,\xi,m)^{th\ZZ/2} @<\textup{incl.}<<  \bA(Y,\xi,m)^{h\ZZ/2}
\endCD
\]
(and a variant with $\bA^\%$ instead of $\bA$), this explains why the two competing definitions of
$\bLA_{\bullet\%}(Y,\xi,m)$, relying on diagrams~(\ref{eqn-defLA}) and~(\ref{eqn-pracdefLA}) respectively,
are consistent. \newline
Our reasons for preferring the second definition of $\bLA_{\bullet\%}$ are strategic.
Quadratic $L$-theory famously serves as a receptacle for relative invariants, such as surgery obstructions
of degree one normal \emph{maps} $X\to Y$. By contrast,
visible symmetric $L$-theory is a mild refinement of symmetric $L$-theory and as such a good receptacle for absolute
invariants: generalized signatures
of Poincar\ee duality spaces, say. In particular a Poincar\ee duality space $Y$ of
formal dimension $m$, and with Spivak normal fibration $\xi$, determines a characteristic element
\[
v_L(Y) \in \Omega^{\infty+m}\bVL^{\bullet}(Y,\xi)
\]
which can be viewed as a refined signature. (We think of it as a point in an
infinite loop space, not a connected component of an infinite loop space.) Refining this some more to pick up
algebraic $K$-theory information, we get
\begin{equation} \label{eqn-hosignature}
\sigma(Y) \in \Omega^{\infty+m}\bVLA^{\bullet}(Y,\xi,m)
\end{equation}
where
\[ \Omega^m\bVLA^{\bullet}(Y,\xi,m) :=
\holim\left(
\CD
 @. \bA(Y,\xi,m)^{h\ZZ/2} \\
 @.   @VV\textup{inclusion} V \\
\Omega^m\bVL^{\bullet}(Y,\xi) @>\Xi>> \bA(Y,\xi,m)^{th\ZZ/2}
\endCD
\right).  \]
This refinement expresses a compatibility between $v_L(Y)$ and the self-dual Euler characteristic
\[ v_K(Y)\in \Omega^\infty(\bA(Y,\xi,m)^{h\ZZ/2})~. \]
We construct $\sigma(Y)$ in section~\ref{sec-hocharsig}.
The construction enjoys continuity properties. It can be applied to the fibers of a fibration $E\to B$ whose fibers
are Poincar\ee duality spaces $E_b$ of formal dimension $m$, so that we obtain a section of a fibration on $B$ whose fibers are
certain infinite loop spaces. \newline
Now suppose that the Poincar\ee duality space $Y$ is a closed manifold of dimension $m$. Then the point $\sigma(Y)$
lifts across the visible $L$-theory and $A$-theory assembly maps to a point
\begin{equation} \label{eqn-excisignature}
\sigma^\%(Y) \in \Omega^{\infty+m}\bVLA^{\bullet\%}(Y,\xi,m)~.
\end{equation}
We construct $\sigma^\%(Y)$ in section~\ref{sec-excicharsig}.
Again this construction enjoys continuity properties: it can be applied it to the fibers of a
fibre bundle $E\to B$ whose fibers are closed manifolds $E_b$ of dimension $m$, so that we obtain a section of a
fibration on $B$ whose fibers are certain infinite loop spaces. (This is very hard to establish, like
the continuity property of excisive Euler characteristics in \cite{DwyerWeissWilliams}.
Relying on \cite{McDuff80}, we reduce to the case of fiber bundles with
discrete structure group.)
\newline
In particular, the space $\sS(M)$ carries a universal bundle $E\to \sS(M)$ of closed manifolds with a fiber
homotopy trivialization $E\simeq \sS(M)\times M$. Therefore each point $(N,f)\in \sS(M)$ determines an
element $f_*\sigma\upr(N)\in \Omega^{\infty+m}\bVLA\ubul\upr(M,\nu_M,m)$, whose image in
$\Omega^{\infty+m}\bVLA\ubul(M,\nu_M,m)$
under assembly comes with a preferred path to $\sigma(M)\in \Omega^{\infty+m}\bVLA\ubul(M,\nu_M,m)$. This gives us the
map~(\ref{eqn-mainmap}): here it comes as a map from $\sS(M)$ to the homotopy fiber, over the point $\sigma(M)$,
of the assembly map
\[ \Omega^{\infty+m}\bVLA\ubul\upr(M,\nu_M,m) \lra \Omega^{\infty+m}\bVLA\ubul(M,\nu_M,m)~. \]
Similar ideas, i.e., a firework of characteristics and signatures,
can be used to show that the map~(\ref{eqn-mainmap}) is highly connected; we give an overview in
section~\ref{sec-outline} before developing the details.

\medskip
This result has many precursors. The most fundamental and best known of these belong
to surgery theory. From the surgery point of view it is very natural
to introduce certain ``block'' structure spaces such as
\[ \widetilde\sS^s(M)~, \quad \widetilde\sS^h(M)\,. \]
These are designed in such a way that $\pi_0\widetilde\sS^s(M)$ and
$\pi_0\widetilde\sS^h(M)$ are identifiable with, respectively, the subset of $\pi_0\sS(M)$
determined by the simple homotopy equivalences, and the quotient set of $\pi_0\sS(M)$ determined
by the $h$-cobordism relation. In addition they have the property
\[ \pi_i\widetilde\sS^s(M) \cong \pi_0\widetilde\sS^s(M\times D^i)~,\quad
\pi_i\widetilde\sS^h(M) \cong \pi_0\widetilde\sS^h(M\times D^i) \,. \]
This is obviously very useful in calculations.
The surgery-theoretic calculations of these spaces are of the form
\begin{equation} \label{eqn-surgcalc}
\begin{split}
\begin{array}{c}
\widetilde\sS^s(M)\,\,\simeq\,\,
\fiber\left[\,\Omega^{\infty+m}\bL^s\lbul\lpr(M,w)\right. \lra \left.8\ZZ\,\right], \rule{0mm}{-3mm}\\
\widetilde\sS^h(M)\,\,\simeq\,\,
\fiber\left[\,\Omega^{\infty+m}\bL^h\lbul\lpr(M,w)\right. \lra \left.8\ZZ\,\right], \rule{0mm}{+4mm}
\end{array}
\end{split}
\end{equation}
where $\bL^s\lbul\lpr$ and $\bL^h\lbul\lpr$ are homotopy invariant functors from
spaces with double coverings to spectra. (In particular $w$ denotes the orientation
covering of $M$.) The functors $\bL^s\lbul\lpr$ and $\bL^h\lbul\lpr$
can be defined entirely in terms of algebraic $L$-theory, again compounded
with assembly. They are therefore fully 4-periodic:
\[ \Omega^4\bL^s\lbul\lpr(X,v)\simeq \bL^s\lbul\lpr(X,v)~, \qquad
\Omega^4\bL^h\lbul\lpr(X,v)\simeq \bL^h\lbul\lpr(X,v)\,. \]
This calculation of $\widetilde\sS^s(M)$ and $\widetilde\sS^h(M)$ is sometimes called the
Casson-Sullivan-Wall-Quinn-Ranicki theorem. An earlier version of it, describing
the homotopy groups of the block structure space(s), is known as the
Casson-Sullivan-Wall long exact sequence. The space level formulation was championed by
Quinn. The complete and final reduction to $L$-theory, at the space level,
is mainly due to the untiring efforts of Ranicki. This took many years.

\medskip
Our calculation of structure spaces $\sS(M)$ in the concordance stable range
is in agreement with the surgery theoretic calculation of block structure
spaces. For example, there is a commutative diagram
\[
\xymatrix@C45pt@M18pt{
\sS^s(M) \ar[r] \ar[d]^-{\textup{incl.}} & \ar[d]^-{\textup{forgetful}}
\fiber[\Omega^{\infty+m}\bLA\lbul\lpr(M,\nu,m)\lra \Wh(\pi_1M)\times 8\ZZ] \\
\widetilde\sS^s(M) \ar[r]^-{\simeq} & \fiber[\Omega^{\infty+m}\bL\lbul^s\lpr(M,w)\lra 8\ZZ].
}
\]
where the upper horizontal arrow is the restriction of the map~(\ref{eqn-mainmap}).
Passing to vertical homotopy fibers, we obtain a highly connected map
\[  \widetilde\TOP(M)/\TOP(M) \lra \Omega^\infty(\bA^s\lpr(M,\nu,m)\ho)~. \]
This is reminiscent of a highly connected map
\[ \widetilde\TOP(M)/\TOP(M) \lra \Omega^{\infty}(\mathbf H^s(M)\ho) \]
constructed in \cite{WW1}; see also \cite{WWsurvey} for notation.
These two maps are intended to be the same, modulo Waldhausen's
identification of the $h$-cobordism spectrum $\mathbf H^s(M)$ with $\bA^s\lpr(M)$.
We do not quite prove that here, but we come close to it. It will be the theme
of another paper in this series.

\medskip
Meanwhile Burghelea and Lashof \cite[cor.\ D]{BurgheleaLashof82} obtained
results on the homotopy type of
$\sS(M)$.
Localizing at odd primes, they were able to construct
a highly connected map
\[ \Omega\sS(M)\,\,\lra \,\,
\Omega\widetilde\sS(M) \times \Omega^{\infty+1}\bA\lpr(M,\nu,m)^{h\ZZ/2}~. \]
(The localization is applied to $\sS(M)$, $\widetilde\sS(M)$ and
$\bA\lpr(M,\nu,m)$ before other operations are carried out: $\Omega$ in both sides, $\Omega^{\infty+1}$
and the homotopy fixed point operation in the right-hand side.)
After localization of $\bA\lpr(M,\nu,m)$ at odd primes, the homotopy  
fixed point spectrum $\bA\lpr(M,\nu,m)^{h\ZZ/2}$ is a wedge summand
of $\bA\lpr(M)$ which depends only on $\nu$ and the parity of $m$.

\medskip
With hindsight, the Burghelea-Lashof result can be explained in terms
of our calculation of $\sS(M)$ described above and the surgery-theoretic calculation of
the block structure space. At odd primes, the six-term diagram~(\ref{eqn-pracdefLA})
simplifies because the Tate constructions in the middle column (again to be applied after the
localization of $\bA\lpr(Y,\xi,m)$) vanish. Therefore at odd primes
\[ \Omega^m\bLA_{\bullet\%}(Y,\xi,m)~\simeq~
\Omega^m\bL_{\bullet\%}(Y,\xi) \vee \bA_\%(Y,\xi,m)^{h\ZZ/2}~.  \]

\bigskip
This paper is a continuation of \cite{WW1} and \cite{WW2}. In another
sense it is a continuation of \cite{DwyerWeissWilliams}.
For technical support, we use a fair amount of controlled topology as
in \cite{ACFP}, the Thurston-Mather-McDuff-Segal discrete approximation theory
\cite{McDuff80} for homeomorphism groups as in \cite{DwyerWeissWilliams},
and Spanier-Whitehead duality theory with its implications for algebraic $K$-theory as in \cite{WWduality}.

\section{Outline of proof} \label{sec-outline}
In the introduction, we gave a rough description of certain invariants of type signature and Euler characteristic
for manifolds and Poincar\ee duality spaces. This led us to a map of the form~(\ref{eqn-mainmap}).
We wish to show that the map is highly connected. The main tools in the proof are
\begin{itemize}
\item[(i)] a controlled version of the Casson-Sullivan-Wall-Quinn-Ranicki (CSWQR) theorem in surgery theory;
\item[(ii)] more invariants of type signature and Euler characteristic
for manifolds and Poincar\ee duality spaces in a controlled setting;
\item[(iii)] a simple downward induction, where the induction beginning relies on (i) while (ii) enables us to
do the induction steps.
\end{itemize}
Let $\sS(M\times\RR^i\cc)$
be the \emph{controlled} structure space of $M\times\RR^i$; here we view $M\times\RR^i$ as an open dense subset of
the join $M*S^{i-1}$. An element of $\sS(M\times\RR^i\cc)$ should be thought of as a pair $(N,f)$ where $N$ is a manifold
of dimension $m+i$, without boundary, and $f\co N\to M\times\RR^i$ is a \emph{controlled} homotopy equivalence \cite{ACFP}.
There is also a controlled block structure space
\[  \widetilde\sS^{cs}(M\times\RR^i\cc) \]
where the decoration $cs$ (controlled simple)
indicates that we allow only structures with vanishing controlled Whitehead torsion.

\smallskip
The homotopy type of $\widetilde\sS^{cs}(M\times\RR^i\cc)$ can be described by a formula which combines the CSWQR ideas
with controlled algebra \cite{ACFP}: namely,
\begin{equation} \label{eqn-CSWQRcontrol}
\widetilde\sS^{cs}(M\times\RR^i\cc)\,\,\simeq\,\,
\fiber\left[\,\Omega^{\infty+m+i}\bL^{cs}_{\bullet}\lpr(M\times\RR^i,\nu\cc) \lra 8\ZZ\,\right]
\end{equation}
where $\bL^{cs}\lbul(M\times\RR^i,\nu\cc)$ is the controlled quadratic $L$-theory (with vanishing controlled
Whitehead torsion) of the control space $(M*S^{i-1},M\times\RR^i)$. Taking $i$ to the limit we have
\[ \colimsub{i\ge 0} \widetilde\sS^{cs}(M\times\RR^i\cc)
\,\,\simeq\,\,
\fiber\left[\,\colimsub{i\ge 0}\Omega^{\infty+m+i}\bL^{cs}_{\bullet}\lpr(M\times\RR^i,\nu\cc) \lra 8\ZZ\,\right]
\]
where the colimits are formed using product with $\RR$ in various shapes. Moreover, it is well-known \cite{WW1} that the inclusions
\[ \colimsub{i\ge 0} \sS(M\times\RR^i\cc) \longleftarrow
\colimsub{i\ge 0} \sS^{cs}(M\times\RR^i\cc) \lra   \colimsub{i\ge 0} \widetilde\sS^{cs}(M\times\RR^i\cc) \]
are homotopy equivalences. Therefore we have
\begin{equation} \label{eqn-industart} \colimsub{i\ge 0} \sS(M\times\RR^i\cc) \,\,\simeq\,\,
\fiber\left[\,\colimsub{i\ge 0}\Omega^{\infty+m+i}\bL^{cs}_{\bullet}\lpr(M\times\RR^i,\nu\cc) \lra 8\ZZ\,\right]
\end{equation}
and this is the starting point for our downward induction. \newline
Next we discuss the induction steps.
Let $(\bar Y,Y)$ be a control space. For the present purposes we can take this to mean
that $\bar Y$ is compact metrizable, and $Y$ is open dense in $\bar Y$. A choice of spherical fibration
$\xi$ on $Y$ and integer $m$ makes the Waldhausen category of locally finitely dominated retractive spaces
over $Y$ into a Waldhausen category with duality (see \cite{DwyerWeissWilliams} for details). By forming
$L$-theory, $K$-theory etc.,
we define spectrum-valued functors
\[ (\bar Y,Y,\xi)\quad\mapsto
\left\{\begin{array}{l} \bL\lbul(Y,\xi\cc), \\
\bVL\ubul(Y,\xi\cc), \\
\bA(Y\cc), \\
\bLA\lbul(Y,\xi,m\cc), \\
\bVLA\ubul(Y,\xi,m\cc)
\end{array} \right.
\]
much as before. (Three of these can be viewed as functors of a general Waldhausen category with duality; the ones
having a $\mathbf V$ in their name use more special features.) The symbol $c$ is a shorthand for control conditions,
allowing us to avoid direct reference to the inclusion $Y\to \bar Y$.
There are natural assembly transformations
\begin{equation} \label{eqn-ctrassemblylist}
\begin{array}{ccc}
\bL\lbul\upr(Y,\xi\cc) \lra \bL\lbul(Y,\xi\cc), \\
\bVL\ubul\upr(Y,\xi\cc) \lra \bVL\ubul(Y,\xi\cc), \\
\bA\upr(Y\cc) \lra \bA(Y\cc),  \\
\bLA\lbul\upr(Y,\xi,m\cc) \lra \bLA\lbul(Y,\xi,m\cc), \\
\bVLA\ubul\upr(Y,\xi,m\cc) \lra \bVLA\ubul(Y,\xi,m\cc),
\end{array}
\end{equation}
where the domain is now designed so that its homotopy groups are the \emph{locally finite}
generalized homology groups of $Y$ with (twisted where appropriate)
coefficients in $\bL\lbul(\pt,\xi)$, $\bVL\ubul\upr(\pt,\xi)$,
$\bA(\pt)$, $\bLA\lbul(\pt,\xi,m)$ and $\bVLA\ubul\upr(\pt,\xi,m)$. Here $\pt$ should be thought
of as a variable point in $Y$, and we restrict $\xi$ from $Y$ to that point where necessary. The homotopy fibers
of the assembly maps~(\ref{eqn-ctrassemblylist}) are denoted by
\begin{equation} \label{eqn-ctrassemblyfib}
\begin{array}{rl}
 & \bL\lbul\lpr(Y,\xi\cc) \\
  \simeq & \bVL\ubul\lpr(Y,\xi\cc), \\
 &  \bA\lpr(Y\cc),  \\
& \bLA\lbul\lpr(Y,\xi,m\cc) \\
\simeq & \bVLA\ubul\lpr(Y,\xi,m\cc),
\end{array}
\end{equation}
respectively. (The homotopy equivalences asserted here are nontrivial; they are established in section~\ref{sec-controlvisible}.)
If $(\bar Y,Y)$ happens to be a \emph{controlled} Poincar\ee duality space of formal dimension $m$ and with
Spivak normal fibration $\xi$, then there is a signature invariant
\begin{equation} \label{eqn-hosignaturectrl}
\sigma(Y) \in \Omega^{\infty+m}\bVLA\ubul(Y,\xi,m\cc)
\end{equation}
which generalizes~(\ref{eqn-hosignature}). This invariant has the expected naturality and continuity
properties. It is constructed in section~\ref{sec-hocharsig}.

If $Y$ happens to be a manifold of dimension $m$ and $\xi=\nu$ is its
normal bundle, then $(\bar Y,Y)$ is automatically a controlled Poincar\ee duality space of formal dimension $m$ and the
signature invariant $\sigma(Y)$ lifts across the assembly map~(\ref{eqn-ctrassemblylist}) to an element
\begin{equation} \label{eqn-excisignaturectrl}
\sigma\upr(Y) \in \Omega^{\infty+m}\bVLA\ubul\upr(Y,\xi,m),
\end{equation}
generalizing~(\ref{eqn-excisignature}). This lift is constructed in section~\ref{sec-excicharsig}.
In particular, the space $\sS(M\times\RR^i\cc)$ carries a universal
bundle where each fiber is an $(m+i)$-manifold $N$ together with a controlled homotopy
equivalence $f\co N\to M\times\RR^i$. We may compactify each fiber $N$ to a control space $\bar N=N\cup S^{i-1}$
in such a way that $N$ is open dense in $\bar N$ and $f$ extends to a map from $\bar N$ to $M*S^{i-1}$.
Therefore each point $(N,f)\in \sS(M\times\RR^i\cc)$ determines an
element
\[ f_*\sigma^\%(N)\in \Omega^{\infty+m+i}\bVLA\ubul\upr(M\times\RR^i,\nu,m+i\cc) \]
whose image in $\Omega^{\infty+m+i}\bVLA\ubul(M\times\RR^i,\nu,m+i\cc)$
under assembly~(\ref{eqn-ctrassemblylist}) comes with a preferred path to $\sigma(M\times\RR^i)$.
If this construction were to enjoy certain continuity properties, it would give us a map
\[
\xymatrix@M10pt{
\sS(M\times\RR^i\cc) \ar@{..>}[r] & \Omega^{\infty+m+i}\bLA\lbul\lpr(M\times\RR^i,\nu,m+i\cc)
}
\]
generalizing~(\ref{eqn-mainmap}),
where we think of the target as the homotopy fiber over the point $\sigma(M\times\RR^i)$ of the
appropriate assembly map in controlled $\bVLA\ubul$ theory of $M\times \RR^i$.
Unfortunately we could not avoid some sacrifices in establishing the continuity properties, and
so we only get a map
\begin{equation} \label{eqn-mainmapctr}
\xymatrix@M10pt{
\sS^\red(M\times\RR^i\cc) \ar[r] & \Omega^{\infty+m+i}\bLA\lbul\lpr(M\times\RR^i,\nu,m+i\cc)
}
\end{equation}
where $\sS^\red(M\times\RR^i\cc)\subset \sS(M\times\RR^i\cc)$ is the union of the connected components
of $\sS(M\times\RR^i\cc)$ which are reducible in the sense that they come from $\pi_0\sS(M)$.
Combining the maps~(\ref{eqn-mainmapctr}) for all $i\ge 0$ results in a commutative ladder
\begin{equation} \label{eqn-mainmapladder}
\begin{split}
\xymatrix@R16pt{
\vdots & \vdots \\
\sS^\red(M\times\RR^{i+1}\cc) \ar[r] \ar[u] & \Omega^{\infty+m+i+1}\bLA\lbul\lpr(M\times\RR^{i+1},\nu,m+i+1\cc) \ar[u] \\
\sS^\red(M\times\RR^i\cc) \ar[r] \ar[u] & \Omega^{\infty+m+i}\bLA\lbul\lpr(M\times\RR^i,\nu,m+i\cc) \ar[u] \\
\vdots \ar[u]  & \vdots \ar[u] \\
\sS^\red(M\times\RR\cc) \ar[r] \ar[u] & \Omega^{\infty+m+1}\bLA\lbul\lpr(M\times\RR^1,\nu,m+1\cc) \ar[u] \\
\sS(M) \ar[r] \ar[u] & \Omega^{\infty+m}\bLA\lbul\lpr(M,\nu,m) \ar[u]
}
\end{split}
\end{equation}
where the vertical arrows are given by product with $\id_\RR$ in the left-hand column, and product with
$\sigma^\%(\RR)$ in the right-hand column. Each vertical arrow in the left-hand column induces
a surjection on $\pi_0$. At the bottom of the ladder we recognize the map~(\ref{eqn-mainmap})
and at the top we recognize with a small effort
(see section~\ref{sec-approxproofs}) the map of~(\ref{eqn-industart}). In particular, all homotopy fibers of the horizontal map at the top of the
ladder are either contractible or empty. We use downward induction to establish a similar property for \emph{all}
horizontal maps in the ladder:
\begin{itemize}
\item[($\dagger$)] for each of these maps, all homotopy fibers are highly connected or empty.
\end{itemize}
It is enough show that in each square
\begin{equation} \label{eqn-mainmapsquare}
\begin{split}
\xymatrix@R16pt{
\sS^\red(M\times\RR^{i+1}\cc) \ar[r] & \Omega^{\infty+m+i+1}\bLA\lbul\lpr(M\times\RR^{i+1},\nu,m+i+1\cc) \\
\sS^\red(M\times\RR^i\cc) \ar[r] \ar[u] & \Omega^{\infty+m+i}\bLA\lbul\lpr(M\times\RR^i,\nu,m+i\cc) \ar[u] \\
}
\end{split}
\end{equation}
of the ladder, all total homotopy fibers are highly connected or empty. \newline
Each \emph{vertical}
homotopy fiber in the left-hand column can be identified with a union of connected components of
a controlled $h$-cobordism space $\sH(N\times\RR^i\cc)$, where $N$ is some closed $m$-manifold homotopy
equivalent to $M$. By an easy calculation carried out mainly in section~\ref{sec-controlsusp}, the vertical homotopy fibers
in the right-hand column have the form
\[ \Omega^\infty\bA\lpr(M\times\RR^i\cc)~\simeq~\Omega^\infty\bA\lpr(N\times\RR^i\cc)~. \]
With these descriptions, the map between matching vertical homotopy fibers in~(\ref{eqn-mainmapsquare})
extends to a controlled form
\begin{equation} \label{eqn-Whsonfibers}
\sH(N\times\RR^i\cc) \lra \Omega^\infty\bA\lpr(N\times\RR^i\cc)
\end{equation}
of Waldhausen's map relating $h$-cobordism spaces to $A$-theory. This is verified in section~\ref{sec-approxproofs}.
The map~(\ref{eqn-Whsonfibers}) is highly connected.
So all its homotopy fibers are highly connected, and so our claim regarding~(\ref{eqn-mainmapsquare}) is proved,
and claim $(\dagger)$ is also established. In particular,
any homotopy fiber of our map
\[ \sS(M) \lra \Omega^{\infty+m}\bLA\lbul\lpr(M,\nu,m) \]
is highly connected of empty. It only remains to show that the nonempty homotopy fibers correspond
to elements of $\Omega^{\infty+m}\bLA\lbul\lpr(M,\nu,m)$ whose connected component
is in the kernel of the local degree homomorphism to $8\ZZ$.

\medskip
For this we use the commutative diagram
\[
\xymatrix{
\pi_0\sS(M) \ar[r] \ar[d]_-{\textup{induced by incl.}} & \pi_m\bLA\lbul\lpr(M,\nu,m) \ar[r] \ar[d]_-{\textup{forget}}
& 8\ZZ \ar[d]_-{=} \\
\pi_0\widetilde\sS^h(M) \ar[r] & \pi_m\bL\lbul\lpr(M,\nu) \ar[r] & 8\ZZ
}
\]
where the lower row is short exact.
The left-hand vertical arrow is onto by definition. Its fibers are
the orbits of an action of $\Wh(\pi_1M)$ on $\pi_0\sS(M)$.
By direct calculation, and almost by construction, the middle vertical arrow (which is a group homomorphism)
is also onto and its kernel is the image of a homomorphism
\begin{equation} \label{eqn-homact}
\pi_0((\bA^h\lpr(M,\nu,m))_{h\ZZ/2}) \lra \pi_m\bLA\lbul\lpr(M,\nu,m).
\end{equation}
Here $\pi_0((\bA^h\lpr(M,\nu,m))_{h\ZZ/2})$ is a quotient of $\Wh(\pi_1M)$. Hence we need to show that
the action of the Whitehead group in the upper left-hand term corresponds in the
upper middle term to a translation action, using the homomorphism~(\ref{eqn-homact}). This is
done in section~\ref{sec-approxproofs}.

\section{Visible $L$-theory revisited} \label{sec-visrevis}
\label{sec-visible} Mishchenko \cite{MishchenkoSig} and Ranicki
\cite{RanickiLMS1}, \cite{RanickiLMS2} introduced ``symmetric
structures'' on certain chain complexes over rings with involution
with a view to understanding signature invariants and product
formulae in surgery theory. For a ring $R$ with involution
(=involutory antiautomorphism) $r\mapsto \bar r$ and a bounded chain
complex $C$ of finitely generated projective left $R$--modules, a
symmetric structure of dimension $m$ on $C$ is a chain map of
$\ZZ[\ZZ/2]$--module chain complexes
\[ \varphi\co \Sigma^mW\lra C^t\otimes_R C\,. \]
Here $C^t$ is $C$ with the right $R$--module structure defined by
$xr=\bar rx$ and $W$ denotes the standard resolution of $\ZZ$ by
free $\ZZ[\ZZ/2]$--modules:
\[
\CD \ZZ[\ZZ/2] @<1-T<< \ZZ[\ZZ/2] @<1+T<\phantom{}< \ZZ[\ZZ/2]
@<1-T<<\ZZ[\ZZ/2] @<1+T<<\cdots\,\,\,. \endCD
\]
The value of $\varphi$ on $1\in W_0$ is an $m$--cycle in
$C^t\otimes_RC$, corresponding to a degree $m$ chain map from the
dual $C^{-*}$ to $C$. If this is a chain homotopy equivalence,
$\varphi$ is called \emph{nondegenerate}. The bordism groups of
objects $(C,\varphi)$ as above, with nondegenerate $\varphi$ of
dimension $m$, are the symmetric $L$--groups $L^m(R)$. (This
definition of $L^m(R)$ is in agreement with \cite{RanickiTopMan} but
in slight disagreement with \cite{RanickiLMS1} because we do not
require that $C_k=0$ for $k\notin \{0,1,\dots,n\}$.) Ranicki's
analogous description of the quadratic $L$--groups $L_n(R)$, in
which the homotopy fixed point construction
$\hom_{\ZZ[\ZZ/2]}(W,\blank)$ is replaced by a homotopy orbit
construction $W\otimes_{\ZZ[\ZZ/2]}\,\blank$, makes it easy to
define multiplication operators
\[  L^m(R_1) \otimes L_n(R_2) \lra L_{m+n}(R_1\otimes R_2). \]
A (connected) Poincar\ee duality space $X$ of formal dimension $m$
with fundamental group $\pi$ and orientation character $w\co \pi\to
\set{\pm1}$ determines
\begin{itemize}
\item an involution on $R=\ZZ\pi$ given by $g\mapsto w(g)\cdot g^{-1}$
for $g\in \pi\subset\ZZ\pi$,
\item a chain complex $C$ over $R=\ZZ\pi$, the singular or cellular chain complex
of the universal cover of $X$, and
\item a nondegenerate $n$-dimensional structure
$\varphi$ on $C$, obtained by evaluating the Eilenberg--Zilber
diagonal chain map on the fundamental class of $X$.
\end{itemize}
The corresponding element in $L^m(R)$ is the \emph{symmetric
signature} $\sigma^*(X)$ of $X$. If $X$ is a closed manifold and
$f\co X'\to Y$ is a degree one normal map from a closed
$n$--manifold to a Poincar\ee duality space of formal dimension $n$,
then
\[ (\id_X\times f)\co X\times X'\to X\times Y\]
is also a degree one normal map. The surgery obstructions
$\sigma_*(f)$ and $\sigma_*(\id_X\times f)$ are related by
\[ \sigma_*(\id_X\times f) = \sigma^*(X)\cdot \sigma_*(f)\in L_{m+n}(\ZZ\pi\otimes\ZZ\pi'), \]
using the above product, with $\pi'=\pi_1(Y)$.

\medskip
A few years later it was found \cite{WeissLMS} that the symmetric
$L$--groups admit a homological description relative to the
quadratic $L$--groups. That is, there is a long exact sequence
\[
\cdots \to L_n(R)\to L^n(R)\to \widehat L^n(R)\to L_{n-1}(R)\to
\cdots
\]
and the calculation of the relative terms $\widehat L^n(R)$ is
``only'' a matter of homological algebra. Efforts to reduce the
homological algebra to an absolute minimum eventually led to the
visible symmetric $L$--groups $VL^n(R)$, defined for group rings
$R=\ZZ\pi$. We now recall their definition, following
\cite{WeissVisible}.
\newline
Let $R=\ZZ\pi$ with involution $g\mapsto w(g)\cdot g^{-1}$ for some
homomorphism $w\co\pi\to \set{\pm1}$. Let $C$ be a chain complex of
f.g. projective left $R$--modules, bounded as before. A symmetric
structure of dimension $n$ on $C$ can be viewed as an $n$--cycle in
\[ ((C^t\otimes C)^{h\ZZ/2})_{\pi} \cong
((C^t\otimes C)_{\pi})^{h\ZZ/2}  \] where the various subscripts and
superscripts indicate orbit constructions and homotopy fixed point
constructions for the appropriate symmetry groups, here $\pi$ and
$\ZZ/2$. (Note that $\pi$ acts diagonally on $C^t\otimes C$.) A
\emph{visible symmetric structure} of dimension $m$ on $C$ is an
$m$--cycle in
\[ ((C^t\otimes C)^{h\ZZ/2})_{h\pi}\]
where $(\blank)_{h\pi}$ means $(\blank \otimes P)_{\pi}$ for a
resolution $P$ of the trivial module $\ZZ$ by projective
$\ZZ\pi$--modules. In contrast to $C^t\otimes C$ and $(C^t\otimes
C)_{h\ZZ/2}$, the chain complex $(C^t\otimes C)^{h\ZZ/2}$ is not
bounded below if $C\ne 0$, so that there is no good reason to think
that the augmentation--induced chain map
\[ ((C^t\otimes C)^{h\ZZ/2})_{h\pi} \lra ((C^t\otimes C)^{h\ZZ/2})_{\pi} \]
should induce an isomorphism in homology. In fact visible symmetric
structures generally carry more information than symmetric
structures. It is sometimes convenient to organise both types of
structures into homotopy classes: the groups of such homotopy
classes are denoted
\[ Q^m(C)= H_m((C^t\otimes C)^{h\ZZ/2})_{\pi})\,,
\qquad VQ^m(C) = H_m((C^t\otimes C)^{h\ZZ/2})_{h\pi}).
\]
With a view to generalizations later on, we mention that there is a
homotopy (co)cartesian square of chain complexes
\[
\CD
((C^t\otimes C)^{h\ZZ/2})_{h\pi} @>>> ((C^t\otimes C)^{th\ZZ/2})_{h\pi} \\
@VVV    @VVV  \\
((C^t\otimes C)^{h\ZZ/2})_{\pi} @>>> ((C^t\otimes C)^{th\ZZ/2})_{\pi} \\
\endCD
\]
where the ${}^{th\ZZ/2}$ superscript denotes a Tate construction
$\hom_{\ZZ[\ZZ/2]}(\widehat W,\blank)$, with
\[
\begin{array}{rcl}
\widehat W_k & = &\ZZ[\ZZ/2] \quad \textup{for all } k\in \ZZ, \\
d_k\co \widehat W_k\to \widehat W_{k-1} & ; & z\mapsto
(1+(-1)^kT)z\,.
\end{array}
\]
This is reflected in a long exact Mayer--Vietoris sequence
\[ \cdots\to \widehat Q^{n+1}(C)\to
VQ^n(C) \to Q^n(C)\oplus V\widehat Q^n(C) \to \widehat Q^n(C) \to
VQ^{n-1}(C) \to\cdots \] where $V\widehat Q^n(C)=H_n(((C^t\otimes
C)^{th\ZZ/2})_{h\pi})$ and $\widehat Q^n(C) = H_n(((C^t\otimes
C)^{th\ZZ/2})_{\pi})$.  \newline The visible hyperquadratic theory
has the property of being invariant under a "change of rings". That
is, for a bounded chain complex $C$ of f.g. left projective
$\ZZ\pi$--modules and a homomorphism $h\co\pi\to \pi'$ we have
\[ V\widehat Q^*(C) \cong V\widehat Q^*(C') \]
where $C'=\ZZ\pi'\otimes_{\ZZ\pi}C$. (It is understood that $\pi$
and $\pi'$ are equipped with homomorphisms to $\set{\pm1}$ and that
$h$ respects these.) We also note that
\[ V\widehat Q^* = \widehat Q^* \]
in the case $\pi=\set{1}$. These two properties, suitably sharpened,
could be used to characterize the visible hyperquadratic theory in
terms of the ordinary hyperquadratic one. But we need not go into
that. \newline A visible symmetric structure on $C$ is considered
nondegenerate if the induced symmetric structure is nondegenerate.
The bordism groups of chain complexes $C$ as above with a
nondegenerate $m$--dimensional visible symmetric structure are the
visible symmetric $L$--groups $VL^m(\ZZ\pi)$. A mild improvement on
the Mishchenko--Ranicki construction of the symmetric signature of a
Poincar\ee duality space $X$ of formal dimension $m$ gives the
\emph{visible symmetric signature}
\[  \sigma^*(X) \in VL^m(\ZZ\pi). \]
Other useful features of the symmetric $L$--groups (such as the
products and the product formula for surgery obstructions) can be
transferred to the visible symmetric $L$--groups by means of the
forgetful homomorphisms $VL^m(\ZZ\pi)\to L^m(\ZZ\pi)$.

\medskip
The main result of \cite{WeissVisible} is a long exact sequence
relating the quadratic $L$--groups of $\ZZ\pi$ to the visible
symmetric $L$--groups, with ``easy'' relative terms:
\[
\cdots \to L_n(\ZZ\pi)\to VL^n(\ZZ\pi)\to \bigoplus_{i+j=n}
H_i(B\pi;\widehat L^j(\ZZ)) \to L_{n-1}(\ZZ\pi) \to \cdots\,\,.
\]
Ranicki \cite{RanickiTopMan} found a generalization of this from the
group ring case to the case of a simplicial group ring $\ZZ[\Omega
X]$, and used it in a revised approach to his total surgery
obstruction theory, a project going back to \cite{RanickiTot}. He
defines a visible symmetric $L$--theory spectrum which we (not he)
denote by $\bV\bL\ubul(\ZZ[\Omega X])$, with homotopy groups
$VL^n(\ZZ[\Omega X])$. There is a long exact sequence
\[
\cdots \to L_n(\ZZ\pi)\to VL^n(\ZZ[\Omega X])\to \bigoplus_{i+j=n}
H_i(X;\widehat L^j(\ZZ)) \to L_{n-1}(\ZZ\pi) \to \cdots\,,
\]
identical with the one above when $X=B\pi$. One of Ranicki's main
results in \cite{RanickiTopMan} states roughly that the closed
manifold structures on an oriented Poincar\ee duality space $X$ of
formal dimension $n$ are in canonical bijection with the connected
components of the homotopy fiber of the assembly map
\[ \Omega^{\infty+m}(X_+\wedge\bV\bL\ubul(\ZZ))
\lra \Omega^{\infty+m}\bV\bL\ubul(\ZZ[\Omega X]) \] over the point
determined by $\sigma^*(X)$, where $\bV\bL\ubul(\dots)$ denotes
visible $L$--theory spectra. (To be more precise, each of these
connected components determines a class in
$\pi_m(\Omega^{\infty}(X_+\wedge\bV\bL\ubul(\ZZ)))$ which has a
\emph{local signature} $d\in 8\ZZ+1$; the condition $d=1$ must be
added.) This is obviously relevant to our program.

\bigskip
We come to a definition of visible symmetric structures in the
setting of retractive spaces and retractive spectra. Let $Y_1$ and
$Y_2$ be finitely dominated retractive spaces over $X$, with
retractions $r_1$ and $r_2$. We recall first the definition of an
`unstable'' Spanier--Whitehead product $Y_1\curlywedge Y_2$ from
\cite[1.A.3]{WWduality}. This is the based space obtained by first
forming the external smash product
\[
\begin{array}{ccc}
Y_1 \,\,{}_X\!\wedge_X\, Y_2 & = & Y_1\times Y_2/\sim
\end{array}
\]
where $\sim$ identifies $(y_1,x)$ with $(r_1(y_1),x)$ and $(x,y_2)$
with $(x,r_2(y_2))$~; then taking the homotopy pullback of
\[
\CD   X @>\textup{diagonal}>> X\times X@<\textup{retraction}<< Y_1
\,\,{}_X\!\wedge_X\, Y_2
\endCD
\]
and then dividing that by the homotopy pullback of
\[
\CD   X @>\textup{diagonal}>> X\times X@<\id<< X\times X \,.
\endCD
\]
We make this unstable $SW$ product ``stable'' essentially by
applying $\Omega^{\infty}\Sigma^{\infty}$. More technically,
however, we have to work in a stable category of retractive spaces
with objects of the form $(Y,k)$ where $Y$ is a (finitely dominated)
retractive space over $X$ and $k\in \ZZ$. The set of morphisms from
$(Y,k)$ to $(Y',\ell)$ in the stable category is
\[  \colim_i\, \mor_{\textup{uns}}(\Sigma^{i-k}Y,\Sigma^{i-\ell}Y') \]
where $\Sigma$ is short for $\Sigma_X$ and $\mor_{\textup{uns}}$
refers to morphisms in the ordinary category of retractive spaces
over $X$. It is worth noting that $(\Sigma^kY,k)$ is isomorphic to
$(Y,0)$ in the stable category, so that $(Y,k)$ can be regarded as a
formal $k$--fold desuspension of $Y$ alias $(Y,0)$.

\begin{defn}
\label{defn-SW} {\rm We let
\[
(Y_1,k)\odot(Y_2,\ell) = \colim_n
\,\Omega^{2n}(\Sigma^{n-k}Y_1\curlywedge\Sigma^{n-\ell}Y_2)\,.
\]
More generally, we let
\[
(Y_1,k)\odot_j(Y_2,\ell) = \colim_n
\,\Omega^{2n}\Sigma^j(\Sigma^{n-k}Y_1\curlywedge\Sigma^{n-\ell}Y_2),
\]
so that $(Y_1,k)\odot(Y_2,\ell)=(Y_1,k)\odot_0(Y_2,\ell)$, and
denote the $\Omega$--spectrum with $j$--th term
$(Y_1,k)\odot_j(Y_2,\ell)$ by
\[ (Y_1,k)\odot\lbul(Y_2,\ell)\,. \]
}
\end{defn}

(By convention $\Omega^{m}Z$, for a based space $Z$ and an integer
$m\ge 0$, is the geometric realization of the simplicial set whose
$n$--simplices are the based maps from the one--point
compactification of $\Delta^n\times\RR^m$ to $Z$. Hence all the
spaces and spectra in definition~\ref{defn-SW} are $CW$ spaces and
$CW$ spectra.) \newline Note that $\odot\lbul$ comes with a
structural symmetry $(Y_1,k)\odot\lbul(Y_2,\ell) \cong
(Y_2,\ell)\odot\lbul(Y_1,k)$ determined by the obvious symmetry of
$\curlywedge$. For $Y_1=Y_2=Y$ and $k=\ell$ we obtain an
$\Omega$--spectrum $(Y,k)\odot\lbul(Y,k)$ with an action of $\ZZ/2$.

\begin{defn}
\label{defn-nonlinearsymandvis} {\rm An \emph{$n$--dimensional
symmetric structure} on $(Y,k)$ is an element of
$\Omega^n((Y,k)\odot(Y,k))^{h\ZZ/2}$. An \emph{$n$--dimensional
visible symmetric structure} on $(Y,k)$ is an element of
$\Omega^n((Y,k)\odot(Y,k))^{\ZZ/2}$.}
\end{defn}

The first part of this definition comes from \cite{WWduality}, but
the second part is new. (We are extremely grateful to John Klein for
suggesting it as an improvement on some earlier attempts of ours.)
It is best understood from the point of view of equivariant homotopy
theory. The $\Omega$--spectrum $(Y,k)\odot\lbul(Y,k)$, with the
action of $\ZZ/2$, turns out to be the ``underlying spectrum'' of a
$\ZZ/2$--spectrum in the sense of the equivariant theory
\cite{Carlsson92}, \cite{May96}. (See also \cite{Adams84}.) We shall
explain this using the following (conservative) language.

\begin{conv} {\rm Let $G$ be a finite group, $W$ the
regular representation of $G$. Let $nW$ be the direct sum of $n$
copies of $W$, with one--point compactification $S^{nW}$. A
\emph{$G$--spectrum} $\bC$ is a family of well--based $G$--spaces
$C_{nW}$, defined for all sufficiently large positive integers $n$,
together with based $G$--maps
\[ S^W\wedge C_{nW}\to C_{(n+1)W}, \]
with the diagonal action of $G$ on $S^W\wedge C_{nW}$. The
\emph{underlying spectrum} $u\bC$ of $\bC$ is the ordinary spectrum
whose $j$--th space is
\[  \colim_n\, \Omega^{nW}\Sigma^jC_{nW} . \]
It is an ordinary $\Omega$--spectrum with a degreewise action of $G$
(coming from the conjugation action of $G$ on
$\Omega^{nW}\Sigma^jC_{nW}$, for each $n$ and $j$). The subspectrum
of fixed points, $(u\bC)^G$, is often called the \emph{fixed point
spectrum of} $\bC$. It is again an $\Omega$--spectrum.}
\end{conv}

\nin\emph{Remark.} The above definition of a $G$--spectrum is
economical in that we only use the representations $nW$ for
bookkeeping. The price for that is a mildly under-motivated
definition of the ``underlying spectrum''. As before, the loop
spaces $\Omega^{nW}$ which appear in the definition of the
underlying spectrum are to be constructed as geometric realizations
of certain simplicial sets, so that the passage to the (co)limit is
safe from the point of view of homotopy theory. Note that if $G$ is
trivial, $G=\{1\}$, then $u\bC$ is simply a $CW$-substitute for
$\bC$. \newline Beware that the expressions \emph{spectrum} and
\emph{$\Omega$--spectrum}, as used here, correspond roughly to
\emph{prespectrum} and \emph{spectrum}, respectively, in the
language of \cite{Carlsson92} and \cite{May96} for example.

\medskip
Returning to definition~\ref{defn-nonlinearsymandvis} now, we have
that $(Y,k)\odot\lbul(Y,k)$ is the underlying spectrum of the
$\ZZ/2$--spectrum given by $nW \, \mapsto \, S^{(n-k)W}\wedge
Y^{\curlywedge 2}$ where $Y^{\curlywedge 2}$ is short for
$Y\curlywedge Y$. (For $G=\ZZ/2$ we like to identify the regular
representation $W$ with the permutation representation on $\RR^2$.)
Clearly this $\ZZ/2$--spectrum (not its underlying spectrum) can be
described as
\[
\bS_{\ZZ/2}^{-kW}\wedge Y^{\curlywedge 2}
\]
where $\bS_{\ZZ/2}^{-kW}$ is a shifted $\ZZ/2$--sphere spectrum,
given by $nW \,\, \mapsto \,\, S^{(n-k)W}$.

\begin{prop}
\label{prop-equifact} For any finite group $G$, any well--based
$G$--space $Z$ which is free away from the base point and any
$G$--spectrum $\bC$, the fixed point spectrum $(u(\bC\wedge Z))^G$
is homotopy equivalent to the homotopy orbit spectrum $(u(\bC\wedge
Z))_{hG}$. Under this identification, the inclusion of $(u(\bC\wedge
Z))^G$ in the homotopy fixed point spectrum $(u(\bC\wedge Z))^{hG}$
corresponds to the norm map.
\end{prop}

\medskip
\begin{cor}
\label{cor-tDalmostsplitting} There is a natural homotopy fiber
sequence of spectra
\[
\CD ((Y,k)\odot\lbul(Y,k))_{h\ZZ/2} @>>>
((Y,k)\odot\lbul(Y,k))^{\ZZ/2} @>J>> \Sigma^{\infty-k}(Y/X)
\endCD \]
where $\Sigma^{\infty-k}(Y/X)$ means a $CW$--substitute for the
suspension spectrum of $(Y/X)$.
\end{cor}

\proof[Proof of the corollary.] Let $T= Y^{\curlywedge 2}$ and
$T'=(E\ZZ/2)_+\wedge T$. Let $f\co T'\to T$ be the projection. The
homotopy cofiber sequence of $\ZZ/2$--spaces $T' \lra T \lra
\cone(f)$ determines a homotopy fiber sequence of spectra
\[ (u(\bC\wedge T'))^{\ZZ/2}
\lra (u(\bC\wedge T))^{\ZZ/2} \lra (u(\bC\wedge \cone(f)))^{\ZZ/2}
\] with $\bC=\bS^{-kW}_{\ZZ/2}$. The middle term in the sequence is
\[ ((Y,k)\odot\lbul(Y,k))^{\ZZ/2} \]
by construction. Proposition~\ref{prop-equifact} identifies the
left--hand term with
\[ (u(\bC\wedge T'))_{h\ZZ/2}
\,\, \simeq \,\,  (u(\bC\wedge T))_{h\ZZ/2}\,=\,
((Y,k)\odot\lbul(Y,k))_{h\ZZ/2}\,. \] The expression $(u(\bC\wedge
\cone(f)))^{\ZZ/2}$ can be identified with $\Sigma^{\infty-k}(Y/X)$
as follows. It is an $\Omega$--spectrum whose $j$--th space is
\[ \colim_n\,\,\map^{\ZZ/2}_*(S^{nW},S^{j\RR\oplus(n-k)W}\wedge\cone(f)) \]
where $\map^{\ZZ/2}_*(...)$ denotes a space of equivariant based
maps and $j\RR$ denotes a $j$--dimensional trivial representation of
$\ZZ/2$. Because $\cone(f)$ is non--equivariantly contractible,
equivariant based maps from $S^{nW}$ to
$S^{j\RR\oplus(n-k)W}\wedge\cone(f)$ are essentially determined by
their restrictions to the fixed point sets. Hence the above
expression for the $j$--th space of $(u(\bC\wedge
\cone(f)))^{\ZZ/2}$ simplifies to
\[ \colim_n \,\, \map_*(S^n,S^{j+n-k}\wedge\cone(f)^{\ZZ/2}).\]
Since $\cone(f)^{\ZZ/2}$ is $T^{\ZZ/2}\simeq Y/X$, this simplifies
even more to $\Omega^{\infty}\Sigma^{\infty+j-k}(Y/X)$, which is the
$j$--th space in the $\Omega$--spectrification of
$\Sigma^{\infty-k}(Y/X)$. \qed

\medskip
\proof[Proof of the proposition.] This is a standard fact from
equivariant homotopy theory. We begin with a preliminary remark
about pathologies. Choose a based $G$--$CW$--space $Z'$ which is
$G$--free away from the base point and a $G$--map $e\co Z'\to Z$
which is a weak equivalence. Because $Z$ is well--based, $v$ induces
weak equivalences
\[ \id\wedge e\co C\wedge Z'\to C\wedge Z \]
for any well--based $G$--space $C$, in particular for $C=C_{nW}$. It
follows that $e$ induces a homotopy equivalence of underlying
spectra,
\[ u(\bC\wedge Z')\to u(\bC\wedge Z), \]
and a homotopy equivalence of the fixed point spectra, $(u(\bC\wedge
Z'))^G\to (u(\bC\wedge Z))^G$. Therefore, without loss of
generality, $Z$ is a based $G$--$CW$--space which is $G$--free away
from the base point. \newline We shall prove the two parts of the
proposition together using a characterization of the norm map as an
``assembly'' transformation. (For this idea we are again indebted to
John Klein.) Let $\sC$ be the category of all based
$G$--$CW$--spaces which are $G$--free away from the base point, with
based $G$--maps as morphisms. Let $\bF$ be a functor from $\sC$ to
spectra which takes homotopy equivalences to weak equivalences and
takes $\pt$ to a weakly contractible spectrum. Then there exists a
natural transformation $\alpha\co \bF^{\%}\to \bF$ where $\bF^{\%}$
is another functor from $\sC$ to spectra, and
\begin{itemize}
\item $\bF^{\%}$ respects weak equivalences,
\item $\bF^{\%}$ respects (weak) homotopy pushout squares
\item $\bF^{\%}$ respects arbitrary wedges up to weak equivalence,
\item $\alpha\co \bF^{\%}(Z)\to \bF(Z)$ is a weak equivalence when $Z=G_+$.
\end{itemize}
The pair $(\bF^{\%}, \alpha)$ is essentially determined by $\bF$ and
$\alpha$ is called the assembly transformation for $\bF$. For the
case $G=\set{1}$, the proof (and a more detailed statement) can be
found in \cite{WWassembly} and the general case follows the same
lines. (One possible definition of $\bF^{\%}(Z)$ for arbitrary $Z$
in $\sC$ is as follows. Take the geometric realization of the
simplicial spectrum
\[ n\mapsto Z(n)\wedge \bF(\Delta^n_+\wedge G_+) \]
where $Z(n)$ is the based set of singular $n$--simplices of $Z$~;
then divide out by the diagonal action of $G$.) \newline Now put
$\bF(Z):=(u(\bC\wedge Z))^{hG}$. This $\bF$ clearly takes homotopy
equivalences to weak equivalences and takes the trivial space $\pt$
to a trivial spectrum. It also respects (weak) homotopy pushout
squares, but it does not satisfy the wedge axiom for infinite
wedges. The norm transformation
\[ (u(\bC\wedge Z))_{hG} \to \bF(Z) \]
satisfies all the properties which characterize the assembly for
$\bF$. Therefore it \emph{is} the assembly for $\bF$. We can now
give our proof by showing that the natural transformation
\[  (u(\bC\wedge Z))^G \lra \bF(Z) \]
given by the inclusion of fixed point spectra in homotopy fixed
point spectra also satisfies all the properties which characterize
the assembly for the functor $\bF$. \newline Of the four properties
listed, three hold by inspection. So it only remains to check that
the inclusion of $(u(\bC\wedge Z))^G$ in $(u(\bC\wedge Z))_{hG}$ is
a homotopy equivalence when $Z=G_+$. From the definitions, an
element of $\pi_n((u(\bC\wedge G_+))^G)$ is represented by a $G$-map
\[   f\co S^{(n+j)\RR\oplus iW} \lra S^{j\RR}\wedge C_{iW} \wedge G_+ \]
with ``large'' $i$ and $j$ (where $j\RR$ for example denotes a
trivial $j$--dimensional representation). We may assume that $f$ is
transverse to $0\times C_{iW}\times G$. The inverse image of
$0\times C_{iW}\times 1$ is then a framed smooth closed
$(n+i|G|)$-dimensional submanifold $M$ of $(n+j)\RR\oplus iW$.
Clearly $M\cap gM=\emptyset$ for $g\in G\smin\set{1}$ and we have a
map $f|M$ from $M$ to $C_{iW}$. Conversely, given any framed smooth
closed $(n+i|G|)$--dimensional submanifold $M$ of $(n+j)\RR\oplus
iW$ with $M\cap gM=\emptyset$ for $g\in G\smin\set{1}$, and a map
$q$ from $M$ to $C_{iW}$, the Pontryagin--Thom construction gives us
an appropriate $f$ for which $M=f^{-1}(0\times C_{iW}\times 1)$ and
$f|M=q$. In the limit, when $j$ and $i$ tend to infinity, the
condition $M\cap gM=\emptyset$ for $G\smin\set{1}$ becomes
irrelevant and so the $n$--th homotopy group under consideration is
identified with the $n$--dimensional framed bordism group of $u\bC$,
i.e., with $\pi_n(u\bC)$. Moreover, this identification clearly
agrees with the homomomorphism induced by the composition
\[ (u(\bC\wedge G_+))^G\lra (u(\bC\wedge G_+))^{hG}\lra
u(\bC\wedge G_+) \lra u(\bC) \] where the second arrow is forgetful
and the third is induced by the (non--equivariant) map $\bC\wedge
G_+\to \bC$ which isolates the summand $\bC\wedge\set{1}_+$. Since
the composition of the last two arrows is a homotopy equivalence,
the first arrow is a homotopy equivalence. \qed

\bigskip
The following continuation of
definition~\ref{defn-nonlinearsymandvis} is suggested by
corollary~\ref{cor-tDalmostsplitting}.

\begin{defn}
\label{defn-nonlinearvishyper} {\rm An $n$--dimensional visible
hyperquadratic structure on $(Y,k)$ is an element in
$\Omega^n\Omega^{\infty}\Sigma^{\infty-k}(Y/X)$. An $n$--dimensional
quadratic structure on $(Y,k)$ is an element of
$\Omega^n\Omega^{\infty}(((Y,k)\odot\lbul(Y,k))_{h\ZZ/2})$.
Alternatively, an $n$--dimensional quadratic structure on $(Y,k)$
can be defined as an element of $\Omega^n\Omega^{\infty}$ of the
homotopy fiber of the natural map $J\co
((Y,k)\odot\lbul(Y,k))^{\ZZ/2} \to \Sigma^{\infty-k}(Y/X)$. }
\end{defn}

\medskip
An $n$--dimensional visible symmetric structure on $(Y,k)$ is
considered nondegenerate if the underlying $n$--dimensional
symmetric structure is nondegenerate. Writing $s\sR(X)$ for the
stable category of finitely dominated retractive spaces over $X$, we
obtain the definition of a visible symmetric $L$--theory spectrum
\[ \bV\bL\ubul(s\sR(X))= \bV\bL\ubul(X) \]
by substituting nondegenerate visible symmetric structures for
nondegenerate symmetric structures throughout in the construction of
the symmetric $L$--theory spectrum $\bL\ubul(s\sR(X))= \bL\ubul(X)$.
See \cite{WWduality}. The standard map
\[  \bL\lbul(X) \to \bL\ubul(X) \]
can be factorized as $\bL\lbul(X)\to \bV\bL\ubul(X) \to
\bL\ubul(X)$. This is clear from
definition~\ref{defn-nonlinearvishyper}.

\medskip
We write  $\bV\widehat\bL\ubul(X)$ for the mapping cone (in the
category of spectra) of the above map $\bL\lbul(X)\to
\bV\bL\ubul(X)$.

\begin{thm}
\label{thm-VhatLexcision} The functor $X\mapsto
\bV\widehat\bL\ubul(X)$ is homotopy invariant and excisive.
\end{thm}

\medskip
``Homotopy invariance'' is intended to mean that the functor takes
weak equivalence to homotopy equivalences, and this is clear. (An
equivalent formulation says that, for each $X$, the maps
$\bV\widehat\bL\ubul(X)\to \bV\widehat\bL\ubul(X\times[0,1])$
induced by $x\mapsto (x,0)$ and $x\mapsto(x,1)$ are homotopic.
They are indeed homotopic because the exact functors which
they induce are related by a chain of natural weak equivalences.)
The excision property means that the functor takes empty space to
a contractible spectrum and takes weak homotopy pushout
squares (also known as cocartesian squares) of spaces to homotopy
pushout squares (equivalently, homotopy pullback squares) of
spectra. Our proof of the excision property relies on three
decomposition lemmas.

\medskip
For the first of these, suppose that $X$ is the union of two closed
subspaces $X_a$ and $X_b$ with intersection $X_{ab}$, such that the
inclusions $X_{ab}\to X_a$ and $X_{ab}\to X_b$ are cofibrations. Let
$r\co E\to X$ be a fibration with section $s$ making $E$ into a
\emph{homotopy finite} retractive space over $X$. Let $Y$ be a
finite retractive space over $X$ with a morphism $f\co Y\to E$ of
retractive spaces over $X$. We assume that $Y$ is decomposed as
\[ Y:= Y_a\cup Y_b \]
where $Y_a$, $Y_b$ and $Y_{ab}=Y_a\cap Y_b$ are finite retractive
spaces over $X_a$, $X_b$ and $X_{ab}$, respectively, with
cofibrations $Y_{ab}\cup X_a\to Y_a$ and $Y_{ab}\cup X_b\to Y_b$.

\medskip
\begin{lem}
\label{lem-retractivedeco} The morphism $f\co Y\to E$ has a
factorization of the form
\[
\CD Y@>f_1>>Z@>g>> E
\endCD
\]
where
\begin{itemize}
\item[(i)] $Z$ is a finite retractive space over $X$
\item[(ii)] $f_1$ is a cofibration
\item[(iii)] $g$ is a weak equivalence
\item[(iv)] the decomposition of $Y$ extends to a similar decomposition
\[ Z:= Z_a\cup Z_b \]
where $Z_a$, $Z_b$ and $Z_{ab}=Z_a\cap Z_b$ are finite retractive
spaces over $X_a$, $X_b$ and $X_{ab}$, respectively, with
cofibrations $Z_{ab}\cup X_a\to Z_a$ and $Z_{ab}\cup X_b\to Z_b$.
\end{itemize}
\end{lem}

\proof Since the inclusions $X_{ab}\to X_a$ and $X_{ab}\to X_b$ are
cofibrations, we can easily reduce to the situation where $X_{ab}$
has collar neighborhoods $X_{ab}\times[-1,0]$ in $X_a$ and
$X_{ab}\times[0,1]$ in $X_b$. Ignoring condition (iv), we can easily
produce a factorization $f=gf_1$ with properties (i), (ii) and
(iii); then $Z$ has a filtration
\[ X=Z^{-1}\subset Z^0 \subset Z^1 \subset\cdots\subset Z^k=Z \]
where $Z^i$ is the relative $i$-skeleton. Suppose now
that $Z^{i-1}$ is already decomposed as in (iv). We can assume that
the attaching data for the cells we must attach to obtain $Z^i$ have
the form of a commutative diagram
\[
\CD
\coprod_{\alpha} S^{i-1} @>>>\coprod_{\alpha} D^i \\
@VVV  @VVV \\
Z^{i-1} @>>> E
\endCD
\]
such that the composition $\coprod D^i\to E\to X$ is transverse to the
subspace $X_{ab}\times\set{0}$ of $X_{ab}\times[-1,+1]\subset X$.
Triangulating the pair $(\coprod D^i,\coprod S^{i-1})$ in such a way that the
inverse image of $X_{ab}\times\set{0}$ is a subcomplex, we can also
arrange that the attaching map $S^{i-1}\to Z^{i-1}$ is cellular for
the chosen triangulation. (Here we are using the assumption that
$E$ is fibered over $X$.)
Using the triangulation cell structure on
the attached $\coprod D^i$, we obtain a ``new'' relative $CW$ structure on
$Z^i$ which extends the $CW$ structure on $Z^{i-1}$ and in which
$Z^i$ decomposes as in (iv). We continue inductively. \qed

\medskip
Keeping the notation of lemma~\ref{lem-retractivedeco}, we define an
$n$--dimensional \emph{quadratic structure on $(Z_a,Z_b,k)$}, for
$n\ge 0$, to be an element in $\Omega^n\Omega^{\infty}$ of the
homotopy pushout of the diagram
\[
\CD
((Z_{ab},k)\odot\lbul(Z_{ab},k))_{h\ZZ/2} @>>> ((Z_{b},k)\odot\lbul (Z_{b},k))_{h\ZZ/2} \\
@VVV   @. \\
((Z_{a},k)\odot\lbul(Z_{a},k))_{h\ZZ/2}\,. @.
\endCD
\]
Here the $\odot\lbul$ products are taken with respect to the base
spaces $X_a$, $X_{ab}$ and $X_b$, respectively. By a similar
generalization process, we arrive at the notions of a \emph{visible
symmetric} structure on $(Z_a,Z_b,k)$, and the notion of a
\emph{visible hyperquadratic} structure on $(Z_a,Z_b,k)$. The
corresponding abelian groups of path classes are denoted
\[ \begin{array}{ccc}
Q_n(Z_a,Z_b,k), & VQ^n(Z_a,Z_b,k), & V\widehat Q^n(Z_a,Z_b,k)
\end{array}
\]
respectively. They are actually defined for \emph{all} $n\in\ZZ$ as
$n$--th homotopy groups of the appropriate homotopy pushout spectra.
Clearly the map $g\co Z\to E$ in lemma~\ref{lem-retractivedeco}
induces homomorphisms from $Q_n(Z_a,Z_b,k)$ to $Q_n(E,k)$ and from
$VQ^n(Z_a,Z_b,k)$ to $VQ^n(E,k)$. We ask whether these homomorphisms
become isomorphisms ``in the (co)limit''. To speak of a colimit we
need an indexing category $\sC$, and in our case this should clearly
have objects of the form
\[  g\co Z\to E \]
where $Z$ is finite, $g$ is a weak equivalence of retractive spaces
over $X$ and $Z$ is decomposed into $Z_a$ and $Z_b$ with
intersection $Z_{ab}$ as before. Morphisms in the category, say from
$(g,Z,Z_a,Z_b)$ to $(g',Z',Z'_a,Z'_b)$, are cofibrations $u\co Z\to
Z'$ with $g'u=g$ and $u(Z_a)\subset Z'_a$, $u(Z_b)\subset Z'_b$.
Lemma~\ref{lem-retractivedeco} implies that $\sC$ is
\emph{directed}. That is, for any two objects $Z$ and $Z'$ in $\sC$
(in shorthand notation) there exists an object $Z''$ and morphisms
$Z\to Z''$, $Z'\to Z''$; and given any two morphisms $u,v\co Z\to
Z'$ there exists a morphism $w\co Z'\to Z'''$ such that $wu=wv$.

\begin{lem}
\label{lem-retractivedecowithstruc} In the above notation, we have
\[
\begin{array}{ccc}
\colimsub{(g,Z,Z_a,Z_b)}\, Q_n(Z_a,Z_b,k) & \stackrel{\cong}{\lra} & Q_n(E,k)\,, \\
\colimsub{(g,Z,Z_a,Z_b)}\, VQ^n(Z_a,Z_b,k) & \stackrel{\cong}{\lra} & VQ^n(E,k)\,, \\
\colimsub{(g,Z,Z_a,Z_b)}\, V\widehat Q^n(Z_a,Z_b,k) &
\stackrel{\cong}{\lra} & V\widehat Q^n(E,k)
\end{array}
\]
where the direct limits are taken over $\sC$.
\end{lem}

\proof The second isomorphism is a consequence of the first and the
third (and the ``five lemma''), because by
corollary~\ref{cor-tDalmostsplitting} there are exact sequences of
type
\[ \cdots \to V\widehat Q^{n+1} \to Q_n \to VQ^n \to
V\widehat Q^n \to Q_{n-1} \to \cdots \,\,.\] The third isomorphism
is obvious from definition~\ref{defn-nonlinearvishyper}. It remains
to establish the first isomorphism. Let $E_a:= E|X_a$, $E_b:=E|X_b$
and $E_{ab}=E|X_{ab}$. These are all fibered retractive spaces over
the appropriate base spaces: $X_a$, $X_a$ and $X_{ab}$,
respectively. Therefore, in forming the $\curlywedge$ product
$E_a\curlywedge E_a$\,, for example, we can proceed more directly
than otherwise by forming the fiberwise smash product of $E_a$ with
$E_a$ over $X_a$, and then dividing out by the zero section $X_a$.
This leads immediately to a homotopy pushout square consisting of
the four spaces $E_a\curlywedge E_a$, $E_b\curlywedge E_b$,
$E_{ab}\curlywedge E_{ab}$ and $E\curlywedge E$, and consequently a
homotopy pushout square of spectra
\[
\CD
((E_{ab},k)\odot\lbul(E_{ab},k))_{h\ZZ/2}  @>>> ((E_{a},k)\odot\lbul(E_{a},k))_{h\ZZ/2} \\
@VVV   @VVV  \\
((E_{b},k)\odot\lbul(E_{b},k))_{h\ZZ/2}  @>>>
((E,k)\odot\lbul(E,k))_{h\ZZ/2} \,.
\endCD
\]
where the $\odot\lbul$ products are taken with respect to the
appropriate base spaces: $X_{ab}$, $X_a$, $X_b$ and $X$. (Note that,
strictly speaking, the $\curlywedge$ product and the $\odot\lbul$
product have only been defined for finitely dominated retractive
spaces. We have no reason to think that $E_a$, $E_b$ and $E_{ab}$
are all finitely dominated, but the definition of $\curlywedge$
extends without difficulties.) We have therefore
\[   Q_n(E_a,E_b,k) \stackrel{\cong}{\lra} Q_n(E,k) \]
for all $n$, with the obvious interpretation of $Q_n(E_a,E_b,k)$.
This reduces our task to showing that
\[ \begin{array}{ccc}
\colimsub{(g,Z,Z_a,Z_b)}\, Q_n(Z_a,Z_b,k) & \stackrel{\cong}{\lra} &
Q_n(E_a,E_b,k)\,.
\end{array}
\]
By a Mayer--Vietoris and five lemma argument, this reduces further
to showing that the homomorphisms
\[ \begin{array}{ccc}
\colimsub{(g,Z,Z_a,Z_b)}\, Q_n(Z_a,k) & \stackrel{\cong}{\lra}
& Q_n(E_a,k) \\
\colimsub{(g,Z,Z_a,Z_b)}\, Q_n(Z_b,k) & \stackrel{\cong}{\lra}
& Q_n(E_b,k) \\
\colimsub{(g,Z,Z_a,Z_b)}\, Q_n(Z_{ab},k) & \stackrel{\cong}{\lra}
& Q_n(E_{ab},k) \\
\end{array}
\]
are all isomorphisms. By lemma~\ref{lem-retractivedeco}, we may now
enlarge the indexing categories to allow objects $(f,Y,Y_a,Y_b)$
where $Y$ is a finite retractive space over $X$ with a decomposition
$Y = Y_a\cup Y_b$ etc., and where $f\co Y\to E$ is any map (not
necessarily a weak equivalence) of retractive spaces over $X$. Then
$E_a$ for example can easily be identified with the homotopy direct
limit of the $Y_a$\,, etc., and $Q_*$ takes the homotopy direct
limits to direct limits, so that the isomorphisms become obvious.
\qed

\medskip
\begin{lem}
\label{lem-dualdeco} Let $X=X_a\cup X_b$ as in
lemma~\ref{lem-retractivedeco}. Let $Z$ be a finite retractive space
over $X$ with a decomposition $Z=Z_a\cup Z_b$ as in
lemma~\ref{lem-retractivedeco}, so that $Z_a$ is retractive over
$X_a$ and $Z_b$ is retractive over $X_b$. Let $k,n\in\ZZ$. Then
there exist a finite retractive space $V$ over $X$ with a
decomposition $V=V_a\cup V_b$ as in lemma~\ref{lem-retractivedeco},
an integer $\ell$ and a nondegenerate element $\eta$ in $\pi_n$ of
the homotopy pushout of
\[  (V_a,\ell)\odot\lbul(Z_a,k)\longleftarrow
(V_{ab},\ell)\odot\lbul(Z_{ab},k) \lra (V_b,\ell)\odot\lbul(Z_b,k)
\] such that the images of $\eta$ in
\[ \begin{array}{cc}
\pi_n((V,\ell)\odot\lbul(Z,k)), &
\pi_n((V_a,\ell)\odot\lbul(Z_a/Z_{ab},k)), \\
\pi_n((V_b,\ell)\odot\lbul(Z_b/Z_{ab},k)), &
\pi_{n-1}((V_{ab},\ell)\odot\lbul(Z_{ab},k))
\end{array}
\]
are all nondegenerate.
\end{lem}

\proof The guiding principle here is the fact that passage from
finite retractive spaces to cellular chain complexes over the
appropriate group(oid) rings respects and detects nondegenerate
pairings. This is due to \cite{Vogellduality}. The relevant
group(oid) rings here are $\ZZ[\pi_1(X_a)]$, $\ZZ[\pi_1(X_b)]$,
$\ZZ[\pi_1(X_{ab})]$ and $\ZZ[\pi_1(X)]$. Note also that a change of
rings, such as the passage from chain complexes over
$\ZZ[\pi_1(X_{ab})]$ to chain complexes over $\ZZ[\pi_1(X_{a})]$ by
means of
\[    \ZZ[\pi_1(X_{a})]\otimes_{\ZZ[\pi_1(X_{ab})]}  \]
respects nondegenerate pairings. It follows that cobase change, such
as the passage from retractive spaces over $X_{ab}$ to retractive
spaces over $X_a$ by means of
\[   X_a\sqcup_{X_{ab}}  \]
respects nondegenerate pairings. Consequently we can construct $V$
and $V_a$, $V_b$ in the following way. We first find an
$(n-1)$--dual for $(Z_{ab},k)$ as a stable retractive space over
$X_{ab}$. This amounts to finding a retractive space $V_{ab}$ over
$X_{ab}$, an integer $\ell$ and a nondegenerate element $\eta_{ab}$
in $\pi_{n-1}((V_{ab},\ell)\odot\lbul(Z_{ab},k))$. Next we find an
$n$--dual for the pair $((Z_a,k),(Z_{ab},k))$ which extends our
chosen $(n-1)$--dual for $(Z_{ab},k)$. This amounts to finding
$V_a$, a cofibration $V_{ab}$ and an element $\eta_a$ in $\pi_n$ of
the mapping cone of
\[ (V_{ab},\ell)\odot\lbul(Z_{ab},k) \lra (V_a,\ell)\odot\lbul(Z_a,k)  \]
whose image in $\pi_n(((V_a,\ell)\odot\lbul(Z_a/Z_{ab},k))$ is
nondegenerate and whose image in
$\pi_{n-1}((V_{ab},\ell)\odot\lbul(Z_{ab},k))$ is $\eta_{ab}$. (It
may be necessary to increase $\ell$.) We proceed similarly with the
pair $((Z_b,k),(Z_{ab},k))$ to obtain $V_b$ and $\eta_b$. Then we
define
\[ V := V_a\sqcup_{V_{ab}}V_b \]
and find $\eta$ in $\pi_n$ of the homotopy pushout of
\[  (V_a,\ell)\odot\lbul(Z_a,k)\longleftarrow
(V_{ab},\ell)\odot\lbul(Z_{ab},k) \lra (V_b,\ell)\odot\lbul(Z_b,k)
\] mapping to $\eta_a$ and $-\eta_b$ under the appropriate
projections. The existence of such an $\eta$ follows from a suitable
Mayer--Vietoris sequence. The image of $\eta$ in the homotopy group
$\pi_n((V,\ell)\odot\lbul(Z,k))$ will automatically be
nondegenerate, by Vogell's chain complex criterion. \qed

\bigskip
As another preliminary for the proof of
theorem~\ref{thm-VhatLexcision}, we offer a lengthy discussion of
how elements in the homotopy group
\[  V\widehat L^n(X) = \pi_n \bV\widehat\bL(X) \]
can be represented. For that discussion we return briefly to the
case of the $L$--theory of a group ring $\ZZ\pi$. Let $\sD$ be the
category of bounded (above and below) chain complexes of f.g. left
projective $\ZZ\pi$--modules. We regard $\sD$ as a Waldhausen
category in which the cofibrations are the chain maps which are
split injective in each dimension. Cofibrations $C\to D$ in $\sD$
are often regarded as pairs $(D,C)$. \newline Let $E$ be an object
of $\sD$ with an $n$--dimensional symmetric structure $\varphi$. The
inclusion of $E$ in the algebraic mapping cone of $\varphi_0\co
E^{n-*}\to E$ classifies an ``extension'' with base $E$ in the shape
of a short exact sequence
\[  0\to C\to D\to E\to 0 \]
where $C\simeq \Sigma^{-1}\cone(\varphi_0)$ and $D\simeq E^{n-*}$.
According to Ranicki \cite{RanickiLMS1}, the symmetric structure
$\varphi$ on $E$ has a preferred lift to an $n$--dimensional
\emph{nondegenerate} symmetric structure
$(\bar\varphi,\partial\bar\varphi)$ on the pair $(D,C)$, so that
$\varphi/\partial\bar\varphi=\varphi$ under the identification
$D/C\cong E$. (\emph{Warning:} $\bar\varphi$ is an $n$--chain in
$\hom_{\ZZ[\ZZ/2]}(W,D^t\otimes_{\ZZ\pi}D)$ with boundary
$\partial\bar\varphi$ in the image of
$\hom_{\ZZ[\ZZ/2]}(W,C^t\otimes_{\ZZ\pi}C)$. Our notation for
symmetric and quadratic structures on pairs deviates from
Ranicki's.) This resolution procedure of Ranicki's leads to a
bijective correspondence between homotopy types of chain complexes
$E$ with an $n$--dimensional symmetric structure, and homotopy types
of chain complex pairs $(D,C)$ with an $n$--dimensional
nondegenerate symmetric structure. There is a similar correspondence
for visible symmetric and quadratic structures. \newline Of
particular interest to us is the mixed case, i.e. the case of
nondegenerate symmetric pairs (or visible symmetric pairs) with
quadratic boundary. For a pair $(D,C)$ with an $n$--dimensional
nondegenerate symmetric structure
$(\bar\varphi,\partial\bar\varphi)$, improving the
$(n-1)$--dimensional symmetric structure $\partial\bar\varphi$ on
$C$ to a quadratic structure on $C$ amounts to ``trivializing'' the
induced $(n-1)$--dimensional hyperquadratic structure
\[  J(\partial\bar\varphi)\in \hom_{\ZZ[\ZZ/2]}(\widehat W,C^t\otimes C) \]
on $C$, in other words, finding an $n$--chain with boundary
$J(\partial\bar\varphi)$. But since the functor
\[  C\mapsto \hom_{\ZZ[\ZZ/2]}(\widehat W,C^t\otimes C) \]
respects homotopy cofiber sequences, finding such a trivialization
is equivalent to finding a trivialization for the suspended
hyperquadratic structure on $\Sigma C$. By Ranicki's correspondence,
if we write $E=D/C$ and $\varphi=\bar\varphi/\partial\bar\varphi$,
the suspension of $J(\partial\bar\varphi)$ can be identified with
the image of $J\varphi$ under the inclusion
\[   E\to \cone(\varphi_0) \,. \]
Summarizing, there is a bijective correspondence between homotopy
types of $n$--dimensional nondegenerate symmetric pairs in $\sD$
with quadratic boundary, and homotopy types of single objects $E$ in
$\sD$ with an $n$--dimensional symmetric structure $\varphi$ and a
trivialization $\tau$ of the image of $J\varphi$ under the inclusion
$E\to \cone(\varphi_0)$. This remains correct if symmetric and
hyperquadratic structures are replaced with visible symmetric and
visible hyperquadratic structures throughout. \newline The
correspondence extends to nullbordisms. For more precision we
suppose again that $((D,C),(\bar\varphi,\partial\bar\varphi))$ is an
$n$--dimensional nondegenerate symmetric pair with quadratic
boundary. Let $(E,\varphi)= (D/C,\bar\varphi/\partial\bar\varphi)$
be its Thom complex. Let $\tau$ be the trivialization of the image
of $J\varphi$ under $E\to \cone(\xi_0)$ determined by the preferred
trivialization of $J(\partial\bar\varphi)$. Then a (nondegenerate)
nullbordism of the nondegenerate pair
$((D,C),(\bar\varphi,\partial\bar\varphi))$ determines a nullbordism
of $(E,\varphi,\tau)$, by which is meant:
\begin{itemize}
\item a chain complex pair $(F,\partial F)$ with $(n+1)$--dimensional symmetric structure
$(\xi,\partial\xi)$, where $\partial F=E$ and
$\partial\xi=\varphi$~;
\item a trivialization $(\kappa,\partial\kappa)$ with $\partial\kappa=\tau$
of the image of $J(\xi,\partial\xi)$ under the inclusion
$(F,\partial F) \to (\cone(\xi_0),\cone(\partial\xi_0))$. Here
$\xi_0$ is regarded as a chain map from $\cone(F^{n-*}\to E^{n-*})$
to $F$.
\end{itemize}
\emph{Conversely}, a nullbordism of $(E,\varphi,\tau)$ determines a
nondegenerate nullbordism of the nondegenerate pair
$((D,C),(\bar\varphi,\partial\bar\varphi))$. Again, this remains
correct if symmetric and hyperquadratic structures are replaced with
visible symmetric and visible hyperquadratic structures throughout.
\newline Returning now to the stabilization $s\sR(X)$ of the
category of finitely dominated retractive spaces over a fixed space
$X$, we remark that these correspondences apply, mutatis mutandis,
in $s\sR(X)$. The fact that there are strictly speaking no
``canonical'' $n$--duals in $s\sR(X)$ does complicate matters
slightly (but only at first). For an object $(Y,k)$ in $s\sR(X)$
with an $n$--dimensional symmetric or visible symmetric structure
$\varphi$, the correct way to determine an $n$--dual $(Y,k-n)^*$
with $n\ge 0$ is to find an object $(Z,\ell)$ in $s\sR(X)$ and a
nondegenerate element $\eta\in \Omega^n((Z,\ell)\odot(Y,k))$. Modulo
a ``trivial'' enlargement of $(Y,k)$, a morphism $f\co (Z,\ell)\to
(Y,k)$ and a path from $f_*(\eta)$ to $\varphi_0\in
\Omega^n((Y,k)\odot(Y,k))$ can then be found. (The trivial
enlargement is an object of $s\sR(X)$ related to $(Y,k)$ by a
morphism which is both a cofibration and a weak equivalence.) The
morphism $f$ then deserves to be regarded as the adjoint of
$\varphi_0$. Therefore, when we write
\[ \varphi_0\co (Y,k-n)^*\lra (Y,k)\,, \]
we mean $f\co (Z,\ell)\to (Y,k)$.

\bigskip
\proof[Proof of theorem~\ref{thm-VhatLexcision}.] We assume that
$X=X_a\cup X_b$ and $X_a\cap X_b=X_{ab}$ as in
lemma~\ref{lem-retractivedeco}. We need to show that the gluing
homomorphism
\[  \alpha_n\co V\widehat L^n(X_a,X_b) \to V\widehat L^n(X) \]
is an isomorphism, where $V\widehat L^n(X_a,X_b)$ is the $n$--th
homotopy group of the homotopy pushout of
\[   \bV\widehat\bL\ubul(X_{a})\leftarrow\bV\widehat\bL\ubul(X_{ab})\to
\bV\widehat\bL\ubul(X_b)\,. \] We establish this only when $n>0$.
The case $n<0$ can be handled in the same way. (Replace $\odot$ by
$\odot_j$ for some $j$ with $j+n\ge 0$ in the argument below.)
\newline Starting with the surjectivity part and assuming $n>0$, we
represent an element of $V\widehat L^n(X)$ by an object $(Z,k)$ in
$s\sR(X)$, an $n$--dimensional visible symmetric structure $\varphi$
on $(Z,k)$ and a ``trivialization'' $\tau$ of the $n$--dimensional
visible hyperquadratic structure on $\cone(\varphi_0)$ obtained by
pushing $J\varphi$ forward along the inclusion
$(Z,k)\to\cone(\varphi_0)$. Here we view $\varphi_0$ as a morphism
from an $n$--dual of $(Z,k)$ to $(Z,k)$. We may assume that $Z$ is
not only finitely dominated, but finite. (Otherwise replace $Z$ by
$Z\wedge \Sigma Z$, which has zero finiteness obstruction; also,
replace $\varphi$ and $\tau$ by their images under appropriate maps
induced by the inclusion $Z\to Z\vee \Sigma Z$.) By
lemma~\ref{lem-retractivedeco} and
lemma~\ref{lem-retractivedecowithstruc} we may then assume that
\[ Z=Z_a\cup Z_b\,, \qquad \varphi=\varphi' + \varphi'' \]
where $Z_a$ and $Z_b$ are as in lemma~\ref{lem-retractivedeco} and
$(\varphi',\partial\varphi')$, $(\varphi'',\partial\varphi'')$ are
visible symmetric structures on the pairs $((Z_a,k),(Z_{ab},k))$ and
$((Z_b,k),(Z_{ab},k))$, respectively, with
$\partial\varphi'=-\partial\varphi''$. (In more detail, if $Z$ and
$\varphi$ do not come equipped with such a decomposition, then we
first use the Serre construction to enlarge $Z$ to a fibered
retractive space $E$ over $X$. The fibered retractive space $E$ can
in turn be approximated as in lemma~\ref{lem-retractivedeco} by
another retractive space $Z'$ over $X$ which is decomposed into
$Z'_a$ and $Z'_b$. Then lemma~\ref{lem-retractivedecowithstruc} can
be applied, etc.) The equations $\varphi=\varphi' + \varphi''$ and
$\partial\varphi'=-\partial\varphi''$ can be more accurately
expressed by saying that $\varphi$ is parametrized by $S^n$ and that
its restrictions to the ``upper'' and ``lower'' hemispheres of $S^n$
define $n$--dimensional visible symmetric structures on the pairs
$((Z_a,k),(Z_{ab},k))$ and $((Z_b,k),(Z_{ab},k))$, respectively.
Under these conditions, $\varphi$ and $\varphi'$, $\varphi''$
represent an element of what we have called $VQ^n(Z_a,Z_b,k)$ in
lemma~\ref{lem-retractivedecowithstruc}. \newline Now by
lemma~\ref{lem-dualdeco} we may assume that we have an $n$--dual
$(V,\ell)$ for $(Z,k)$ which is also decomposed, $V=V_a\cup V_b$.
Hence we have a decomposition in the shape of a pushout square
\[
\CD \cone((V_{ab},\ell)\,\,
\stackrel{\partial\varphi_0'\,\,\,\,}\lra (Z_{ab},k)) @>>>
\cone((V_b,\ell) \stackrel{\varphi_0''}{\lra} (Z_b,k)) \\
@VVV @VVV  \\
\cone((V_a,\ell) \stackrel{\varphi_0'}{\lra} (Z_a,k)) @>>>
\cone((V,\ell) \stackrel{\varphi_0}{\lra} (Z,k))
\endCD
\]
and the inclusion of $Z$ in $\cone(\varphi_0)$ respects the
decompositions. Hence, finally, the trivialization $\tau$
\emph{automatically} decomposes in the same manner. This completes
the solution of our decomposition problem and so establishes the
surjectivity part of the proof. \newline A relative version (which
we will not write out in detail) of the argument shows that, if the
original representative $((Z,k),\varphi,\tau)$ is nullbordant, in
the sense which we gave to the word ``nullbordant'' earlier, then
the lift across
\[  \alpha_n\co V\widehat L^n(X_a,X_b) \to V\widehat L^n(X) \]
which we have constructed is also nullbordant. Hence our
surjectivity proof amounts to a homomorphism of bordism groups which
is right inverse to $\alpha_n$. By a straightforward inspection, it
is also left inverse to $\alpha_n$. \qed

\section{The hyperquadratic $L$--theory of a point}
The $L$--theory of a point is, in our terminology, the $L$--theory
of the (stabilization of) the category of finite based $CW$--spaces
with the standard notion of Spanier--Whitehead duality. In this
chapter we ``calculate'' the homotopy types of the spectra
\[  \widehat\bL\ubul(\pt), \bV\widehat\bL\ubul(\pt) \,. \]
The calculations will not be used for anything else in this paper,
but they are interesting for a number of reasons. In particular we
shall see that the inclusion
\[
\bV\widehat\bL\ubul(\pt) \to \widehat\bL\ubul(\pt)
\]
is not a homotopy equivalence (which spoils the analogy with the
linear version of visible symmetric $L$--theory, outlined above).
But in fact it deviates very little from being a homotopy
equivalence, and the source turns out to be nothing more or less
than a cleaned--up version of the target.

\medskip
Understanding $\bV\widehat\bL\ubul(\pt)$ and $\widehat\bL\ubul(\pt)$
has a lot to do with understanding the ``homology theories''
\[  (Y,k)\mapsto\Sigma^{\infty-k}Y\,,\qquad
(Y,k)\mapsto ((Y,k)\odot\lbul(Y,k))^{th\ZZ/2}, \] and the natural
transformation from the first to the second which is implicit in
corollary~\ref{cor-tDalmostsplitting}. This natural transformation
can be made explicit by (re)defining $\Sigma^{\infty-k}Y$ as the
homotopy cofiber of the improved norm map
\[ ((Y,k)\odot\lbul(Y,k))_{h\ZZ/2} \lra ((Y,k)\odot\lbul(Y,k))^{\ZZ/2} \]
of proposition~\ref{prop-equifact} and
corollary~\ref{cor-tDalmostsplitting}, and (re)defining
$((Y,k)\odot\lbul(Y,k))^{th\ZZ/2}$ as the homotopy cofiber of the
ordinary norm map
\[ ((Y,k)\odot\lbul(Y,k))_{h\ZZ/2} \lra ((Y,k)\odot\lbul(Y,k))^{h\ZZ/2}\,.\]
Because we are dealing with homology theories, we can simplify
through a chain of natural weak homotopy equivalences,
\[
\begin{array}{ccc}
((Y,k)\odot\lbul(Y,k))^{th\ZZ/2} & \simeq \cdots \simeq &
(-k)\textup{--fold shift of } Y\wedge(S^0\odot\lbul S^0)^{th\ZZ/2}.
\end{array}
\]
and more obviously $\Sigma^{\infty-k}Y\simeq
Y\wedge\Sigma^{\infty-k}S^0$. This is explained in
\cite{WWassembly}. Now it is easy to identify $S^0\odot\lbul S^0$
with the sphere spectrum $\bS$ through a chain of equivariant
homotopy equivalences (using the flip action of $\ZZ/2$ on
$S^0\odot\lbul S^0$, and the trivial action of $\ZZ/2$ on $\bS$).
Hence what we need to understand is $\bS^{th\ZZ/2}$.

\medskip
The Segal conjecture \cite{Carlsson84} for a single point with the
(trivial) action of $\ZZ/2$ means that $\bS^{h\ZZ/2}$ is homotopy
equivalent to a certain completion of the fixed point spectrum of
the equivariant sphere spectrum $\bS_{\ZZ/2}$. The fixed point
spectrum of the equivariant sphere spectrum can be identified with
the $K$--theory of the symmetric monoidal category of finite
$\ZZ/2$--sets \cite{Carlsson92} and therefore breaks up as
\[ \bS \vee \bS_{h\ZZ/2} \]
where the summands correspond to the isomorphism types of
irreducible $\ZZ/2$--sets. If we identify the Burnside ring
$\pi_0(\bS \vee \bS_{h\ZZ/2})$ with $\ZZ \oplus \ZZ$, then the
augmentation ideal $\sI$ consists of the elements of the form
$(2z,-z)$. We have to complete at $\sI$. It is therefore to our
advantage to reconsider the splitting of the equivariant fixed point
spectrum: write
\[  (u(\bS_{\ZZ/2}))^{\ZZ/2} \simeq \bS\vee \Gamma(-\textup{transfer}) \]
where $\Gamma(-\textup{transfer})$ is the graph of the negative of
the transfer from $\bS_{h\ZZ/2}$ to $\bS$. Then the $\sI$ is the
$\pi_0$ of the second summand. Its powers are the ideals $2^n\sI$.
Therefore, indicating completion at $2$ by a left--hand superscript
$c$, we have
\[ \bS^{h\ZZ/2} \simeq \bS \vee \cp\Gamma(-\textup{transfer})\,. \]
In this decomposition, the norm map
\[  \bS_{h\ZZ/2} \to \bS^{h\ZZ/2} \]
has first component equal to the transfer and second component equal
to the identity (followed by completion). In calculating the
homotopy cofiber, we may replace the source by its $2$--completion
and $2$--complete the first summand of the target as well; the
homotopy cofiber remains the same. We summarize the result in the
following

\begin{lem} The homology theory
$(Y,k)\mapsto ((Y,k)\odot\lbul(Y,k))^{th\ZZ/2}$ has coefficient
spectrum $\cp\bS$. \qed
\end{lem}

Evaluating the natural transformation $\Sigma^{\infty-k}Y\to
((Y,k)\odot\lbul(Y,k))^{th\ZZ/2}$ just constructed on the object
$(Y,k)=(S^0,0)$, we have a map from $\Sigma^{\infty-k}Y=\bS$ to
$((Y,k)\odot\lbul(Y,k))^{th\ZZ/2} \simeq \cp\bS$.

\begin{lem}
The map under consideration is the inclusion $\bS \to \cp\bS$.
\end{lem}

\proof The homotopy fiber sequence of
corollary~\ref{cor-tDalmostsplitting} splits when $(Y,k)=(S^0,0)$.
Our map can therefore be obtained from the composition
\[ \bS\vee \bS_{h\ZZ/2}\,\simeq\,
(u(\bS_{\ZZ/2}))^{\ZZ/2}\to \bS^{h\ZZ/2} \to \bS^{th\ZZ/2}, \] which
we analyzed earlier, by restricting to the summand $\bS$. \qed

\medskip
We now recall, following \cite{WeissLMS}, how the chain bundle
method for determining hyperquadratic $L$--theory (and certain
variations on that) works in the linear case, and then transport the
technology to the nonlinear situation. \newline Let $R$ be a ring
with involution $R$, let $B$ be a bounded (below and above) chain
complex of f.g. projective left $R$--modules and let $\gamma$ be a
$0$--dimensional cycle in
\[ \hom_{\ZZ[\ZZ/2]}(\hat W,(B^{-*})^t\otimes_R B^{-*})
\,\,=\,\, ((B^{-*})^t\otimes_R B^{-*})^{th\ZZ/2}\,. \] Such a thing
is called a \emph{chain bundle} on $B$ and will be treated as a
chain complex analogue of a spherical fibration.
\newline
In particular, let $(C,\varphi)$ be a symmetric Poincar\ee chain
complex over $R$, of formal dimension $n$. Then $C$ has a
\emph{normal} chain bundle $\nu$, which comes together with an
$(n+1)$--chain $\tau$ in
\[ (C^t\otimes C)^{th\ZZ/2} \]
whose boundary is the difference between $J\varphi$ and
$(\varphi_0)_*(\Sigma^n\nu)$. Here $\Sigma^n\nu$ is the $n$--fold
homological suspension of $\nu$, an $n$--cycle in
$(C^{n-*})^t\otimes C^{n-*})^{th\ZZ/2}$. Because $\varphi_0$ is
invertible up to chain homotopy, the pair consisting of $\nu$ and
$\tau$ is sufficiently unique. \newline Given $B$ and a chain bundle
$\gamma$ on $B$, a \emph{$(B,\gamma)$}--structure on a symmetric
Poincar\ee chain complex $(C,\varphi)$ of formal dimension $n$
consists of a chain map $f\co C\to B$ and an $(n+1)$--chain $\tau$
in $(C^t\otimes_R C)^{th\ZZ/2}$ whose boundary is the difference of
$J\varphi$ and $(\varphi_0)_*\Sigma^n(f^*\gamma)$. The chain $\tau$
gives an identification of $f^*\gamma$ with the normal chain bundle
of $(C,\varphi)$. \newline Let $\bL\ubul(R;B,\gamma)$ be the
algebraic bordism spectrum constructed from the bordism theory of
symmetric algebraic Poincar\ee complexes $(C,\varphi)$ over $R$ with
a $(B,\gamma)$--structure. In the case where $B=0$ this is the
quadratic $L$--theory of $R$, and in the general case there is a
comparison map
\[ \bL\ubul(R) \to   \bL\lbul(R;B,\gamma) \]
with homotopy cofiber $\widehat\bL\lbul(R;B,\gamma)$. The main
result of \cite{WeissLMS} is a long exact sequence
\[
\CD \cdots @>>> \widehat L^n(R;B,\gamma) @>>> Q^n(B)@>J_{\gamma}>>
\widehat Q^n(B) @>>> \widehat L^{n-1}(R;B,\gamma) \cdots
\endCD
\]
where $J_{\gamma}[\varphi]:=
J[\varphi]-(\varphi_0)_*(\Sigma^n[\gamma])$ with
$\Sigma^n[\gamma]\in \widehat Q(B^{n-*})$. There is also a stunted
version of this. To get that we make the changes
\[  W \,\,\leadsto\,\, W_{\le 0} \,,\qquad\widehat W \,\,\leadsto\,\, \widehat W_{\le 0} \]
in the above (passing to $0$--skeletons). More practically we define
\[ Q^-_{n}(B):= H_n(W^-\otimes_{\ZZ[\ZZ/2]}(B^t\otimes B)) \]
where $W^-$ is the dual of $\widehat W_{\le 0}$, or alternatively,
the mapping cone of the (chain) map $W\to \ZZ[\ZZ/2]$ which takes
$1\in W_0$ to $1+T$. This gives an obvious inclusion-induced map
$\iota\co H_n(B^t\otimes_RB)\to Q^-_n(\Sigma B)$. The stunted
version of the above long exact sequence is another long exact
sequence
\[
\CD \cdots\widehat L^n(R;B,\gamma) @>>> H_n(B^t\otimes B)
@>\iota_{\gamma}>> Q^-_n(B) @>>> \widehat L^{n-1}(R;B,\gamma) \cdots
\endCD
\]
in which $\iota_{\gamma}[f]:= \iota[f]-f_*(\Sigma^n[\gamma])$ for a
chain map $f\co \Sigma^{n-*}B\to B$. The long exact sequences can be
set up as the homotopy group sequences associated with certain
homotopy fiber sequences of spectra. Beware that the crucial maps of
spectra which induce $J_{\gamma}$ and $\iota_{\gamma}$ are not
$\bH\ZZ$--module maps, although their sources and targets are
$\bH\ZZ$--module spectra. \newline Suppose now that $B=B$ and
$\gamma=\gamma(u)$ are universal; that is, the natural map
$H_0(\hom_R(C,B))\to Q^0(C^{-*})$ given by $f\mapsto f^*\gamma$ is
an isomorphism for every chain complex $C$ (bounded above and below,
f.g. projective in each degree). Then it is not hard to identify the
groups $L^n(R;B,\gamma)$ with the ordinary symmetric $L$--groups.
The above long exact sequence therefore specializes to
\[
\CD \cdots @>>> \widehat L^n(R) @>>> H_n(B^t\otimes B)
@>\iota_{\gamma}>> Q^-_n(B(u)) @>>> \widehat L^{n-1}(R) \cdots\,\,.
\endCD
\]
It must be said that this is a little harder to justify and use,
because in most cases $B$ can no longer be chosen to be bounded
above and below and f.g. in each degree. However, $B$ can always be
constructed as a direct limit of chain complexes $B'$ satisfying
these finiteness assumptions, and $\gamma$ can be constructed as an
element in the inverse limit of the chain bundle groups associated
with these (sub)complexes. The above long exact sequence for the
universal $B$ and $\gamma$ is then obtained as a direct limit for
the long exact sequences associated with the subcomplexes $B'$ and
chain bundles $\gamma|B'$. \newline In the case where $B$ is not
universal, we always have a comparison chain map to the universal
specimen.  We may think of $B$ as a classifying object for another
cohomology theory which comes with a natural transformation to
ordinary chain bundle theory.

\bigskip
Returning to finite spectra, we see that to calculate
$\widehat\bL\ubul(\pt)$ and $\bV\widehat\bL\ubul(\pt)$, we must
replace $B$ by $\bS^{th\ZZ/2}\simeq \cp\bS$ and by $\bS$,
respectively, in the above. Since $\bS$ has better finiteness
properties than $\cp\bS$, the visible case is easier and we begin
with that. By analogy with one of the long exact sequences just
described (the ``stunted version''), we obtain a homotopy fiber
sequence of spectra
\[
\CD \bV\widehat\bL\ubul(\pt) @>>> \bS\wedge \bS @>\iota_{\gamma}>>
(\bS\wedge\bS)^-_{h\ZZ/2}.
\endCD
\]
Here $(\bS\wedge\bS)^-_{h\ZZ/2}$ is the homotopy cofiber of the
transfer from $(\bS\wedge\bS)_{h\ZZ/2}$ to $\bS\wedge\bS$.
(\emph{Aside.} For the present purposes, the ``right'' notion of
smash product of two spectra $\bE$ and $\bF$ would be the spectrum
with $i$--th term $\Omega^i(E_i\wedge F_i)$, where the loop
coordinates are associated with the antidiagonal of
$\RR^i\times\RR^i$. This is ``commutative'' but neither associative
nor unital, so it is one of many naive smash products.)

\medskip
The spectrum $\bV\widehat\bL\ubul(\pt)$ is a ring spectrum. There is
no need to make that very precise here, since all we need is the
unit map for the ring structure $\bS\lra \bV\widehat\bL\ubul(\pt)$.
Evaluation on $\pi_0$ shows that this unit map is a homotopy right
inverse for the map $\bV\widehat\bL\ubul(\pt)\to \bS\wedge \bS$ in the
homotopy fiber sequence just above. Therefore $\iota_{\gamma}$ has a
preferred nullhomotopy and we have

\begin{thm}
\label{thm-VhatLsymofpoint} $\qquad\displaystyle
\bV\widehat\bL\ubul(\pt) \, \simeq \, \bS \vee \Omega(\bS^-_{h\ZZ/2})
\,=\, \bS \vee \RR\bP^{\infty}_{-1}\,. $
\end{thm}

\medskip
This is surprising. We have shown that the unit map for
$\bV\widehat\bL\ubul(\pt)$ is the injection of a wedge summand $\bS$,
up to homotopy equivalence. In particular, multiplication by 8 does
not annihilate its homotopy class. It follows immediately that the
standard homomorphism $L_0(\pt) \to VL^0(\pt)$ does \emph{not} send
the (signature 8) generator to 8 times the unit of $VL^0(\pt)$.

\begin{cor}
\label{thm-VhatLsymofX} For a space $X$ with $CW$--approximation
$X'\to X$, we have
\[ \bV\widehat\bL\ubul(X) \, \simeq \,X'_+\wedge\left( \bS \vee
\RR\bP^{\infty}_{-1}\right). \]
\end{cor}

\proof[Comment.] This is a formal consequence of
theorem~\ref{thm-VhatLsymofpoint} and the excision
theorem~\ref{thm-VhatLexcision}. Beware that the definition of of
$\bV\widehat\bL\ubul(X)$ which we use here depends on a specific $SW$
product in $s\sR(X)$. There are ``twisted'' versions which will be
considered later. \qed

\bigskip Next we calculate $\widehat\bL\ubul(\pt)$.
The only new aspect in this calculation is that our basic ``homology
theory'' is now $(Y,k)\mapsto ((Y,k)\odot\lbul(Y,k))^{h\ZZ/2}$ and
the representing object is $\cp\bS$. From the point of view of chain
bundle theory, generalized to the nonlinear setting, this means that
our calculation of the hyperquadratic $L$--theory of a point is
going to be almost identical with that of the visible hyperquadratic
$L$--theory of a point. The difference is that $\bS$ has to be
replaced by $\cp\bS$ where applicable. Noting that
$(\cp\bS\wedge\cp\bS)^-_{h\ZZ/2} \simeq (\bS\wedge\bS)^-_{h\ZZ/2}$~,
we obtain a homotopy fiber sequence of spectra
\[
\CD \widehat\bL\ubul(\pt) @>>> \cp\bS\wedge \cp\bS
@>\iota_{\gamma}>> (\bS\wedge\bS)^-_{h\ZZ/2}\,.
\endCD
\]
The map $\iota_{\gamma}$ in this case can immediately be understood
by comparison with the visible case. It must be zero because its
restriction to $\bS\wedge\bS$ is zero and the homotopy groups of the
target are all $2$--torsion. Therefore:

\begin{thm} $\qquad\displaystyle
\widehat\bL\ubul(\pt) \, \simeq \, (\cp\bS\wedge\cp\bS) \vee
\Omega(\bS^-_{h\ZZ/2})\, =\, (\cp\bS\wedge\cp\bS) \vee
\RR\bP^{\infty}_{-1}. $ \qed
\end{thm}

\section{Excision and restriction in controlled $L$--theory}
\label{sec-excision} We start with the Waldhausen category
$\sR\ld(\bar Q,Q)$ of \cite[dfn.7.1]{DwyerWeissWilliams}. Here $\bar Q$
is locally compact Hausdorff, $Q$ is open in $\bar Q$ and we add the
assumption that $\bar Q$ has a countable base. We recall that the objects of $\sR\ld(\bar Q,Q)$
are retractive spaces over $Q$ which are dominated (in a controlled homotopy sense) by locally
finite and finite dimensional retractive spaces with a controlled CW-structure.
\newline
Often we stabilize with respect to the suspension functor $\Sigma$ and write
the result as $s\sR\ld(\bar Q,Q)$. Objects in the stable category
can be written as $(Y,k)$ for some $Y$ in $\sR\ld(\bar Q,Q)$ and
$k\in\ZZ$. In the stabilized category, we want to introduce a
Spanier--Whitehead (external) product in the sense of
\cite[dfn.1.1]{WWduality}. (This has been done in
\cite[1.A.7]{WWduality}, but it will not hurt to present it from a
slightly different angle.)

\begin{defn}
\label{defn-strangequot} {\rm Let $\bar Q\ubul$ be the one-point compactification
of $\bar Q$. Let $Y$ be a retractive space over $Q$, with retraction
$r\co Y\to Q$. We write
\[ Y\cup_Q\bar Q\ubul \]
for the union of $Y$ and $\bar Q\ubul$ along $Q$, equipped with the coarsest topology
such that the inclusion $Y\to Y\cup_Q\bar Q\ubul$ embeds $Y$ as an open subset,
and the retraction $r\cup\id\co Y\cup_Q\bar Q\ubul\to \bar Q\ubul$ is continuous.
(This means that a subset $V$ of $Y\cup_Q\bar Q\ubul$ is a neighborhood
of some $z\in \bar Q\ubul\smin Q$ in $Y\cup_Q\bar Q\ubul$ iff $V$ contains
$(r\cup\id)^{-1}(W)$ for some neighborhood $W$ of $z$ in $\bar Q\ubul$.)
Let $Y/\!\!/Q$ be the topological quotient of $Y\cup_Q\bar Q\ubul$ by the
subspace $\bar Q\ubul$,
\[ Y/\!\!/Q = \frac{Y\cup_Q\bar Q\ubul}{\rule{0mm}{3.5mm}\bar Q\ubul}\,. \]
}
\end{defn}

\medskip\nin
\emph{Remark.} In the important special case where $Y$ has a locally
finite controlled $CW$-structure relative to $Q$, the special quotient
$Y/\!\!/Q$ can be described directly in terms of the ordinary quotient
$Y/Q$, which is a based $CW$-space. Namely, $Y/\!\!/Q$ is the topological
inverse limit of the based $CW$-spaces $Y/Y'$ where $Y'$ runs through
the cofinite based $CW$--subspaces of $Y$. (Here ``cofinite'' means that
$Y\smin Y'$ is a union of finitely many cells.) In general, the homotopy groups
of $Y/\!\!/Q$ should be regarded as ``locally finite''
variants of the homotopy groups of $Y/Q$.

\medskip\nin
\emph{Example.} Let $\bar Q=[0,1]$ and $Q=[0,1[\,$. Let
$T=\{1-2^{-i}\,|\,i=0,1,2,3,\dots\}$ and $Y=Q\amalg_TQ$. With the inclusion
of the first copy of $Q$ as the zero section, $Y$ becomes a retractive
space over $Q$. It has
an obvious locally finite controlled $CW$-structure relative to $Q$. The ordinary
quotient $Y/Q$ is a wedge of infinitely many circles. Its fundamental group
is free on generators $g_1,g_2,g_3,\dots$ corresponding to the 1-cells of $Y/Q$.
In particular it is countably infinite.
But $Y/\!\!/Q$ is homeomorphic to the Hawaiian earring. Its fundamental
group is an inverse limit of finitely generated free groups, and it is uncountable.
Similarly, for the suspension $\Sigma_QY$ (taken in the category of retractive
spaces over $Q$), we have
\[  \pi_2(\Sigma_QY/Q)\cong\bigoplus_{i=1}^{\infty}\ZZ~,\qquad
\pi_2(\Sigma_QY/\!\!/Q)\cong\prod_{i=1}^{\infty}\ZZ\,. \]

\medskip\nin
\emph{Remark.} Because of \cite{DwyerWeissWilliams} we are stuck with
the notation $(\bar Q,Q)$ for control spaces, even though we do not
require that $Q$ be dense in $\bar Q$. We will consequently try to
avoid the overline notation for topological closures. (The overline
notation is also used in section~\ref{sec-discrete} for something
completely unrelated.)

\begin{defn}
\label{defn-SWproduct} {\rm Let $Y$ and $Z$ be objects of
$\sR\ld(\bar Q,Q)$. To define their $SW$ product $Y\odot Z$, we
introduce first an unstable form $Y\curlywedge Z$ of it. We define
it as the geometric realization of a based simplicial set. An
$n$--simplex of this simplicial set is a pair $(f,\gamma)$ where
\begin{itemize}
\item[(i)] $f$ is a continuous map from
the standard $n$--simplex $\Delta^n$ to $Y/\!\!/Q\,\wedge Z/\!\!/Q$~;
\item[(ii)] $\gamma$ is a continuous assignment $c\mapsto \gamma_c$
of paths in $Q$, defined for $c\in \Delta^n$ with $f(c)$ not equal
to the base point $\pt\,$.
\end{itemize}
The paths $\gamma_c$ are to be parametrized by $[-1,+1]$ and must
satisfy $\gamma_c(-1)=r_Yf_Y(c)$ and $\gamma_c(+1)=r_Zf_Z(c)$, where
$r_Y, r_Z$ are the retractions and $f_Y(c)$, $f_Z(c)$ are the
coordinates of $f(c)$. Finally there is a control condition:
\begin{itemize}
\item[] For $z$ in $\bar Q\ubul\smin Q$ and any neighborhood $V$ of $z$ in
$\bar Q\ubul$, there exists a smaller neighborhood $W$ of $z$ in
$\bar Q\ubul$ such that, for any $c\in \Delta^n$ with $f(c)\ne \pt$,
the path $\gamma_c$ either avoids $W$ or runs entirely in $V$.
\end{itemize}
}
\end{defn}

\medskip

\begin{defn}
\label{defn-controlledSW} {\rm For $Y$ and $Z$ in $\sR\ld(\bar Q,Q)$
and integers $k,\ell\in\ZZ$, let
\[
(Y,k)\odot(Z,\ell) = \colim_n \,\Omega^{2n}(\Sigma^{n-k}Y\curlywedge
\Sigma^{n-\ell}Z)\,.
\]
More generally let $(Y,k)\odot\lbul(Z,\ell)$ be the
$\Omega$--spectrum with $j$--th space \[(Y,k)\odot_j(Z,\ell) =
\colim_n \,\Omega^{2n}\Sigma^j(\Sigma^{i-k}Y\curlywedge
\Sigma^{i-\ell}Z) \,.\] }
\end{defn}

\medskip
\nin\emph{Remark.} We have $\sR\ld(\bar Q,Q)= \sR\ld(\bar Q\ubul,Q)$
(an equality of Waldhausen categories). The meaning of $Y\odot Z$ is
the same in both categories. But there is a difference between
passage to germs near $\bar Q\smin Q$ (which we consider next) and
passage to germs near $\bar Q\ubul\smin Q$ (which we are not
interested in).

\bigskip
Next we work in the germ category $\sR\sG\ld(\bar Q,Q)$ of
\cite[dfn.7.1]{DwyerWeissWilliams} and its stable form,
$s\sR\sG\ld(\bar Q,Q)$. Let $Y$ and $Z$ be objects of
$\sR\sG\ld(\bar Q,Q)$. Note that $Y$ and $Z$ are honest retractive
spaces over $Q$. Again, to define their $SW$ product $Y\odot Z$ in
the germwise setting (recycled notation), we begin with an unstable
form $Y\curlywedge Z$ (also recycled notation) which is the
geometric realization of a simplicial set. An $n$--simplex in this
simplicial set is a germ of triples $(U,f,\gamma)$ where
\begin{itemize}
\item[(i)] $U=\bar U\cap Q$ for an open neighborhood $\bar U$
of $\bar Q\smin Q$ in $\bar Q$~;
\item[(ii)] $f$ is a continuous map from $\Delta^n$ to
$(Y_U/\!\!/U)\wedge(Z_U/\!\!/U)$, where $Y_U=r_Y^{-1}(U)$ and
$Z_U=r_Z^{-1}(U)$~;
\item[(ii)] $\gamma$ is a continuous assignment of paths ... (as before).
\end{itemize}
We impose the same control condition on $\gamma$ as
before. In (ii), we regard $Y_U$ and $Z_U=r_Z^{-1}(U)$ as retractive spaces
over $U$, and $U$ is the nonsingular part of the control space
$(\bar U,U)$. Note that $\bar U$ is the union of $U$ and
the singular set $\bar Q\smin Q$. (It is not defined as the closure
of $U$ in $\bar Q$.) Passage to germs is achieved by taking
the direct limit over all possible $U$. (It is a direct limit but the indexing is
contravariant, i.e., we approach it by making $U$ smaller and
smaller.)

\begin{defn}
\label{defn-germSWproducts} {\rm Put $(Y,k)\odot(Z,\ell) := \colim_n
\,\Omega^{2n}(\Sigma^{n-k}Y\curlywedge \Sigma^{j-\ell}Z)$. More
generally let $(Y,k)\odot\lbul(Z,\ell)$ be the $\Omega$--spectrum
with $j$--th space \[(Y,k)\odot_j(Z,\ell) = \colim_n
\,\Omega^{2n}\Sigma^j(\Sigma^{i-k}Y\curlywedge
\Sigma^{i-\ell}Z)\,.\] }
\end{defn}

\medskip
\nin\emph{Remark.} Later we will have to consider twisted versions
of the above, depending on a spherical fibration on $Q$.

\medskip
It is straightforward to verify that the above definitions of
$\odot$ and $\odot\lbul$ in the stable categories $s\sR\ld(\bar
Q,Q)$ and $s\sR\sG\ld(\bar Q,Q)$ satisfy the conditions of
\cite[\S1]{WWduality} for $SW$ products. It is less straightforward
to verify that they also satisfy the axioms of
\cite[\S2]{WWduality}, which are about existence and uniqueness of
``duals''. But this has been verified in \cite[\S2.A]{WWduality}.
Hence there are associated quadratic $L$-theory spectra
\cite{WWduality} which we denote by $\bL\lbul((\bar Q,Q))$ and
$\bL\lbul((\bar Q,Q)_{\infty})$, respectively. Also, visible symmetric
structures on objects of $s\sR\ld(\bar Q,Q)$ and $s\sR\sG\ld(\bar Q,Q)$
can be defined by analogy with definition~\ref{defn-nonlinearsymandvis}.
Hence there are visible symmetric $L$-theory spectra denoted
by $\bV\bL\ubul((\bar Q,Q))$ and $\bV\bL\ubul((\bar Q,Q)_{\infty})$,
respectively.

\medskip
We now specialize to the case $(\bar Q,Q)=\JJ
X=(X\times[0,1],X\times[0,1[\,)$ where $X$ is an ENR. In fact we
think of $X\mapsto L\ubul(\JJ X_{\infty})$ and $X \mapsto L\lbul(\JJ
X_{\infty})$ as covariant functors on the category $\sE\ubul$ whose
objects are the ENR's and where a morphism from $X_1$ to $X_2$ is a
based map $X_1\ubul\to X_2\ubul$ of the one--point
compactifications. (This is the same thing as a proper map from an
open subset of $X_1$ to $X_2$.)

\begin{thm}
\label{thm-homolothy} The spectrum valued functor $X\mapsto \bE(X)$,
where $\bE(X)$ means $\bL\lbul(\JJ X_{\infty})$, is homotopy
invariant and excisive. In detail:
\begin{itemize}
\item The projection from $X\times[0,1]$ to $X$ induces a homotopy
equivalence of $\bE(X\times[0,1])$ with $\bE(X)$.
\item
For an open subset $V$ of $X$, the collapse map $j\co X\ubul\to
V\ubul$ and the inclusion $i\co X\smin V\to X$ determine a homotopy
fiber sequence of spectra
\[\bE(X\smin V) \stackrel{i_*}{\lra} \bE(X) \stackrel{j_*}\lra
\bE(V). \]
\item For a disjoint union $X=\coprod_{i=1}^{\infty}X_i$ of ENR's, the
projections $X\ubul\to X_i\ubul$ induce an isomorphism
\[  \pi_*\bE(X) \lra \prod_{i=1}^{\infty}\pi_*\bE(X_i)~.\]
\end{itemize}
The coefficient spectrum $\bE(\pt)$ is homotopy equivalent to
$\Sigma\bL\lbul(\pt)$.
\end{thm}

\medskip
\nin\emph{Remark.} The spectrum $\bL\lbul(\pt)$ can be viewed as the
quadratic $L$--theory spectrum of the sphere spectrum, where the
latter is regarded as a (brave new) ring with involution. By the
$\pi$--$\pi$--theorem, $\bL\lbul(\pt)$ is homotopy equivalent to the
quadratic $L$--theory spectrum of the ring with involution $\ZZ$,
also known as the quadratic $L$--theory spectrum of the trivial
group. See \cite{WW2}.

\medskip
Statements similar to theorem~\ref{thm-homolothy} have been proved
in many places. See \cite{CarlssonPedersenVogell98} and \cite{ACFP},
for example. Our proof below imitates the proof of the analogous
statement for $A$--theory (= algebraic $K$--theory of spaces) in
\cite[\S\S 6--9]{Weissexci}. This requires some preparations.

\bigskip
We will work with based $CW$--spaces, which we generally view as
$CW$--spaces relative to $\pt$~. On the set of cells (not including
$\pt$) of such a $Y$, there is a partial ordering: $e\ge e'$ if the
smallest based $CW$--subspace containing $e$ also contains $e'$. We
say that $Y$ is \emph{dimensionwise locally finite}
\cite[Dfn.6.1]{Weissexci} if, for every cell $e$ in $Y$ (not
allowing $\pt$) and every $j\ge 0$ there are only finitely many
$j$--cells in $Y$ which are $\ge e$. For example, a wedge of
infinitely many based compact $CW$--spaces is dimensionwise locally
finite. \newline We replace the categories $\sR\sG\ld(\bar Q,Q)$ by
more tractable ones, denoted $\sR(\pt~;\bar Q,Q)_{\infty}$ in
\cite[\S6]{Weissexci}. An object of $\sR(\pt~;\bar Q,Q)_{\infty}$ is
a dimensionwise locally finite based CW--space $Y$ where the set of
cells (excluding $\pt$) is equipped with a map to $Q$. This map must
satisfy the usual control condition: given $n\ge 0$ and $z\in \bar
Q\smin Q$ and a  neighborhood $V$ of $z$ in $\bar Q$, there exists a
smaller neighborhood $W$ of $z$ in $\bar Q$ such that the closure of
any $n$-cell with label in $W$ is contained in a compact based
$CW$--subspace for which the cell labels are all in $V$. In
addition, for any $n\ge0$ and any compact region of $Q$, the set of
$n$-cells of $Y$ with labels in that compact region is required to
be finite. (There is also a finite domination condition to which we
return in a moment.) A morphism $Y\to Z$ is a sequence $(f_n)$ of
compatible cellular map germs between the skeletons, $f_n\co
Y^n_U\to Z^n$, where $f_n$ need only be defined on the cells of
$Y^n$ with labels in some open $U\subset Q$, where $U=\bar U\cap Q$
for some open neighborhood $\bar U$ of the singular set.
The maps $f_n$ are subject to
a straightforward control condition formulated in terms the of cell
labels. There is a good notion of ``controlled homotopy'' in the
category $\sR(\pt~;\bar Q,Q)_{\infty}$, so that the weak
equivalences in $\sR(\pt~;\bar Q,Q)_{\infty}$ can simply be defined
as the morphisms which are invertible up to controlled homotopy. The
cofibrations are, by definition, those morphisms $Y\to Z$ whose
underlying $CW$ map germ is a composition of $CW$ isomorphisms and
$CW$ subspace inclusions. It remains to make the finite domination
condition on objects $Y$ explicit. This is automatically satisfied
if $Y=Y^n$ for some $n$. In general it means that for some $n$ and
all $m\ge n$, the inclusion $Y^n\to Y^m$ admits a (controlled)
homotopy right inverse, so that $Y^m$ is a homotopy retract of $Y^n$
(in the ``germ'' sense). See \cite[\S6]{Weissexci} for more details.
\newline We come to the definition of $Y\curlywedge Z$ (again
recycled notation) for objects $Y$ and $Z$ of $\sR(\pt;\bar
Q,Q)_{\infty}$. Again this is defined as the geometric realization
of a simplicial set. An $n$--simplex in this simplicial set
corresponds to a germ of certain pairs $(U,f)$. Here
\begin{itemize}
\item $U=\bar U\cap Q$ for an open neighborhood $\bar U$ of
the singular set in $\bar Q$~;
\item $f$ is a continuous map from $\Delta^n$ to
$(Y_U/\!\!/U)\wedge(Z_U/\!\!/U)$.
\end{itemize}
Here $Y_U$ and $Z_U$ are the largest based $CW$-subspaces of $Y$ and
$Z$, respectively, containing only cells with labels in $U$. We
impose the usual control condition:
\begin{itemize}
\item[]
For $z\in \bar Q\smin Q$ and any neighborhood $V$ of $z$ in $\bar
Q$, there exists a smaller neighborhood $W$ of $z$ in $\bar Q$ such
that, for any $c\in \Delta^n$ with $f(c)\ne \pt$, either both $f_Y(c)$ and
$f_Z(c)$ are in cells with labels in $V$, or both are in cells with
labels outside $W$.
\end{itemize}
(Note the absence of ``paths''.) We pass to germs by taking the
direct limit over all possible $R$. Let
\[
\begin{array}{ccc}
(Y,k)\odot (Z,\ell)\,\,:=\,\,
\colim_n \,\Omega^{2n}(\Sigma^{n-k}Y\curlywedge \Sigma^{n-\ell}Z)\,, \\
(Y,k)\odot_j (Z,\ell)\,\,:=\,\, \colim_n
\,\Omega^{2n}\Sigma^j(\Sigma^{n-k}Y\curlywedge \Sigma^{n-\ell}Z)\,.
\end{array}
\]
Then $\odot$ and $\odot\lbul$ satisfy the axioms for an
$SW$--product listed in \cite{WWduality}.

\bigskip Now we specialize to the situation(s) where $(\bar Q,Q)=\JJ
X$ for some $X$ in $\sE\ubul$. We abbreviate $\bE'(X)=
\bL\lbul(\sR(\pt~;\bar Q,Q)_{\infty})$. Again we want to view the
assignment
\[ X \mapsto \bE'(X) \]
as a covariant functor on $\sE\ubul$. Indeed, every morphism $X_1\to
X_2$ in $\sE\ubul$, alias based map $f\co X_1\ubul\to X_2\ubul$, has
a factorization
\[ X_1\ubul \to V\ubul \to X_2\ubul \]
where $V=X_1\smin f^{-1}(\infty)$. In this factorization, the second
morphism $V\ubul\to X_2\ubul$ is induced by a proper map $V\to X_2$
and this determines in a straightforward way a map $\bE'(V)\to
\bE'(X_2)$. The first morphism $X_1\ubul\to V\ubul$ induces an exact
functor from $\sR(\pt~;\JJ X)_{\infty}$ to $\sR(\pt~;\JJ
V)_{\infty}$, hence a map $\bE'(X_1)\to \bE'(V)$, roughly as
follows. For an object $Y$ of $\sR(\pt~;\JJ X)_{\infty}$, the
largest based $CW$--subspace of $Y$ having all its cell labels in
$V\times[0,1\,[\,$ is an object of $\sR(\pt~;\JJ V)_{\infty}$.
\newline Later we will show, following \cite[\S9]{Weissexci}, that
$\bE'(X)$ is related to $\bE(X)$ in theorem~\ref{thm-homolothy} by a
chain of natural weak equivalences. But first we will prove the
analogue of theorem~\ref{thm-homolothy} for the functor $\bE'$.
\newline The main ingredients in this proof are certain
approximation statements, related to Waldhausen's approximation
theorem \cite{WaldRutgers}. To state these we fix $X$ (an ENR) and
an open $V\subset X$. On the category $\sR(\pt~;\JJ X)_{\infty}$ we
have, in addition to the standard notion of weak equivalence, a
coarser one denoted by $\omega$. Namely, a morphism is regarded as a
weak $\omega$--equivalence if the induced morphism in $\sR(\pt~;\JJ
V)_{\infty}$ is a weak equivalence. We write $\sR_{\omega}(\pt~;\JJ
X)_{\infty}$ for $\sR(\pt~;\JJ X)_{\infty}$ equipped with the coarse
notion  of weak equivalence. We write $\sR^{\omega}(\pt~;\JJ
X)_{\infty}$ for the full subcategory of $\sR(\pt~;\JJ X)_{\infty}$
consisting of the objects which are weakly $\omega$--equivalent to
the zero object, and this is equipped with the standard notion of
weak eqivalence inherited from $\sR(\pt~;\JJ X)_{\infty}$. \newline

\begin{lem}
\label{lem-app} The functors of stable categories determined by the
inclusion functor from $\sR(\pt~;\JJ(X\smin V))_{\infty}$ to
$\sR^{\omega}(\pt~;\JJ X)_{\infty}$ and the restriction from
$\sR_{\omega}(\pt~;\JJ X)_{\infty}$ to $\sR(\pt~;\JJ V)_{\infty}$
satisfy the hypotheses of Waldhausen's approximation theorem.
\end{lem}

\proof In the first case, the hypotheses are verified in
\cite[\S3,\S7]{Weissexci}, and this works even without
stabilization. In the second case, a closely related statement is
also proved in \cite[\S3,\S7]{Weissexci}, with more general
assumptions. Specialized to our situation this says that the induced
functor
\[ \sR_{\omega}(\pt~;\JJ X)\lf_{\infty}
\lra \sR(\pt~;\JJ V)\lf_{\infty} \] between the full
subcategories of finite dimensional objects satisfies the hypotheses
of the approximation theorem. Given that all weak equivalences in
the categories $\sR(\pt~;\JJ X)_{\infty}$ and $\sR(\pt~;\JJ
V)_{\infty}$ are invertible up to homotopy, it is easy to extend
this result from the full subcategories of finite dimensional
objects to the ambient categories, at the price of stabilizing, by
means of the next lemma. \qed

\begin{lem}
\label{lem-fd} Every object of $\sR(\pt~;\JJ X)_{\infty}$ becomes
weakly equivalent to a finite dimensional object after at most two
suspensions.
\end{lem}

\proof The excision and homotopy invariance theorem for the
algebraic $K$--theory functor $X\mapsto K(\sR(\pt~;\JJ X)_{\infty})$
is proved in \cite[7.1, 7.2]{Weissexci}. The coefficient spectrum is
analyzed in \cite[8.2,8.3]{Weissexci} and it is found to have a
vanishing $\pi_0$. In particular, the $K_0$ group of $\sR(\pt~;\JJ
X)_{\infty}$ is zero. Therefore, by standard finiteness obstruction
theory, all objects of $\sR(\pt~;\JJ X)_{\infty}$ are weakly
equivalent to finite dimensional ones after two suspensions. (This
is more fully explained in the proof of \cite[9.5]{Weissexci},
especially in the statement labelled $(**)$. The point is that the
general case can be reduced to the situation where an object is a
homotopy retract of another object whose cells are all concentrated
in one dimension. If that dimension is at least two, the homotopy
retraction can be linearized without any loss of information.) \qed

\medskip
For an object $Y$ of $\sR(\pt~;\JJ X)_{\infty}$ let $\mu(Y)$ be the
monoid of endomorphisms of $Y$ which are mapped to the identity by
the restriction functor from $\sR(\pt~;\JJ X)_{\infty}$ to
$\sR(\pt~;\JJ V)_{\infty}$. We note that, for objects $Y$ and $Z$ of
$\sR(\pt~;\JJ X)_{\infty}$, the product monoid $\mu(Y)\times\mu(Z)$
acts on the $SW$--product $Y\odot Z$. We like to think of
$\mu(Y)\times\mu(Z)$ as a category with one object. The action is a
functor on that category. Hence there is a canonical map
\[
\begin{array}{ccc}
\hocolimsub{\mu(Y)\times\mu(Z)} \, Y\odot Z & \lra & j_*(Y)\odot
j_*(Z)
\end{array}
\]
where $j_*\co \sR(\pt~;\JJ X)_{\infty}\to \sR(\pt~;\JJ V)_{\infty}$
is the restriction functor. \newline The next approximation lemma
about $SW$--products and its corollary (about quadratic structures)
are adaptations of \cite[2.7, 14.1]{RanickiLower}.

\begin{lem}
\label{lem-SWapp} For a finite dimensional object $Y$ of
$\sR(\pt~;\JJ X)_{\infty}$, the monoid $\mu(Y)$ is directed in the
following sense: given $f_1,f_2\in \mu(Y)$, there is $f_3\in \mu(Y)$
such that $f_3f_1=f_3=f_3f_2$. For two finite dimensional objects
$Y$ and $Z$ in $\sR(\pt~;\JJ X)_{\infty}$, the canonical map of
unstable $SW$ products
\[
\begin{array}{ccc}
\hocolimsub{\mu(Y)\times\mu(Z)} \, Y\curlywedge Z & \lra &
j_*(Y)\curlywedge j_*(Z)
\end{array}
\]
is a homotopy equivalence.
\end{lem}

\proof The statement about directedness is a consequence of the
following observations. For every $f\in \mu(Y)$, there exist a
neighborhood $U(f)$ of $X\times\set{1}$ in $X\times[0,1]$ and a
neighborhood $W(f)$ of $V\times\set{1}$ in $V\times[0,1]$, with
$U(f)\supset W(f)$, such that
\begin{itemize}
\item[] $f$ has a representative which is defined on every cell with label
in $U(f)$ and which is the identity on any cell $e$ whose label is
in $U(f)$ and whose image $f(e)$ has nonempty intersection with some
cell having label in $W(f)$.
\end{itemize}
Conversely, given any neighborhood $W$ of $V\times\set{1}$ in
$V\times[0,1]$, there exists $g\in \mu(Y)$ such that some
representative of $g$ is (undefined or) zero on all cells of $Y$
whose labels are not in $W$. (This is best proved by downward
induction on the dimension of $Y$. Assume $Y=Y^n$. Choose a
representative of an endomorphism of $Y^n/Y^{n-1}$ which is zero on
cells with labels outside $W$, and which belongs to
$\mu(Y^n/Y^{n-1})$. There is a smaller neighborhood $W'$ of
$V\times\set{1}$ in $V\times[0,1]$ such that this representative is
the identity on all cells with labels in $W'$. Next, choose a
representative of an endomorphism of $Y^{n-1}$ which belongs to
$\mu(Y^{n-1})$ and is zero on cells with labels outside $W'$. The
two representatives then combine to give an endomorphism of $Y^n=Y$
with the required property.) Combining these two observations, we
can choose $f_3\in \mu(Y)$ in such a way that it vanishes on all
cells with labels outside $W(f_1)\cap W(f_2)$, and then clearly
$f_3f_1=f_3=f_3f_2$. \newline Now for the statement about
$SW$--products: it is already clear from the foregoing that we have
an identification of (geometric realizations of) simplicial sets
\[
\begin{array}{ccc}
\colimsub{\mu(Y)\times\mu(Z)} Y\curlywedge Z & \cong &
j_*(Y)\curlywedge j_*(Z)
\end{array}
\]
where $\curlywedge$ denotes the unstable form of the $SW$--product.
As the colimit is a colimit of based $CW$--spaces and based cellular
maps over a directed category, we may replace it by a homotopy
colimit.  \qed

\medskip
\begin{cor}
\label{cor-murestriction} For a finite dimensional object $Y$ of
$\sR(\pt~;\JJ X)_{\infty}$ and $k\in\ZZ$, there is a canonical
homotopy equivalence of spectra
\[
\begin{array}{ccc}
\hocolimsub{\mu(Y)} ((Y,k)\odot\lbul(Y,k))_{h\ZZ/2} & \lra &
((j_*Y,k)\odot\lbul (j_*Y,k))_{h\ZZ/2}.
\end{array}
\]
\end{cor}

\bigskip
\proof[Proof of theorem~\ref{thm-homolothy}, excision part, with
$\bE'$ instead of $\bE$.] Writing $i_*$ for the inclusion functor
\[ \sR(\pt~;\JJ(X\smin V))_{\infty}
\lra \sR^{\omega}(\pt~;\JJ X)_{\infty} \] we have natural homotopy
equivalences $Y\curlywedge Z\lra i_*(Y)\curlywedge i_*(Z)$.
Consequently the homotopy classification of (nondegenerate)
quadratic structures is the same for an object of
$s\sR(\pt~;\JJ(X\smin V))_{\infty}$ and its image in
$s\sR^{\omega}(\pt~;\JJ X)_{\infty}$. Therefore and by the first
part of lemma~\ref{lem-app}, the map
\[ i_*\co\bL\lbul(\sR(\pt~;\JJ(X\smin V))_{\infty})
\lra \bL\lbul(\sR^{\omega}(\pt~;\JJ X)_{\infty}) \] is a homotopy
equivalence. For the rest of the argument, we use an $L$--theoretic
precursor, due to Ranicki, of Waldhausen's fibration theorem in
algebraic $K$--theory \cite{WaldRutgers}. Applied to our situation
this gives a homotopy fiber sequence of spectra
\[
\bL\lbul(\sR^{\omega}(\pt~;\JJ X)_{\infty}) \lra
\bL\lbul(\sR(\pt~;\JJ X)_{\infty}) \lra \bL\lbul(\sR(\pt~;\JJ
X)_{\infty},\sR^{\omega}(\pt~;\JJ X)_{\infty}) \] where
$\bL\lbul(\sR(\pt~;\JJ X)_{\infty},\sR^{\omega}(\pt~;\JJ
X)_{\infty})$ denotes the bordism theory of objects in the (stable
category of) $\sR(\pt~;\JJ X)_{\infty}$ equipped with a quadratic
structure which is nondegenerate modulo $\sR^{\omega}(\pt~;\JJ
X)_{\infty}$. See e.g. \cite[\S3]{RanickiTopMan} and
\cite{VogelNouvelle}~; see also remark~\ref{rem-Vogelstuff} below.
Therefore it only remains to show that the map
\[
\bL\lbul(\sR(\pt~;\JJ X)_{\infty},\sR^{\omega}(\pt~;\JJ X)_{\infty})
\lra \bL\lbul(\sR(\pt~;\JJ V)_{\infty})
\]
induced by $j$ is a homotopy equivalence. We verify that the induced
maps of homotopy groups $L_n(\dots)$ are isomorphisms for all $n\in
\ZZ$. For the surjectivity part, fix an object $(Y',k)$ in
$s\sR(\pt~;\JJ V)_{\infty}$ and an $n$--dimensional nondegenerate
quadratic structure $\psi'$ on it. By lemma~\ref{lem-app}, we may
assume that $Y'=j_*Y$ for some $Y$ in $\sR(\pt~;\JJ X)_{\infty}$. By
lemma~\ref{lem-SWapp}, there is an $n$--dimensional quadratic
structure $\psi$ on $(Y,k)$ such that $j_*\psi$ is homotopic to
$\psi'$. Then $\psi$ is automatically nondegenerate modulo
$s\sR^{\omega}(\pt~;\JJ X)_{\infty}$, so that $((Y,k),\psi)$
represents a class in $L_n(\sR(\pt~;\JJ
X)_{\infty},\sR^{\omega}(\pt~;\JJ X)_{\infty})$ which maps to the
class of $((Y',k),\psi')$ in $L_n(\sR(\pt~;\JJ V)_{\infty})$. For
the injectivity part, fix an object $(Z,k)$ in $s\sR(\pt~;\JJ
X)_{\infty}$ with a quadratic structure $\varphi$ which is
nondegenerate modulo $s\sR(\pt~;\JJ V)_{\infty}$, and assume that
$((Z',k),\varphi'):=((j_*Z,k),j_*\varphi)$ is nullbordant. Then
there exist a cofibration $u'\co (Z',k)\to (T',\ell)$ in
$s\sR(\pt~;\JJ V)_{\infty}$ and a nullhomotopy $\tau'$ of
$u'_*\varphi'$ such that $((Z',k)\to (T',\ell),
(\tau',\partial\tau'))$ with $\partial\tau'=\varphi'$ is a
nondegenerate quadratic pair in $s\sR(\pt~;\JJ V)_{\infty}$. By
lemma~\ref{lem-app}, we may assume that $u'$ is obtained from a
cofibration $u\co (Z,k)\to (T,\ell)$ in $s\sR^{\omega}(\pt~;\JJ
X)_{\infty}$ by applying $j_*$. By lemma~\ref{lem-SWapp}, on
composing $u$ with an appropriate endomorphism of $T$ (and restoring
the cofibration property by means of a mapping cylinder
construction), we may also assume that $\tau'$ is obtained from a
nullhomotopy $\tau$ for $u_*\varphi$ by applying $j_*$. Then
$((Z,k)\to (T,\ell), (\tau,\partial\tau))$ is a quadratic pair in
$s\sR(\pt~;\JJ X)_{\infty}$ which is nondegenerate modulo
$s\sR^{\omega}(\pt~;\JJ X)_{\infty}$. Hence $((Z,k),\varphi)$
represents the zero class. \qed

\proof[Proof of theorem~\ref{thm-homolothy}, homotopy invariance
part, with $\bE'$ instead of $\bE$.] It is enough to show that the
inclusion $i\co X\times\set{0}\to X\times[0,1]$ induces a homotopy
equivalence $i_*\co \bE'(X)\to \bE'(X\times[0,1])$. By the excision
property which we just established, it is also enough to show that
$\bE'(X\times\,]0,1])$ is contractible. This uses an Eilenberg swindle.
The details are as in \cite[\S4]{Weissexci}, except for a correction
to \cite[\S4]{Weissexci} in remark~\ref{rem-corr} below.
\qed

\medskip
\proof[Proof of theorem~\ref{thm-homolothy}, disjoint union axiom,
with $\bE'$ instead of $\bE$.] We leave this to the reader as a matter
of inspection.

\medskip
\proof[Proof of theorem~\ref{thm-homolothy}, coefficient spectrum
part, with $\bE'$ instead of $\bE$.] This is similar to the excision
part. We take $X=\pt$~. We introduce a Waldhausen category
$\sR(\pt~;\JJ\pt)$, defined like $\sR(\pt~;\JJ\pt)_{\infty}$ but
without the germ relation. Thus an object of $\sR(\pt~;\JJ\pt)$ is a
based $CW$--space with a map from the set of cells (excluding the
base point) to $[0,1[\,$. A morphism $Y\to Z$ in $\sR(\pt~;\JJ\pt)$
is a based cellular map (not a germ of such maps) from $Y$ to $Z$,
subject to the usual control condition. There is a finite domination
condition on objects $Y$, which says that for some $n$, each $Y^m$
with $m\ge n$ is a homotopy retract of $Y^n$ in the appropriate
controlled homotopy category. The weak equivalences are defined as
the morphisms which are invertible in the controlled homotopy
category. \newline In addition to the standard notion of weak
equivalence in $\sR(\pt~;\JJ\pt)$, we have a coarse notion $\omega$
of weak equivalence. Namely, a morphism in $\sR(\pt~;\JJ\pt)$ is a
weak $\omega$--equivalence if the induced morphism in
$\sR(\pt~;\JJ\pt)_{\infty}$ is a weak equivalence. As in the
excision part, we obtain from general principles a homotopy fiber
sequence of spectra
\[
\bL\lbul(\sR^{\omega}(\pt~;\JJ\pt)) \lra \bL\lbul(\sR(\pt~;\JJ\pt))
\lra \bL\lbul(\sR(\pt~;\JJ\pt),\sR^{\omega}(\pt~;\JJ\pt)).
\]
The Waldhausen category $\sR^{\omega}(\pt~;\JJ\pt)$ has an exact
subcategory consisting of those objects $Y$ which have only finitely
many cells. This is equivalent to the category of based finite
$CW$--spaces, so that its $L$--theory spectrum is $\bL\lbul(\pt)$.
The inclusion of this exact subcategory in
$\sR^{\omega}(\pt~;\JJ\pt)$ satisfies the conditions of the
approximation theorem; for the proof, see \cite[8.3]{Weissexci}. The
homotopy classification of (nondegenerate) quadratic structures on
an object in the (stabilized) subcategory is the same whether we
classify in the subcategory or in the ambient category. Consequently
we have
\[ \bL\lbul(\sR^{\omega}(\pt~;\JJ\pt))\,\,\simeq\,\, \bL\lbul(\pt) \,. \]
It remains to show that the ``passage to germs'' functor from
$\sR(\pt~;\JJ\pt)$ to $\sR(\pt~;\JJ\pt)_{\infty}$ induces a homotopy
equivalence of spectra
\[
\bL\lbul(\sR(\pt~;\JJ\pt),\sR^{\omega}(\pt~;\JJ\pt)) \lra
\bL\lbul(\sR(\pt~;\JJ\pt)_{\infty}) \] and this can be done by
considering the homotopy groups. We need to know that the functor of
stable categories determined by $\sR_{\omega}(\pt~;\JJ\pt)\to
\sR(\pt~;\JJ\pt)_{\infty}$ satisfies the conditions of the
approximation theorem; for this, see again \cite[8.3]{Weissexci} and
make use of lemma~\ref{lem-fd} above. The other ingredient is an
approximation lemma for quadratic structures analogous to
lemma~\ref{lem-SWapp}, but applicable to the ``passage to germs''
functor from $\sR(\pt~;\JJ\pt)$ to $\sR(\pt~;\JJ\pt)_{\infty}$. We
leave the remaining details to the reader. \qed

\proof[Proof of theorem~\ref{thm-homolothy}: comparing $\bE'$ and
$\bE$.] Recall that $\bE(X)=\bL\lbul(\sR\sG\ld(\bar Q,Q))$ and
$\bE'(X)=\bL\lbul(\sR(\pt~;\bar Q,Q)_{\infty})$ where $(\bar
Q,Q)=\JJ X$. The Waldhausen categories $\sR\sG\ld(\bar Q,Q)$ and
$\sR(\pt~;\bar Q,Q)_{\infty}$ are related, for a general control
space $(\bar Q,Q)$, by exact functors
\[
\CD \sR\sG\ld(\bar Q,Q) @<\textup{ inclusion }<< \sR\sG\lf(\bar Q,Q) @>v>> \sR(\pt~;\bar
Q,Q)_{\infty}\,.
\endCD
\]
Here $\sR\sG\lf(\bar Q,Q)$ is defined very much like $\sR\sG\ld(\bar Q,Q)$,
but the objects $Y$ come equipped with a finite dimensional controlled
$CW$--structure relative to $Q$ and morphisms are required to be
cellular relative to $Q$. See the proof of \cite[9.5]{Weissexci} for
details. (Except for a homotopy finiteness condition, which
is unimportant in our setting thanks to lemma~\ref{lem-fd}, the
category $\sR\sG\ld(\bar Q,Q)$ is identical with something denoted
$t\sR(\pt~;\bar Q,Q)_{\infty}$ in that proof, and
$\sR\sG\lf(\bar Q,Q)$ is denoted $\sB$ there.) It is also proved in
\cite[\S9]{Weissexci} that the two exact functors in the chain,
viewed as functors of the associated stable categories, satisfy the
conditions of the approximation theorem (again modulo
lemma~\ref{lem-fd}). Finally $v$ respects the $SW$-products,
in the strong sense that we have a
binatural homotopy equivalence
$Y\odot Z\to v(Y)\odot v(Z)$, for $Y$ and $Z$ in $\sR\sG\lf(\bar Q,Q)$. It follows
that $u$ and $v$ induce homotopy equivalences of the associated
quadratic $L$--theory spectra. \qed

\bigskip Our next goal is to formulate an excision
theorem similar to theorem~\ref{thm-homolothy} for
the visible symmetric $L$--theory. (We do not have, and
we do not need, an analogue of the excision theorem for a controlled
version of ordinary symmetric $L$--theory.) This is straightforward
modulo chapter~\ref{sec-visible}.

\begin{defn}
\label{defn-controllednonlinearvis} {\rm An \emph{$n$--dimensional
visible symmetric structure} on an object $(Y,k)$ in $\sR\sG\ld(\bar
Q,Q)$ is an element of $\Omega^n((Y,k)\odot(Y,k))^{\ZZ/2}$, with
$(Y,k)\odot(Y,k)$ as in definition~\ref{defn-germSWproducts}. An
\emph{$n$--dimensional visible symmetric structure} on an object
$(Y,k)$ in $\sR(\pt~;\bar Q,Q)_{\infty}$ is an element of
$\Omega^n((Y,k)\odot(Y,k))^{\ZZ/2}$, with the appropriate definition
of $\odot\lbul$ for the category $\sR(\pt~;\bar Q,Q)_{\infty}$. }
\end{defn}

\medskip
Again $(Y,k)\odot\lbul(Y,k)$ turns out to be the underlying spectrum
of a $\ZZ/2$--spectrum which we can describe as a shifted suspension
spectrum
\[  \bS_{\ZZ/2}^{-kW}\wedge Y^{\curlywedge 2}. \]
(The meaning of $Y^{\curlywedge 2}=Y\curlywedge Y$ depends on the
category, which may be $\sR\sG\ld(\bar Q,Q)$ or $\sR(\pt~;\bar
Q,Q)_{\infty}$.) The analogues of
corollary~\ref{cor-tDalmostsplitting} hold in both categories, and
they are still corollaries of proposition~\ref{prop-equifact}.
\newline Let $Y$ be an object of $\sR\sG\ld(\bar Q,Q)$ or
$\sR(\pt~;\bar Q,Q)_{\infty}$. Let $\bar U$ be an open neighborhood
of the singular set in $\bar Q$ and put $U=\bar U\cap Q$.
Recall that $Y_U$ means $r^{-1}(U)$ for $Y$ in
$\sR\sG\ld(\bar Q,Q)$~; for $Y$ in $\sR(\pt~;\bar Q,Q)_{\infty}$ it
means the largest based $CW$-subspace of $Y$ having all its cell
labels in $U$. In the second case we also introduce the notation $Y_U/\!\!/\pt$
for the topological inverse limit of the based spaces $Y_U/Y_U'$ where $Y_U'$
runs though all cofinite based $CW$-subspaces of $Y_U$.

\begin{cor}
\label{cor-germtDalmostsplitting} In the setting of
definition~\ref{defn-controllednonlinearvis}, there is a natural
homotopy fiber sequence of spectra
\[
\CD ((Y,k)\odot\lbul(Y,k))_{h\ZZ/2} @>>>
((Y,k)\odot\lbul(Y,k))^{\ZZ/2} @>J>>
\!\!\hocolimsub{U}\!\! \Sigma^{\infty-k}(Y_U/\!\!/U).
\endCD \]
\end{cor}

\begin{cor}
\label{cor-moregermtDalmostsplitting} For any object $Y$ of
$\sR(\pt~;\bar Q,Q)_{\infty}$ and $k\in\ZZ$, there is a homotopy
fiber sequence of spectra
\[
\CD ((Y,k)\odot\lbul(Y,k))_{h\ZZ/2} @>>>
((Y,k)\odot\lbul(Y,k))^{\ZZ/2} @>J>>
\hocolimsub{U}\!\! \Sigma^{\infty-k}(Y_U/\!\!/\pt).
\endCD \]
\end{cor}

In these two corollaries, $U$ runs through the open subsets of $Q$
of the form $U=\bar U\cap Q$ where $\bar U$ is an open neighborhood
of the singular set in $\bar Q$.
The homotopy colimits are reduced (taken in the based
category) and the suspension spectrum construction $\Sigma^{\infty-k}$
is meant to include a $CW$-approximation mechanism.

\medskip
There is also an analogue of corollary~\ref{cor-murestriction}. We
keep the assumptions and notation of that corollary to state the
analogue:

\begin{cor}
\label{cor-newmurestriction} For a finite dimensional object $Y$ of
$\sR(\pt~;\JJ X)_{\infty}$ and $k\in\ZZ$, there is a canonical
homotopy equivalence of spectra
\[
\begin{array}{ccc}
\hocolimsub{\mu(Y)} ((Y,k)\odot\lbul(Y,k))^{\ZZ/2} & \lra &
((j_*Y,k)\odot\lbul (j_*Y,k))^{\ZZ/2}.
\end{array}
\]
\end{cor}

\proof By corollaries~\ref{cor-murestriction}
and~\ref{cor-moregermtDalmostsplitting}, and a five lemma argument,
it is enough to verify that the canonical map
\[
\begin{array}{ccc}
\hocolimsub{\mu(Y)}\hocolimsub{U}\Sigma^{\infty-k}(Y_U/\!\!/U) &
\lra & \hocolimsub{W}\Sigma^{\infty-k}((j_*Y)_W/\!\!/\pt)
\end{array}
\]
is a homotopy equivalence. But this is obvious.\qed

\bigskip
With these tools available, the
visible symmetric $L$-theory version of theorem~\ref{thm-homolothy},
which we are about to state,
can be proved in strict analogy with
the original (quadratic $L$--theory) version.

\begin{thm}
\label{thm-visiblehomolothy} The spectrum valued functor $X\mapsto
\bE(X)$, where $\bE(X)$ means $\bV\bL\ubul(\JJ X_{\infty})$, is
homotopy invariant and excisive. The coefficient spectrum $\bE(\pt)$
is homotopy equivalent to $\Sigma \bV\bL\ubul(\pt)$. \qed
\end{thm}

\section{Control and visible $L$-theory}
\label{sec-controlvisible}
In this section our goal is to generalize
theorem~\ref{thm-VhatLexcision} to a controlled setting, as far as
possible.

\medskip
Fix a compact Hausdorff space $S$. For most of this section the only
control spaces we shall be interested in are of the form $(\bar
X,X)$, with compact $\bar X$ and an identification $\bar X\smin
X\cong S$. The only morphisms $f\co (\bar X_1,X_1)\to (\bar
X_2,X_2)$ between such control spaces that we shall be interested in
are those which are relative to $S$. These objects and morphisms
form a category $\sK^S$. \newline We can also speak of homotopies
between morphisms in $\sK^S$. These will also be relative to $S$,
and they allow us to define a homotopy category $\sH\sK^S$. A
morphism in $\sK^S$ is a \emph{cofibration} if it is an embedding
which has the homotopy extension property, for such homotopies.

\begin{thm}
\label{thm-controlledVhatLexcision} On $\sK^S$, the functor
$(\bar X,X) \mapsto \bV\widehat\bL\ubul((\bar X,X))$
is homotopy invariant, excisive and satisfies a strong ``wedge'' axiom.
\end{thm}

This needs a few explanations. The homotopy invariance property
means that the functor takes homotopy equivalences in
$\sK^S$ to homotopy equivalences of spectra. Here \emph{homotopy
equivalences in $\sK^S$} refers to morphisms in $\sK^S$ which become
invertible in $\sH\sK^S$. \newline
For the excision property, suppose given a pushout diagram
\[
\CD
(\bar X_{ab},X_{ab}) @>>> (\bar X_a,X_a) \\
@VVV @VVV \\
(\bar X_b,X_b) @>>> (\bar X,X)
\endCD
\]
in $\sK^S$ where all the arrows are cofibrations (and without loss
of generality, all are inclusions). It is being claimed that in such a case
\[
\CD
\bV\widehat\bL\ubul((\bar X_{ab},X_{ab})) @>>> \bV\widehat\bL\ubul((\bar X_a,X_a)) \\
@VVV @VVV \\
\bV\widehat\bL\ubul((\bar X_b,X_b)) @>>> \bV\widehat\bL\ubul((\bar X,X))
\endCD
\]
is homotopy cocartesian, and also that $\bV\widehat\bL\ubul$ applied
to the initial object $(S,\emptyset)$ of $\sK^S$
is a weakly contractible spectrum. \newline
For the strong wedge
axiom, suppose that $(\bar X,X)$ is in $\sK^S$ and $X$ is a topological
disjoint union of subspaces $X_{\alpha}$, where $\alpha$ runs through
some (countable) set. Let $\bar X_{\alpha}$ be the union of $X_{\alpha}$
and the singular set $\bar X\smin X$. Then the embedding
\[ (\bar X_{\alpha},X_{\alpha})\to (\bar X,X) \]
is a morphism in $\sK^S$~,
for every $\alpha$. It follows from the ordinary excision property
(just above) that the induced homomorphisms
\[ V\widehat L^n((\bar X_{\alpha},X_{\alpha}))\lra V\widehat L^n((\bar X,X)) \]
are split injective with a preferred splitting, since $(\bar X,X)$
is the coproduct in $\sK^S$
of $(\bar X_{\alpha},X_{\alpha})$ and $(\bar X,X\smin X_{\alpha})$.
The decomposition of $(\bar X,X)$ into the
$(\bar X_{\alpha},X_{\alpha})$ could be regarded as a generalized wedge
decomposition (because it is when $S$ is a point). It is not in general
a coproduct decomposition. Nevertheless, it is being claimed that
the projections
$V\widehat L^n((\bar X,X))\to
V\widehat L^n((\bar X_{\alpha},X_{\alpha}))$
induce an isomorphism
\[ V\widehat L^n((\bar X,X)) \lra
\prod_{\alpha} V\widehat L^n((\bar X_{\alpha},X_{\alpha})). \]

\bigskip
We turn to the proofs. The homotopy invariance property in
theorem~\ref{thm-controlledVhatLexcision} can be proved by the
same argument as the homotopy invariance property in
theorem~\ref{thm-VhatLexcision}. The excision property for the
special case of a coproduct
(that is, $(\bar X,X)$ in $\sK^S$ with $X=X_1\amalg X_2$)
and the strong ``wedge axiom''
are valid by inspection. This leaves the general excision
property. It is correct to say that the proof of the excision
property in theorem~\ref{thm-VhatLexcision} carries over, but some
clarifications are nevertheless in order. The difficulty is that
in lemma~\ref{lem-retractivedeco}, lemma~\ref{lem-retractivedecowithstruc}
etc., which were part of the proof of theorem~\ref{thm-VhatLexcision},
we made essential use of the concept of \emph{fibration}. Here we
will need a corresponding concept of \emph{controlled fibration} and this
is not completely obvious.

\medskip
Fix a control space $(\bar X,X)$ with compact $\bar X$ for simplicity,
but not necessarily in $\sK^S$. Let $p\co E\to X$ be any map. Suppose that
a collection of open subspaces $E_\lambda\subset E$ has been specified,
where $\lambda$ runs though a directed set; suppose also that the
indexing is monotone, so that $\lambda_1<\lambda_2$ implies $E_{\lambda_1}
\subset E_{\lambda_2}$. We assume that $E$ is the union of the $E_{\lambda}$.

\begin{defn}
\label{defn-controlledfibration}     
 {\rm The map $p\co E\to X$ together with the directed
system of open subspaces $\{E_{\lambda}\}$ such that
$E=\bigcup E_{\lambda}$ is a \emph{controlled Serre fibration system} if the following
holds. For every controlled finite dimensional locally finite $CW$-space
$Y$ over $X$, every controlled map $f\co Y\to E_{\lambda}$ and
controlled homotopy $h\co Y\times[0,1]\to X$ starting with $p\circ f$~, there
exists $\kappa>\lambda$ and a homotopy $Y\times[0,1]\to E_{\kappa}$
which lifts $h$ and starts with $f$.
}
\end{defn}

\medskip\nin
\emph{Remark.} Controlled fibration systems can be pulled back along
maps of control spaces. More precisely, if $f\co (\bar X_1,X_1)\to (\bar X_2,X_2)$
is a map of control spaces (compact $\bar X_1$ and $\bar X_2$), and
$\big(p\co E\to X_2,\{E_{\lambda}\}\big)$ is a controlled Serre fibration
system over $X_2$, then the pullback $f^*E$ with the subspaces $f^*E_{\lambda}$ is a
controlled Serre fibration system over $X_1$.

\begin{lem}
\label{lem-controlledSerre}
Let $Z$ be a retractive space over $X$, with retraction $r\co Z\to X$.
Then there exists a controlled Serre fibration system
$\big(p\co E\to X,\{E_{\lambda}\}\big)$
and an embedding $Z\to \bigcap E_{\lambda}$ over $X$, inducing controlled
homotopy equivalences $Z\to E_{\lambda}$ for all $\lambda$.
\end{lem}

\proof This is supposed to be a controlled version of the Serre construction
in ordinary fibration theory. Let $W$ be an open neighborhood of the diagonal
in $X\times X$, invariant under permutation of the two factors $X$.
We say that $W$ is \emph{controlled} if its closure in
$\bar X\times\bar X$ is disjoint from $X\times(\bar X\smin X)$.
(Equivalently, $W$ is controlled if, for every
$x\in \bar X\smin X$ and every neighborhood $V$ of $x$ in $\bar X$, there
exists a smaller neighborhood $V'$ of $z$ in $\bar X$ such that for
any $(y_1,y_2)\in W$ with $y_1\in V'$ we have $y_2\in V$.) Ordered by inclusion,
these controlled neighborhoods form a directed system $\Lambda$.
For $W$ in $\Lambda$ let $E_W$ be the space of pairs $(z,\omega)$ where
$z\in Z$ and $\omega\co [0,1]\to X$ is a path in $X$ with $\omega(0)=r(z)$, subject to
the control condition determined by $W$. (Namely, all points of the form
$(\omega(s),\omega(t))$ are in $W$.) Let $p_W\co E_W\to X$ be defined
by $p_W(z,\omega)=\omega(1)$. The inclusion $Z\to E_W$ is clear. It is a map
over $X$ and as such it is a controlled homotopy equivalence.
The inclusions $E_W\to E_{W'}$ for $W'>W$ in $\Lambda$ are also clear.
We let $E=\bigcup E_W$. The remaining details are left to the reader. \qed

\bigskip
Returning to $\sK^S$ and the proof of theorem~\ref{thm-controlledVhatLexcision}, it
only remains to say that lemma~\ref{lem-dualdeco} carries over to the
controlled situation without essential changes. In the controlled version
of lemma~\ref{lem-retractivedeco}, the fibration $E\to X$ should be replaced
by a controlled Serre fibration system as in definition~\ref{defn-controlledfibration};
the morphisms $f$ should land in some $E_{\lambda_1}$ by assumption and the morphism
$g$ should be constructed to land in some $E_{\lambda_2}$ where
$\lambda_2\ge \lambda_1$. In the controlled version of
lemma~\ref{lem-retractivedecowithstruc}, the groups $Q_n(E;k)$,
$VQ^n(E;k)$ and $V\widehat Q^n(E,k)$ should be replaced by the direct
limits over $\lambda$ of $Q_n(E_{\lambda};k)$,
$VQ^n(E_{\lambda};k)$ and $V\widehat Q^n(E_{\lambda},k)$, respectively.
Constructions like $E_a=E|X_a$ should be read as restrictions or pullbacks
of controlled fibration systems. These two lemmas (in the controlled version)
feed into the proof of theorem~\ref{thm-controlledVhatLexcision} via
the controlled Serre construction, lemma~\ref{lem-controlledSerre}.
The proof then proceeds as in the situation without control.
\qed

\section{Control, stabilization and change of decoration} \label{sec-controlsusp} 
Our main point in this chapter is to make connections between the $L$-theory and algebraic $K$-theory
of a control space $(\bar Q,Q)$ on one side, and the $L$-theory and algebraic $K$-theory of
$(\bar Q\times\bar\RR,Q\times\RR)$ on the other side, where $\bar\RR=[-\infty,+\infty]$.
We are going to be more specific, as follows. Let $X$ be a compact Hausdorff space. Fix an
integer $i\ge 0$ and form the control space $(X*S^{i-1},X*S^{i-1}\smin S^{i-1})$.
We usually identify $X*S^{i-1}\smin S^{i-1}$ with $X\times\RR^i$ and so write
$(X*S^{i-1},X\times\RR^i)$. Less formally still, we like to refer to $X\times\RR^i$ only
and use a letter $c$ as in $\bL\lbul(X\times\RR^i\cc)$ to indicate that
$X\times\RR^i$ is the nonsingular part of a control space $(X*S^{i-1},X\times\RR^i)$.
Our goal is then to compare constructions like $\bL\lbul(X\times\RR^i\cc)$ and $\bA(X\times\RR^i\cc)$
with $\bL\lbul(X\times\RR^{i+1}\cc)$ and $\bA(X\times\RR^{i+1}\cc)$. \newline
In particular, there is an exact functor
$\times\RR$ from retractive spaces over $X\times\RR^i$ (subject to various conditions) to
retractive spaces over $X\times\RR^{i+1}$. Because of the geometric applications that we have in mind,
it is important for us to come to grips with the induced maps in controlled $L$-theory and $A$-theory.
From an algebraic viewpoint, this is not the best starting point. Instead
we start by using fibration theorems from algebraic $K$-theory and $L$-theory to
see connections between the controlled algebraic $K$ or $L$-theory of $X\times\RR^i$ and the
controlled algebraic $K$ or $L$-theory of $X\times\RR^{i+1}$. These connections are of the well-known type.
We have suppressed some of the more mechanical details in the proofs.

\medskip
Let $\bL\lbul^h(X\times\RR^i\cc)$ be the controlled
$L$-theory spectrum of $(X*S^{i-1},X\times\RR^i)$, constructed
using locally finite and finite dimensional retractive spaces (with a
controlled relative $CW$-structure) over $X\times\RR^i$ rather
than locally finitely dominated ones. Similarly,
$\bA^h(X\times\RR^i\cc)$ is the controlled $A$-theory
constructed using locally finite and finite dimensional retractive
spaces (with a controlled relative $CW$-structure) over $X\times\RR^i$.

\begin{thm}
\label{thm-changeofdeco} We have
\[
\begin{array}{ccc}
\bL\lbul(X\times\RR^i\cc) & \simeq &
\Omega \bL\lbul^h(X\times\RR^{i+1}\cc), \\
\bV\bL\ubul(X\times\RR^i\cc) & \simeq &
\Omega {\bV\bL\ubul}^h(X\times\RR^{i+1}\cc), \\
\bA(X\times\RR^i\cc) & \simeq &
\Omega \bA^h(X\times\RR^{i+1}\cc).
\end{array}
\]
Indeed there is a homotopy cartesian square of inclusion-induced maps
\begin{equation} \label{eqn-Lsuspend}
\begin{split}
\xymatrix{
\bL\lbul(X\times\RR^i\cc) \ar[r] \ar[d] & \ar[d] \bL\lbul(X\times\RR^i\times~]-\infty,0]\cc) \\
\bL\lbul(X\times\RR^i\times[0,\infty[~\cc) \ar[r] & \bL\lbul^h(X\times\RR^{i+1}\cc)
}
\end{split}
\end{equation}
with contractible off-diagonal terms; and there are analogous homotopy cartesian squares for $\bV\bL\ubul$ and $\bA$.
\end{thm}

\proof  We concentrate on the quadratic
$L$-theory case, the other two cases being very similar.
The first step is to replace the control conditions by stronger ones.
Let $p\co X\times\RR^i\to \RR^i$ be the projection.
Given retractive spaces $Y_1$ and $Y_2$
over $X\times\RR^i$ and a controlled map $f\co Y_1\to Y_2$ (which
is understood to be relative to $X\times\RR^i$ but need not respect the
retractions to $X\times\RR^i$), we say that $f$ is \emph{bounded}
if there exists a real number $a\ge0$ such that
$\|p(y)-pf(y)\|\le a$ for all $y\in Y_1$. Similarly, there is a
notion of bounded and controlled homotopy between bounded and controlled
maps (between retractive spaces over $X\times\RR^i$).
A morphism $f\co Y_1\to Y_2$ between
retractive spaces over $X\times\RR^i$ (i.e., a map which respects
both the zero sections and the retractions) is a bounded map for
trivial reasons; we call it a \emph{weak equivalence} (in the
bounded sense) if it is invertible up to bounded
homotopies relative to $X\times\RR^i$. A finite dimensional
controlled $CW$--structure on a retractive space $Y$ over $X\times\RR^i$~,
relative to $X\times\RR^i$~, is \emph{bounded} if there exists $a>0$
such that the image of each cell in $\RR^i$ has diameter $\le a$.
We use all that to introduce a
Waldhausen category
\[ \sR\ld(X\times\RR^i;b) \]
similar to $\sR\ld(X\times\RR^i\cc)$~, but with all
control conditions replaced by the corresponding boundedness
conditions. There is also a stable version
\[ s\sR\ld(X\times\RR^i;b) \]
obtained by adjoining formal desuspensions. There is also a
preferred $SW$-product on $s\sR\ld(X\times\RR^i;b)$, which can be
constructed roughly as in
definition~\ref{defn-SWproduct}. (Replace $Q$ there by
$X\times\RR^i$. Replace the control condition on $\gamma$ there by
the condition that, for some $a>0$, the images in $\RR^i$ of
all paths $\gamma_c$ have diameter $\le a$.) \newline
There is an inclusion of Waldhausen categories
\[ s\sR\ld(X\times\RR^i;b) \lra s\sR\ld(X\times\RR^i\cc)\,.\]
This is compatible with the $SW$-products and preserves nondegenerate
pairings. (The claim about nondegeneracy can be reduced to the case of a
paring between two objects which have all their cells in a single dimension.)
Hence we have an induced map in $L$-theory
\[
\begin{array}{ccc}
\bL\lbul(s\sR\ld(X\times\RR^i;b)) & \lra &
\bL\lbul(s\sR\ld(X\times\RR^i\cc)). \\
\end{array}
\]
It is a key fact, and one whose proof we want to skip, that
this map is a homotopy equivalence. See \cite{ACFP}.
Similarly there is an inclusion map
\[
\begin{array}{ccc}
\bL\lbul(s\sR\lf(X\times\RR^{i+1};b)) & \lra &
\bL\lbul(s\sR\lf(X\times\RR^{i+1}\cc)). \\
\end{array}
\]
of $L$-theory spectra, where the superscript \emph{lf} refers to retractive
spaces with a locally finite, finite dimensional and bounded $CW$-structure
relative to $X\times\RR^{i+1}$. This is again a homotopy equivalence. \newline
For the remainder of this proof we work in the bounded setting.
We use the Waldhausen fibration theorem \cite{WaldRutgers}
or rather its analogue in $L$-theory due to
Ranicki/Vogel. The category to which we will apply it is
\[ s\sR\lf(X\times\RR^{i+1};b) \]
and for the purposes of this proof we abbreviate this to $s\sR$.
In addition to the default notion of weak
equivalence in $s\sR$, we introduce two coarser notions involving
passage to germs. Let $Y_1$ and  $Y_2$ be retractive
spaces over $X\times\RR^{i+1}$. A $u$-germ of bounded and controlled maps
from $Y_1$ to $Y_2$
is represented by a bounded and controlled map from the portion of $Y_1$ lying over
$X\times\RR^i\times[a,\infty[~$, for some $a\in\RR$, to $Y_2$. Two such representatives
define the same $u$-germ if they are both defined on $X\times\RR^i\times[a',\infty[~$
for some $a'\in\RR$ and agree there. Similarly, a $v$-germ of
bounded and controlled maps from $Y_1$ to $Y_2$
is represented by a bounded and controlled map from the portion of $Y_1$ lying over
$X\times\RR^i\times\,]-\infty,a[~$, for some $a\in\RR$, to $Y_2$.
Call a morphism in $s\sR$ a $u$-equivalence if its mapping cone
is weakly equivalent to the zero object in the $u$-germ sense.
Call it a $v$-equivalence if its mapping cone
is weakly equivalent to zero in the $v$-germ sense. Using the notation which Waldhausen
uses in his formulation of the fibration theorem, we obtain a
commutative diagram of Waldhausen categories with SW-duality
\[
\CD
s\sR^u\cap s\sR^v @>>> s\sR^u  \\
@VVV   @VVV  \\
s\sR^v @>>> s\sR \,.
\endCD
\]
In the resulting commutative diagram of $L$-theory spectra
\begin{equation} \label{eqn-Lpresuspend}
\begin{split}
\CD
\bL\lbul(s\sR^u\cap s\sR^v) @>>> \bL\lbul(s\sR^u) @>>> \bL\lbul(s\sR^u,s\sR^u\cap s\sR^v) \\
@VVV   @VVV @VVV \\
\bL\lbul(s\sR^v) @>>> \bL\lbul(s\sR) @>>> \bL\lbul(s\sR,s\sR^v)
\endCD
\end{split}
\end{equation}
the rows are homotopy fibration sequences by the fibration theorem.
The terms $\bL\lbul(s\sR^u)$ and $\bL\lbul(s\sR^v)$ are contractible
because the Waldhausen categories involved are flasque. (More precisely
there is an Eilenberg swindle argument for contractibility, as follows.
The translation $x\mapsto x-n$
acting on $\RR$ induces an endofunctor $\kappa_n$ of $s\sR^u$.
The sum
\[ \tau= \bigvee_{n\ge 0} \kappa_n \]
is again an endofunctor of $s\sR^u$. While it is not strictly true that
$\tau\cong\id\vee\tau$, it is easy to relate $\tau$ and $\id\vee\tau$
by a chain of natural equivalences, preserving $SW$-duality. Hence the
identity map of $\bL\lbul(s\sR^u)$ is contractible.) Next, the right-hand
vertical arrow in diagram~(\ref{eqn-Lpresuspend}) is a homotopy equivalence by
an easy application of the approximation theorem. (What we have in mind
here is an $L$-theoretic cousin of Waldhausen's approximation theorem
in algebraic $K$-theory. This is very neatly formulated in \cite{VogelNouvelle}.)
Hence the left-hand square in diagram~(\ref{eqn-Lpresuspend}) is homotopy cartesian and so
\[ \bL\lbul(s\sR^u\cap s\sR^v)\simeq \Omega\bL\lbul(s\sR).\]
Finally the projection $X\times\RR^i\times\RR\to X\times\RR^i$
induces a map from
$\bL\lbul(s\sR^u\cap s\sR^v)$ to $\bL\lbul((X*S^{i-1},X\times\RR^i))$ which is a
homotopy equivalence by another (but easier) application of the approximation
theorem. This leads to a situation where we can map the $b$-variant ($b$ for bounded) of
the commutative square~(\ref{eqn-Lsuspend})
to the commutative square~(\ref{eqn-Lpresuspend}) by a map which is a termwise homotopy
equivalence. (Here we adopt a generous interpretation of
$s\sR\lf(X\times\RR^{i+1};b)$ as the full Waldhausen subcategory of $s\sR\ld(X\times\RR^{i+1};b)$
consisting of all objects which can be related to locally finite and finite dimensional objects
by a chain of weak equivalences.) Indeed, on lower right-hand
terms our map is an identity map; on upper left-hand terms it is a homotopy equivalence by the
observation about $\bL\lbul(s\sR^u\cap s\sR^v)$ just made; and the remaining terms, in both squares,
are contractible. Therefore the square~(\ref{eqn-Lsuspend}) is also homotopy cartesian.
\qed

\begin{cor} \label{cor-stabLandVL}
The maps
\[ \begin{array}{ccc}
\times\RR\co \bL\lbul(X\times\RR^i\cc) &\lra& \Omega\bL^h\lbul(X\times\RR^{i+1}\cc)~, \\
\times\RR\co \bVL\ubul(X\times\RR^i\cc) &\lra& \Omega{\bVL\ubul}^h(X\times\RR^{i+1}\cc)
\end{array}
\]
induced by the exact functor $\times\RR$ are homotopy equivalences.
\end{cor}

\begin{rem} \label{rem-apologies}
{\rm Before proving corollary~\ref{cor-stabLandVL} we need to clarify its meaning.
The notation $\times\RR$ is self-explanatory to the extent that it describes an exact functor from
$s\sR\ld(X\times\RR^i\cc)$ to $s\sR\lf(X\times\RR^{i+1}\cc)$, with the usual generous interpretation of the
$lf$ superscript. In addition, we need to know that
it takes $n$-dualities to $(n+1)$-dualities. More precisely, there is a binatural transformation
\[   Y\odot_n Z \lra (Y\times\RR)\odot_{n+1}(Z\times\RR), \]
for $Y$ and $Z$ in $s\sR\ld(X\times\RR^i\cc)$, which commutes with the symmetry actions of $\ZZ/2$ and
preserves nondegenerate pairings.
In this section we use only some formal properties of $\times\RR$, which is why we defer a more detailed
description to sections~\ref{sec-hocharsig} and~\ref{sec-excicharsig}. From the point of view of section~\ref{sec-hocharsig},
the two maps in corollary~\ref{cor-stabLandVL} are given by external product with
$\sigma(\RR)\in \Omega^{\infty+1}{\bVL\ubul}^h(\RR\cc)$, the
controlled (visible symmetric) signature of $\RR$ as a manifold with control map $\id\co\RR\to\RR$.  \newline
We also need to clarify the meaning of the $\Omega$ prefix in corollary~\ref{cor-stabLandVL}.
Let $\bE$ be a spectrum, e.g.~in the sense of \cite{Adams74}, with $n$-th space $E_n$.
It is generally a good idea to make a distinction between the spectrum $\Omega\bE$
whose $n$-th space is $\Omega E_n$ and the spectrum $\bE[-1]$ whose $n$-th space is $E_{n-1}$. The canonical maps
$E_{n-1}\to \Omega E_n$ do not fit together to make a spectrum map $\bE[-1]\to \Omega\bE$.
The two spectra can nevertheless be related by a chain of weak homotopy equivalences. Also,
it is easy to identify $\Omega^{\infty+m}\bE$ with $\Omega^\infty\bE[-m]$ because it is easy to
identify $\bS^0[m]$ with $S^m\wedge \bS^0$. In any case,
we do in this paper sometimes write $\Omega\bE$, or $\Omega^m\bE$, when we ought to write $\bE[-1]$
or $\bE[-m]$. It was difficult to avoid. In the statement of theorem~\ref{thm-changeofdeco}, the prefix
$\Omega$ is an honest $\Omega$. By contrast in the statement of corollary~\ref{cor-stabLandVL} it should
be read as a shift operator. \newline
}
\end{rem}

\proof[Proof of corollary~\ref{cor-stabLandVL}] We concentrate on the quadratic $L$-theory case.
Despite remark~\ref{rem-apologies}, it is meaningful to say that the map $\times\RR$ which we are
discussing induces the same homomorphism on homotopy groups as the map
\[  \bL\lbul(X\times\RR^i\cc) \lra
\Omega \bL\lbul^h(X\times\RR^{i+1}\cc) \]
resulting from the homotopy cartesian square~(\ref{eqn-Lsuspend}). This is enough and it is
what we shall prove. Start with a representative $(Y,\varphi)$ of an element $\alpha_0$
of $\pi_n\bL\lbul(X\times\RR^i\cc)$, so that $Y$ is a retractive space over $X\times\RR^i$.
In order to see the corresponding element $\alpha_1$ of $\pi_{n+1}\bL\lbul(X\times\RR^{i+1}\cc)$, according to~(\ref{eqn-Lsuspend}),
we proceed as follows. We choose an $(n+1)$-dimensional nullbordism for $(Y,\varphi)$
as a quadratic Poincar\ee object in \[\sR\ld(X\times\RR^i\times~]-\infty,0]\cc)\]
and another for
$(Y,\varphi)$ as a quadratic Poincar\ee object in
\[\sR\ld(X\times\RR^i\times[0,\infty[~\cc)\]
We glue the two together along the common boundary $(Y,\varphi)$ to obtain an $(n+1)$-dimensional quadratic
Poincar\ee object in $\sR\ld(X\times\RR^{i+1}\cc)$. This represents $\alpha_1$. Now the most obvious choices
for the two $(n+1)$-dimensional nullbordisms are $(Y,\varphi)\times~]-\infty,0]$ and $(Y,\varphi)\times[0,\infty)$,
where we are using informal notation. Then the representative for $\alpha_1$ which we get is
$(Y,\varphi)\times\RR$.
\qed

\begin{lem}
\label{lem-degenerateprod} The map $\times\RR\,\co \bA(X\times\RR^i\cc)\lra
\bA(X\times\RR^{i+1}\cc)$
is nullhomotopic.
\end{lem}

\proof
For each retractive space $Y$ over $X\times\RR^i$ we have a
cofibration sequence of retractive spaces over $X\times\RR^{i+1}$, as follows:
\[   Y_0 \lra Y\times\RR \lra Y_\lambda\amalg Y_\rho~. \]
Here $Y_0$ is the union of $Y\times 0$ and $X\times\RR^i\times\RR=X\times\RR^{i+1}$;
it is clear that the cofiber of the inclusion $Y_0\to Y\times\RR$, taken in the category of retractive spaces
over $X\times\RR^{i+1}$, breaks up into a coproduct of two retractive spaces
$Y_\lambda$ and $Y_\rho$ over $X\times\RR^{i+1}$, in such a way that $Y_\lambda$
is trivial over $X\times\RR^i\times[0,\infty[~$ and $Y_\rho$ is trivial over
$X\times\RR^i\times~]-\infty,0]$.
By the additivity theorem, the map $\times\RR$ is therefore homotopic (by a preferred homotopy) to the sum
of three maps induced by the exact functors taking $Y$ to $Y_0$, $Y_\lambda$ and $Y_\rho$~,
respectively. The functors taking $Y$ to $Y_0$ and $Y_\lambda$ induce maps from
$\bA(X\times\RR^i\cc)$ to $\bA(X\times\RR^{i+1}\cc)$ which clearly factor through the contractible spectrum
\[  \bA(X\times\RR^i\times~]-\infty,0]\cc). \]
The functor taking $Y$ to $Y_\rho$ induces a map from
$\bA(X\times\RR^i\cc)$ to $\bA(X\times\RR^{i+1}\cc)$ which factors through the contractible spectrum
$\bA(X\times\RR^i\times[0,\infty[~\cc)$. \qed

\begin{cor} \label{cor-ultimateAstab} There is a commutative diagram of spectra with action of $\ZZ/2$,
\[
\xymatrix@R=50pt@C=30pt{
\bA(X\times\RR^i,n\cc) \ar[r]^-{\times\RR} \ar[d]^-{\textup{incl.}} & \bA^h(X\times\RR^{i+1},n+1\cc) \\
S^1_!\wedge \bA(X\times\RR^i,n\cc) \ar@{..>}[ur]^-{\simeq}
}
\]
where $S^1_!$ denotes the based space $S^1$ with the conjugation action of $\ZZ/2$ (fixed point set $S^0$)
and $S^1_!\wedge \bA(X\times\RR^i,n\cc)$ has the diagonal action of $\ZZ/2$.
\end{cor}

\proof[Interpretation and proof.] The vertical arrow is induced by the inclusion $S^0\to S^1_!$, which is
a $\ZZ/2$-map with the trivial action of $\ZZ/2$ on $S^0$. Of course
we identify $\bA(X\times\RR^i,n\cc)$ with $S^0\wedge \bA(X\times\RR^i,n\cc)$. For the horizontal arrow,
we interpret $\bA^h(X\times\RR^{i+1},n+1\cc)$ as the $0$-connected cover
of $\bA(X\times\RR^{i+1},n+1\cc)$. Therefore the map $\times\RR$ of lemma~\ref{lem-degenerateprod}
factors through $\bA^h(X\times\RR^{i+1},n+1\cc)$, with $\ZZ/2$-action.
(Indeed, it induces the zero homomorphism on $\pi_0$ as it is nullhomotopic.) The dashed arrow
remains to be constructed. \newline
A statement equivalent to the above is that the map $\times\RR$, top horizontal arrow, admits a
(non-equivariant) nullhomotopy $H$
such that the concatenation of $H$ and the conjugate $tHt$ (where $t\in\ZZ/2$ is the generator)
produces a map
\[ \bA(X\times\RR^i,n\cc) \lra \Omega\bA^h(X\times\RR^{i+1},n+1\cc) \]
which is a homotopy equivalence.
To establish that, it is enough to produce a nullhomotopy $H$ such that the concatenation of $H$ and $tHt$
produces a map from $\bA(X\times\RR^i,n\cc)$ to $\Omega\bA^h(X\times\RR^{i+1},n+1\cc)$ which is homotopic to
the map of theorem~\ref{thm-changeofdeco} (up to sign). We are going to show that the nullhomotopy $H$ constructed in
lemma~\ref{lem-degenerateprod} has this property. \newline
We are comparing two (non-equivariant) maps from $\bA(X\times\RR^i,n\cc)$ to
\[ \Omega\bA^h(X\times\RR^{i+1},n+1\cc)~\simeq~\Omega\bA(X\times\RR^{i+1},n+1\cc) ~. \]
By inspection, both maps can be described by means of external products. The first (concatenation
of $H$ and $tHt$) is given by external product with an element $a$ of $\pi_1\bA(\RR,1\cc)$ and the second
(from theorem~\ref{thm-changeofdeco}) is given by product with an element $b$ of $\pi_1\bA(\RR,1\cc)$.
To describe $a$, we use the map
\[ \CD \bA(\pt,0) @>\times\RR>> \bA(\RR,1\cc) \endCD \]
of spectra with action of $\ZZ/2$. From lemma~\ref{lem-degenerateprod}, we have a
preferred nullhomotopy $H$ for it. The concatenation of $H$ and $tHt$ induces
a map from $\bA(\pt)$ to $\Omega\bA(\RR\cc)$, and so a homomorphism from
$\pi_0\bA(\pt)$ to $\pi_1\bA(\RR\cc)$. Let $a$ be the image of $1\in \pi_0\bA(\pt)$
under that map. Define $b$ as the image of $1\in \pi_0\bA(\pt)$ under
the isomorphism $\pi_0\bA(\pt,0) \cong \pi_1\bA(\RR,1\cc)$
implied by theorem~\ref{thm-changeofdeco}. Now it remains only to show that $a=\pm b$.
This is an exercise in duality which we leave to the
reader. \qed

\begin{cor} \label{cor-moreAstab}
The following is homotopy cartesian,
\[
\xymatrix{
{\bA(X\times\RR^i,n\cc)\hf} \ar[r]^-{\times\RR} \ar[d]^-{\textup{incl$\,\circ\,$forget}} &
{\bA^h(X\times\RR^{i+1},n+1\cc)\hf} \ar[d]^-{\textup{forget}} \\
{\cone~\bA(X\times\RR^i\cc)} \ar[r]^-{\times\RR} & {\bA^h(X\times\RR^{i+1}\cc)}
}
\]
where lemma~\ref{lem-degenerateprod} is used to extend
$\times\RR\co \bA(X\times\RR^i\cc)\to \bA(X\times\RR^{i+1}\cc)$ to the cone on $\bA(X\times\RR^i\cc)$.
\end{cor}

\proof By corollary~\ref{cor-ultimateAstab}, this amounts to saying that
\[
\xymatrix{
\bE\hf \ar[r]^-{\textup{incl}} \ar[d]^-{\textup{incl$\,\circ\,$forget}} &
(S^1_!\wedge \bE)\hf \ar[d]^-{\textup{forget}} \\
S^1_\lambda\wedge \bE \ar[r]^-{\textup{incl}} & S^1\wedge \bE
}
\]
is homotopy cartesian, where $\bE$ is $\bA(X\times\RR^i,n\cc)$ and $S^1_\lambda$ is the lower half of $S^1$
(a closed interval with boundary equal to $S^0$). \qed

\bigskip
For the next lemma, let $\sR$ be any Waldhausen category with an
SW-product $\odot\lbul$ satisfying the axioms of \cite[\S2]{WWduality}.
We define $\sR^{(1)}$ as above, so $\sR^{(1)}$ is $\sR$ with a new SW-product which is a shift of the
old one. To describe $K$-theoretic and $L$-theoretic relationships between $\sR$ and $\sR^{(1)}$,
we introduce the Waldhausen category $\sP\sR$ of pairs in $\sR$: an object is a
cofibration $C_0\to C_1$ in $\sR$. Compare \cite[\S1]{WWduality}.
There is a canonical SW product on $\sP\sR$
making the inclusion and forgetful functors
\[  \sR^{(1)} \lra \sP\sR \lra \sR \]
given by $C\mapsto (*\to C)$ and $(C_0\to C_1)\mapsto C_0$
duality-preserving. The resulting diagram of spectra with $\ZZ/2$-action
\begin{equation} \label{eqn-notoriousfib}
\bK(\sR^{(1)}) \lra \bK(\sP\sR) \lra \bK(\sR)
\end{equation}
is a homotopy fiber sequence of spectra by the additivity theorem. It admits a non-equivariant homotopy
splitting $u\co \bK(\sR)\to \bK(\sP\sR)$ induced by the exact functor $C\mapsto(C,C)$; by the aforementioned
notorious fact from \cite{WWduality} and \cite{WWprodua}, the map $u\vee tu$ from
$\bK(\sR)\vee\bK(\sR)$ to $\bK(\sP\sR)$ is a homotopy equivalence (and a $\ZZ/2$-map).
It follows that in the homotopy fiber sequence
\begin{equation} \label{eqn-notoriousfibthf}
\bK(\sR^{(1)})\thf \lra \bK(\sP\sR)\thf \lra \bK(\sR)\thf
\end{equation}
determined by~(\ref{eqn-notoriousfib}), the middle term is contractible.
For much more elementary reasons,
\begin{equation} \label{eqn-notoriousfibthfL}
\bL\lbul(\sR^{(1)}) \lra \bL\lbul(\sP\sR) \lra \bL\lbul(\sR)
\end{equation}
is also a homotopy fiber sequence with contractible middle term. Putting those observations together,
we have a commutative diagram
\begin{equation} \label{eqn-notoriousdiagram}
\begin{split}
\xymatrix{
\bL\lbul(\sR^{(1)}) \ar[r] \ar[d]^-{\Xi} & \bL\lbul(\sP\sR)  \ar[r] \ar[d]^-{\Xi} & \ar[d]^-{\Xi}  \bL\lbul(\sR) \\
\bK(\sR^{(1)})\thf \ar[r] &  \bK(\sP\sR)\thf \ar[r] & \bK(\sR)\thf
}
\end{split}
\end{equation}
where the rows are homotopy fiber sequences and the terms in the middle column are contractible.
This implies a relationship between left-hand column
and right-hand column which may be loosely described as the \emph{shift invariance} of $\Xi$.

\medskip For a basic application of this, let $H$ be a
subgroup of the group $K_0(\sR)$ which is invariant under the involution
determined by the $SW$-product. Let $\sR^H\subset \sR$ be the full Waldhausen
subcategory consisting of all objects $C$ whose class $[C]\in K_0(\sR)$ belongs to
$H$. The $SW$-product on $\sR$ restricts to one on $\sR^H$ which also satisfies
the axioms of \cite[\S2]{WWduality}.

\begin{lem}
\label{lem-rothenberg} The following square is homotopy cartesian:
\[
\CD
\bL\lbul(\sR^H)   @>\Xi>>  \bA(\sR^H)\thf \\
@VV\textup{ inclusion }V @VV\textup{ inclusion }V \\
\bL\lbul(\sR)   @>\Xi>>  \bA(\sR)\thf.
\endCD
\]
\end{lem}

\proof It is well known \cite{Shanesonproduct} that the relative $n$-th homotopy
group of the inclusion map $\bL\lbul(\sR^H)\to \bL\lbul(\sR)$ is
isomorphic to the Tate cohomology group
\[ \hat H^{-n}(\ZZ/2;K_0/H)\]
where $K_0=K_0(\sR)$.
More precisely, let an element in this relative homotopy group be
represented by a pair $C\to D$ in $\sR$ with a nondegenerate
quadratic structure of formal dimension $n$, and with $C$ in
$\sR^H$. Poincar\'e duality implies that $[C]\in K_0(\sR)/H$ is
in the $(-1)^n$ eigensubgroup of the standard involution. Hence $[C]$
determines an element in the above-mentioned Tate cohomology group. This
describes the isomorphism. --- But now it is also clear that the relative $n$-th homotopy
group of the inclusion map $\bA(\sR^H)\thf\to \bA(\sR)\thf$
is isomorphic to
\[ \hat H^{-n}(\ZZ/2;K_0/H). \]
Indeed this relative homotopy group can be identified with $\pi_n$
of the spectrum $(\bA(\sR)/\bA(\sR^H))\thf$, and here $\bA(\sR)/\bA(\sR^H)$
is an Eilenberg-MacLane spectrum with (at most) one nonzero homotopy
group,
\[ \pi_0(\bA(\sR)/\bA(\sR^H)) = K_0/H\,. \]
Using these identifications, the relative $n$--th homotopy group for
the two columns in the above commutative square are already abstractly identified and
we want to show that $\Xi$ induces that identification of the $\pi_n$ groups. This is clear in
the case $n=0$ by inspection. The case of arbitrary $n$ follows from the case $n=0$
by the shift invariance of $\Xi$,
described in diagram~(\ref{eqn-notoriousdiagram}). \qed

\begin{prop}
\label{prop-timesRandXi} 
The following commutative squares are homotopy cartesian:
\[
\CD
\Omega^n\bL\lbul(X\times\RR^i\cc) @>\Xi >> \bA(X\times\RR^i,n\cc)\thf \\
@VV{\times\RR}V @VV{\times\RR}V \\
\Omega^{n+1}\bL\lbul(X\times\RR^{i+1}\cc) @>\Xi >> \bA(X\times\RR^{i+1},n+1\cc)\thf~,
\endCD
\]
\[
\CD
\Omega^n\bV\bL\ubul(X\times\RR^i\cc) @>\Xi >> \bA(X\times\RR^i,n\cc)\thf \\
@VV{\times\RR}V @VV{\times\RR}V \\
\Omega^{n+1}\bV\bL\ubul(X\times\RR^{i+1},n+1\cc) @>\Xi >> \bA(X\times\RR^{i+1},n+1\cc)\thf~.
\endCD
\]
\end{prop}

\medskip
\proof[Interpretation and proof] Let $\sR=\sR\ld(X\times\RR^i\cc)$ and let $\sR^{(n)}$ be the same with
a shifted SW-product. Where we have $\Omega^n\bL\lbul(X\times\RR^i\cc)=\Omega^n\bL\ubul(\sR)$ etc.~in the statement above,
we really mean
\[ \bL\lbul(\sR^{(n)}) \]
and justify that as in remark~\ref{rem-apologies}.  \newline
We prove the first statement, the other being similar. The diagram can be enlarged to
\begin{equation} \label{eqn-clevercut}
\CD
\Omega^n\bL\lbul(X\times\RR^i\cc) @>\Xi >> \bA(X\times\RR^i,n\cc)\thf \\
@VV{\times\RR}V @VV{\times\RR}V \\
\Omega^{n+1}\bL\lbul^h(X\times\RR^{i+1}\cc) @>\Xi >> \bA^h(X\times\RR^{i+1},n+1\cc)\thf \\
@VV\textup{incl.}V @VV\textup{incl.}V \\
\Omega^{n+1}\bL\lbul(X\times\RR^{i+1}\cc) @>\Xi >> \bA(X\times\RR^{i+1},n+1\cc)\thf~.
\endCD
\end{equation}
In diagram~(\ref{eqn-clevercut}), the lower square is homotopy cartesian by lemma~\ref{lem-rothenberg},
and in the upper square the left-hand column is a homotopy equivalence by corollary~\ref{cor-stabLandVL}.
It is therefore enough to show that the right-hand column
in the upper square of~(\ref{eqn-clevercut}) is also a homotopy equivalence. This follows easily from
corollary~\ref{cor-ultimateAstab}.
\qed

\bigskip
We turn to a slightly different but related theme: the homotopy
invariance properties of constructions such as $\bL\lbul$,
$\bV\bL\ubul$ and $\bA$ when applied to control spaces such as
$(X*S^{i-1},X\times\RR^i)$ where $X$ is compact Hausdorff. It is
easy to show homotopy invariance in the
variable $X$, but we are also interested in the other variable,
the sphere.

\begin{lem} \label{lem-homotopyinvariance} For a compact Hausdorff space $X$, and $I=[0,1]$, the projection
\[ (X*S^{i-1})\times I \to X*S^{i-1} \]
induces homotopy equivalences
\[ \begin{array}{ccc}
\bA^h(((X*S^{i-1})\times I,X\times\RR^i\times I)) &\lra & \bA^h((X*S^{i-1},X\times\RR^i))  \\
\bL\lbul^h(((X*S^{i-1})\times I,X\times\RR^i\times I)) &\lra & \bL\lbul^h((X*S^{i-1},X\times\RR^i))  \\
{\bVL\ubul}^h(((X*S^{i-1})\times I,X\times\RR^i\times I)) &\lra & {\bVL\ubul}^h((X*S^{i-1},X\times\RR^i)).
\end{array}
\]
\end{lem}

\begin{cor}
\label{cor-homotopyinvariance}
The constructions $\bL\lbul$, $\bV\bL\ubul$ and $\bA$ applied to
control spaces of the form $(X*S^{i-1},X\times\RR^i)$, with compact
Hausdorff $X$, are homotopy invariant as functors of the sphere variable.
\end{cor}

\proof[Proof of the corollary modulo the lemma]
What we mean here by \emph{homotopy invariant} is that
homotopic maps $f,g$ from $S^{i-1}$ to $S^{j-1}$ induce homotopic maps
\[ f_*,g_*\co \bA((X*S^{i-1},X\times\RR^i)) \lra
\bA((X*S^{j-1},X\times\RR^j)) \]
etc.~To prove it, we use theorem~\ref{thm-changeofdeco}; this increases $i$ and $j$ by 1, and replaces
$\bA$ by $\bA^h$. Now we can apply the lemma. \qed

\medskip
\proof[Proof of lemma~\ref{lem-homotopyinvariance}]
We concentrate on the $\bA$-theory case (and return to the explicit notation
where control spaces are shown with their singular parts). \newline
Now we have several choices of Waldhausen category
with an algebraic $K$-theory spectrum that deserves to be called
$\bA^h((X*S^{i-1},X\times\RR^i))$. Among these we choose one which
is more ``algebraic'' than the one which we normally prefer. This is described
in \cite[\S6]{Weissexci}. The objects of this Waldhausen category, call it
$\sR$ for now, are certain retractive $CW$-spaces $Y$ over $X$ with a
finite dimensional
relative $CW$-structure and a map from the set of
cells to $\RR^i$. This labelling map must satisfy certain local finiteness
and control conditions. (The control conditions are expressed in terms
of the compactification $D^i$ of $\RR^i$.)
The morphisms in $\sR$ are retractive cellular maps (over and under $X$)
which, again, satisfy certain control conditions (expressed in terms of the
cell labels).
The relationship between this Waldhausen category $\sR$
and the usual one, $\sR\lf((X*S^{i-1},X\times\RR^i))$~, is essentially given
by \emph{cobase change along the projection $X\times\RR^i\to X$},
regarded as a functor
\[ \sR\lf((X*S^{i-1},X\times\RR^i)) \to \sR. \]
Note that the objects in $\sR\lf((X*S^{i-1},X\times\RR^i))$ are retractive
spaces over $X\times\RR^i$ with a relative $CW$-structure.
See \cite[\S9]{Weissexci} for related ideas. \newline
We introduce a similar
Waldhausen category $\sQ$ with an algebraic $K$-theory spectrum that
deserves to be called $\bA^h(((X*S^{i-1})\times I,X\times\RR^i\times I))$. Its
objects are retractive spaces over $X$ with a relative $CW$-structure
and a map from the set of cells to $\RR^i\times I$, subject to
certain local finiteness and control conditions (which are formulated
using the compactification
$D^i\times I$ of $\RR^i\times I$). It is convenient to stabilize
the two categories by introducing formal desuspensions; then we have
$s\sQ$ and $s\sR$. Now our task is to show that the functor
\[ s\sQ\to s\sR \]
induced by the projection $\RR^i\times I\to \RR^i$ induces a homotopy
equivalence of the algebraic $K$-theory spectra. Note that this functor
does not satisfy the first hypothesis of the approximation theorem
(which is about ``detection'' of weak equivalences).
We use the fibration theorem instead.
In addition to the standard notion of weak equivalence in $s\sQ$, we
therefore introduce a coarser notion. Let $I'=\,]0,1]$ and call a morphism
in $s\sQ$ a $v$-equivalence
if the germ near $S^{i-1}\times I'$ of its mapping cone is weakly
equivalent to zero in the usual controlled sense. Note here that $S^{i-1}\times I'$
is part of the singular set $S^{i-1}\times I$ of the control space
$((X*S^{i-1})\times I,X\times\RR^i\times I))$. By the fibration theorem there is
a homotopy fibration sequence of $K$-theory spectra,
\[ \bK(s\sQ^v)\to \bK(s\sQ) \to \bK(s\sQ_v)\,.\]
Hence it is enough to verify that
\begin{itemize}
\item[(i)] $s\sQ_v$ is ``flasque'' enough so that $\bK(s\sQ_v)$ is contractible~;
\item[(ii)] the composition of exact functors $s\sQ^v\to s\sQ\to s\sR$ satisfies
the hypotheses of the approximation theorem, so that $\bK(s\sQ^v)\simeq \bK(s\sR)$.
\end{itemize}
As regards (i), the situation is really quite similar to the one
which we have considered in the ``homotopy invariance'' part of the
proof of theorem~\ref{thm-homolothy}. To make this connection
clearer we note that $s\sQ_v$ can be replaced by a simpler category
$s\sQ_{(v)}$ with the same objects, where the morphisms
themselves are germs of retractive maps over $X\times\RR^i\times I$ defined
near (i.e., over an arbitrarily small neighborhood of) $S^{i-1}\times I'$.
Indeed, another application of the approximation theorem shows that
the forgetful functor $s\sQ_v\to s\sQ_{(v)}$ induces a
homotopy equivalence of the $K$-theory spectra. Up to an equivalence
of categories, $s\sQ_{(v)}$ however depends only on
\[   (W,W\smin(S^{i-1}\times I')) \]
where $W$ can be any neighborhood of $S^{i-1}\times I'$ in $(X*S^{i-1})\times
I$. We can take $W$ to be homeomorphic to the mapping cylinder of
the projection
\[ X\times S^{i-1}\times I' \to S^{i-1}\times I'. \]
In particular when $X$ is a point, the control pair $(W,W\smin(S^{i-1}\times I'))$
is what we have called $\JJ(S^{i-1}\times I')$. The same arguments for
``flasqueness'' as in the homotopy invariance part of the proof of
theorem~\ref{thm-homolothy} apply to the category $s\sQ_{(v)}$~; this includes
remark~\ref{rem-corr} below. \newline
As regards (ii), it is hard to verify directly that the hypotheses
of the approximation theorem are satisfied by the functor
$s\sQ^v\to s\sR$. It seems wiser therefore to introduce a full
Waldhausen subcategory $s\sQ^{(v)}$ of $s\sQ^v$ consisting of the
objects $Y$ whose set of cell labels avoids a neighborhood
(depending on $Y$) of $S^{i-1}\times I'$ in $(X*S^{i-1})\times I$. Then it is
clear that the composition
\[ s\sQ^{(v)} \to s\sQ^v \to s\sQ \to s\sR \]
is an equivalence of (Waldhausen) categories. Hence it remains only
to verify that the inclusion $s\sQ^{(v)} \to s\sQ^v$ satisfies the
hypotheses of the approximation theorem. The first condition
(about detection of weak equivalences) is clearly satisfied. For the
second condition, suppose given a cofibration
\[  f\co Y_1\to Y_2 \]
in $s\sQ^v$~, with $Y_1$ in $s\sQ^{(v)}$. From the definitions, it
is easy to construct a factorization
\[ Y_1\to Y_3\to Y_2 \]
of $f$, where $Y_1\to Y_3$ is still a cofibration, $Y_3$ is also
in $s\sQ^{(v)}$ and the morphism from $Y_3$ to $Y_2$ is a ``domination''
(i.e. there exists a controlled map $Y_2\to Y_3$, not necessarily a
morphism, which is right inverse to $Y_3\to Y_2$ up to a controlled
homotopy). Attaching additional cells to $Y_3$ where necessary, one can
then easily improve $Y_3\to Y_2$ to a weak equivalence in $s\sQ^v$. (Solve the
corresponding problem in $s\sR$ first and use that solution as a model.) \qed

\section{Spherical fibrations and twisted duality}
\label{sec-twisted} It is well known \cite{Vogellduality} that a
spherical fibration $\xi$ on a space $X$ determines a twisted $SW$
product in the category of finitely dominated retractive spaces on
$X$. This is compatible with the so--called $w$--twisted involution
on $\ZZ[\pi_1X]$ where $w\co \pi_1(X)\to \ZZ$ is $w_1(\xi)$. We
recall some of the details, following \cite[1.A.9]{WWduality} rather
more than \cite{Vogellduality}. \newline We assume that $\xi\co E\to
X$ is a fibration with fibers homotopy equivalent to $S^d$ and with
a preferred section which is a fiberwise cofibration \cite{James89}.
(What we have in mind is, for example, the fiberwise one--point
compactification of a vector bundle on $X$, with the preferred
section which picks out the point at infinity in each fiber.) Let
$Y_1$ and $Y_2$ be finitely dominated retractive spaces over $X$,
with retractions $r_1$ and $r_2$. Again we start by defining an
unstable Spanier--Whitehead product $Y_1\curlywedge Y_2$. This is
the based space obtained by first forming the external smash product
\[
\begin{array}{ccc}
Y_1 \,\,{}_X\!\wedge_X\, Y_2 & = & Y_1\times Y_2/\sim
\end{array}
\]
where $\sim$ identifies $(y_1,x)$ with $(r_1(y_1),x)$ and $(x,y_2)$
with $(x,r_2(y_2))$~; then forming the homotopy pullback of
\[
\CD   E @>\textup{\,\,diagonal$\,\circ\,$proj.\,\,}>> X\times
X@<\textup{\,\,retraction\,\,}<< Y_1 \,\,{}_X\!\wedge_X\, Y_2~;
\endCD
\]
then collapsing the subspace consisting of all elements in the
homotopy pullback which are mapped to the base point under the
forgetful projection to
\[ (E/X)\wedge (Y_1/X)\wedge(Y_2/X)\,. \]
To make this unstable $SW$ product into an $SW$ product on the
stable category $s\sR(X)$ of finitely dominated retractive spaces
over $X$, we proceed much as in \S2.

\begin{defn}
\label{defn-twistedSW} {\rm We let
\[
(Y_1,k)\odot(Y_2,\ell) = \colim_n
\,\Omega^{2n+d}(\Sigma^{n-k}Y_1\curlywedge\Sigma^{n-\ell}Y_2).
\]
More generally we let
\[
(Y_1,k)\odot_j(Y_2,\ell) = \colim_n
\,\Omega^{2n+d}\Sigma^j(\Sigma^{n-k}Y_1\curlywedge\Sigma^{n-\ell}Y_2),
\]
so that $(Y_1,k)\odot(Y_2,\ell)=(Y_1,k)\odot_0(Y_2,\ell)$, and
denote the $\Omega$--spectrum with $j$--th term
$(Y_1,k)\odot_j(Y_2,\ell)$ by $(Y_1,k)\odot\lbul(Y_2,\ell)$. }
\end{defn}

\medskip
This comes with a structural symmetry $(Y_1,k)\odot\lbul(Y_2,\ell)
\cong (Y_2,\ell)\odot\lbul(Y_1,k)$ determined by the obvious
symmetry of $\curlywedge$. For $Y_1=Y_2=Y$ and $k=\ell$ we obtain an
$\Omega$--spectrum $(Y,k)\odot\lbul(Y,k)$ with an action of $\ZZ/2$.
\newline Note the $d$ in definition~\ref{defn-twistedSW} which is
the formal fiber dimension of $\xi$. This causes a slight
disagreement with the conventions of \cite[1.A.9]{WWduality}, but it
is convenient here.

\medskip
Definition~\ref{defn-twistedSW} is essentially insensitive to a
stabilization of $\xi$ by fiberwise suspension. More precisely, a
spherical fibration $\xi$ on $X$ as above determines an $SW$ product
$\odot_{\lbul}$ on $s\sR(X)$, and the fiberwise suspension
$\Sigma\xi$ determines another, which we (temporarily) denote by
$\odot'_{\lbul}$. There is a natural homotopy equivalence
\[ (Y_1,k)\odot\lbul(Y_2,\ell)\lra (Y_1,k)\odot'\lbul(Y_2,\ell) \]
which respects the canonical involutions. Note also that the
$\xi$--twisted version of $(Y_1,k)\odot\lbul(Y_2,\ell)$ is naturally
homeomorphic to the standard one (definition~\ref{defn-SW}) if $\xi$
is a trivial sphere bundle $S^0\times X\to X$.

\medskip
Definition~\ref{defn-nonlinearsymandvis} can be re--used with the
$\xi$--twisted interpretation of the $SW$ product. The
$\xi$--twisted versions of corollary~\ref{cor-tDalmostsplitting} and
definition~\ref{defn-nonlinearvishyper} take a slightly more
complicated form. For a finitely dominated retractive space $Y$ over
$X$, let $Y^{\xi}$ be the fiberwise smash product of $Y$ and
$E=E(\xi)$ over $X$. As before, suppose that the formal fiber
dimension of $\xi$ is $d$.

\begin{cor}
\label{cor-twistedtDalmostsplitting} There is a natural homotopy
fiber sequence of spectra
\[
\CD ((Y,k)\odot\lbul(Y,k))_{h\ZZ/2} @>>>
((Y,k)\odot\lbul(Y,k))^{\ZZ/2} @>J>>
\Sigma^{\infty-k-d}(Y^{\xi}/X)\,.
\endCD \]
\end{cor}

\begin{defn}
\label{defn-twistednonlinearvishyper} {\rm An $n$--dimensional
visible hyperquadratic structure on $(Y,k)$ is an element in
$\Omega^n\Omega^{\infty}\Sigma^{\infty-k-d}(Y^{\xi}/X)$. An
$n$--dimensional quadratic structure on $(Y,k)$ is an element of
$\Omega^n\Omega^{\infty}(((Y,k)\odot\lbul(Y,k))_{h\ZZ/2})$.
Alternatively, an $n$--dimensional quadratic structure on $(Y,k)$
can be defined as an element of $\Omega^n\Omega^{\infty}$ of the
homotopy fiber of the natural map $J\co
((Y,k)\odot\lbul(Y,k))^{\ZZ/2} \to \Sigma^{\infty-k-d}(Y^{\xi}/X)$.
}
\end{defn}

\medskip
We write $\bL\lbul(X,\xi)$, $\bV\bL\ubul(X,\xi)$ and $\bL\ubul(X,\xi)$
for the $L$--theory spectra determined by the $\xi$--twisted $SW$
product on the stable category of finitely dominated retractive
spaces over $X$. Theorem~\ref{thm-VhatLexcision} remains correct
(with essentially the same proof) in this more general setting and
can informally be regarded as a statement for spaces $X$ over $BG$,
the classifying space for stable spherical fibrations. A \emph{weak
homotopy equivalence}, in that setting, is a map of spaces over $BG$
which is a weak homotopy equivalence of spaces. A \emph{cocartesian
square}, in that setting, is a commutative square of spaces over
$BG$ which is cocartesian as a square of spaces.

\bigskip
We come to an outline of a calculation of $\bV\widehat\bL\ubul(X,\xi)$
(which will not be used elsewhere in the paper). Let $\Th(X,\xi)$ be
the Thom spectrum of $X$ and $\xi$. For convenience or otherwise, we
index that in such a way that the $(d+k)$--th space is
$S^k\wedge(E/X)$. There is a natural map
\[ \hat\sigma\co \Th(X,\xi) \lra \bV\widehat\bL\ubul(X,\xi)~, \]
the \emph{visible hyperquadratic signature map}. \newline
In more detail,
$\pi_n\Th(X,\xi)$ can be identified with the group of bordism
classes of \emph{normal spaces} of formal dimension $n$ over
$(X,\xi)$. A normal space of formal dimension $n$ over $(X,\xi)$
consists of a finitely dominated space $Y$, a map $f\co Y\to X$ and
a stable (pointed) map $\eta\co S^{n+d}\to E(f^*\xi)/Y$. (The image
under the Thom isomorphism of the homology class carried by $[\eta]$
is a class in $H_n(Y;\ZZ^w)$, where $\ZZ^w$ is the twisted integer
coefficient system determined by $w_1(f^*\xi)$. It is called the
fundamental class of the normal space, but it is not subject to any
Poincar\ee duality condition. If it \emph{does} satisfy the
Poincar\ee duality condition as formulated by Wall \cite{Wall70},
then that makes $Y$ into a Poincar\ee duality space with Spivak
normal bundle $f^*\xi$.) Therefore normal spaces generalize
Poincar\ee duality spaces, and in fact the visible hyperquadratic signature
map is a variant of a better known and easier--to--understand map
$\sigma$ from the Poincar\ee duality bordism spectrum
$\Bo_{\textup{PD}}(X,\xi)$ of $(X,\xi)$ to the visible symmetric
$L$--theory spectrum $\bV\bL\ubul(X,\xi)$. The two maps fit into a
commutative square
\[
\CD \Bo_{\textup{PD}}(X,\xi) @>\sigma>> \bV\bL\ubul(X,\xi) \\
@VV\textup{forgetful} V @VV\textup{inclusion}V \\
\Th(X,\xi) @>\hat\sigma>> \bV\widehat\bL\ubul(X,\xi).
\endCD
\]
\newline
In both the normal bordism theory and the visible hyperquadratic
theory, there are external products. They have the form
\[
\begin{array}{c}
\Th(X,\xi)\wedge \Th(X',\xi') \lra \Th(X\times X',\xi\times \xi')~, \\
\rule{0mm}{4mm}\bV\widehat\bL\ubul(X,\xi)\wedge
\bV\widehat\bL\ubul(X',\xi') \lra \bV\widehat\bL\ubul(X\times
X',\xi\times\xi')
\end{array}
\]
where $\xi\times\xi'$ is the external fiberwise smash product of
$\xi$ and $\xi'$. The composition
\[
\CD \Th(X,\xi)\wedge \bV\widehat\bL\ubul(\pt) @>\sigma\wedge\,\id >>
\bV\widehat\bL\ubul(X,\xi)\wedge \bV\widehat\bL\ubul(\pt) @>\textup{ext.
prod.}>> \bV\widehat\bL\ubul(X,\xi)
\endCD
\]
is a natural transformation between excisive functors on the
category of spaces over $BG$. It is an equivalence when $X$ is a
point (mapping to the base point of $BG$). Hence it is always an
equivalence and we have

\begin{thm}
\label{thm-twistedVhatLsymofpoint} $\,\,\bV\widehat\bL\ubul(X,\xi)
\,\,\simeq\,\, \Th(X,\xi) \,\wedge\,(\bS\vee\RR\bP^{\infty}_{-1})$.
\end{thm}

Moreover, the fact that $\hat\sigma$ commutes with external products
immediately leads to a ``calculation'' of $\hat\sigma$:

\begin{prop}
\label{prop-whatissigmahat} The hyperquadratic signature map
\[ \hat\sigma\co \Th(X,\xi)\lra \bV\widehat\bL\ubul(X,\xi) \,\,\simeq\,\,
\Th(X,\xi) \,\wedge\,(\bS\vee\RR\bP^{\infty}_{-1}) \] is homotopic
to the inclusion of the wedge summand $\Th(X,\xi) \,\wedge\,\bS
\simeq \Th(X,\xi)$.
\end{prop}

\bigskip
Next we need a $\xi$--twisted version of section~\ref{sec-excision}.
We begin with a control space $(\bar Q,Q)$ as
in~\ref{defn-controlledSW}, and a spherical fibration $\xi\co E\to
Q$ of formal fiber dimension $d$, with a distinguished section.
\newline Let $Y$ and $Z$ be objects of $\sR\ld(\bar Q,Q)$. To define
their $SW$ product $Y\odot Z$, we introduce first an unstable form
$Y\curlywedge Z$ of it. We define it as the geometric realization of
a based simplicial set. An $n$--simplex of this simplicial set is a
pair $(f,\gamma)$ where
\begin{itemize}
\item[(i)] $f$ is a continuous map from
the standard $n$--simplex $\Delta^n$ to the topological inverse
limit of the spaces $(Y/Q)\wedge (E^P/Q)\wedge (Z^P/Q)$, where $P$ runs
through the large closed subsets of $Q$~;
\item[(ii)] $\gamma$ is a continuous assignment $c\mapsto \gamma_c$
of paths in $Q$, defined for $c\in \Delta^n$ with $f(c)$ not equal
to the base point $\pt\,$.
\end{itemize}
The paths $\gamma_c$ are to be parametrized by $[-1,+1]$ and must
satisfy
\[ \gamma_c(-1)=r_Yf_Y(c)~, \qquad
\gamma_c(+1)=r_Zf_Z(c)~, \qquad \gamma_c(0)= \xi f_E(c)\,. \] There
is a control condition:
\begin{itemize}
\item[] For $q\in \bar Q\smin Q$
and any neighborhood $V$ of $q$ in $\bar Q$, there exists a smaller
neighborhood $W$ of $q$ in $\bar Q$ such that, for any $c\in
\Delta^n$ with $f(c)\ne \pt$, the path $\gamma_c$ either avoids $W$
or runs entirely in $V$.
\end{itemize}

\medskip
\begin{defn}
\label{defn-twistedcontrolledSW} {\rm For $Y$ and $Z$ in
$\sR\ld(\bar Q,Q)$ and integers $k,\ell\in\ZZ$, let
\[
(Y,k)\odot(Z,\ell) = \colim_n
\,\Omega^{2n+d}(\Sigma^{n-k}Y\curlywedge \Sigma^{n-\ell}Z)\,.
\]
More generally let $(Y,k)\odot\lbul(Z,\ell)$ be the spectrum with
$j$--th space
\[(Y,k)\odot_j(Z,\ell) =
\colim_n \,\Omega^{2n+d}\Sigma^j(\Sigma^{n-k}Y\curlywedge
\Sigma^{n-\ell}Z). \] }
\end{defn}

\medskip
This is the appropriate definition for the category $s\sR\ld(\bar
Q,Q)$. There is also a germ version, for the category
$s\sR\sG\ld(\bar Q,Q)$, which generalizes
definition~\ref{defn-germSWproducts}. With that generalized
definition, theorem~\ref{thm-homolothy} generalizes as follows:

\begin{thm}
\label{thm-twistedhomolothy} The spectrum valued functor $X\mapsto
\bE(X)$ is homotopy invariant and excisive. Here $X$ can be a space
over $BG$ or a space with a spherical fibration $\xi$, and $\bE(X)$
means $\bL\ubul(\JJ X_{\infty})$, defined using the $\xi$--twisted
$SW$ product.
\end{thm}

We leave the detailed formulation to the reader. The proof is
essentially identical with the proof of theorem~\ref{thm-homolothy}.
\newline Furthermore definition~\ref{defn-controllednonlinearvis}
can be re--used in the twisted setting and leads to a generalization
of theorem~\ref{thm-visiblehomolothy}:

\begin{thm}
\label{thm-twistedvisiblehomolothy} The spectrum valued functor
$X\mapsto \bE(X)$ is homotopy invariant and excisive. Here $X$ can
be a space over $BG$ or a space with a spherical fibration $\xi$,
and $\bE(X)$ means $\bV\bL\ubul(\JJ X_{\infty})$, defined using the
$\xi$--twisted $SW$ product.
\end{thm}

\medskip The results of sections~\ref{sec-controlvisible}
and~\ref{sec-controlsusp} also have
generalizations to the twisted setting. We do not formulate them here
explicitly. They will however be used in sections~\ref{sec-approxdetails}
and~\ref{sec-approxproofs}.

\medskip
\section{Homotopy invariant characteristics and signatures} \label{sec-hocharsig}
This section is formally analogous to parts of
\cite[\S6]{DwyerWeissWilliams}. We begin with a space $X$, a
spherical fibration $\xi$ on $X$ (with a distinguished section which
is a fiberwise cofibration), and an integer $n\ge 0$. From these
data we produce a spectrum
\[ \bV\bLA\ubul(X,\xi,n). \]
In the case where $X$ is a finitely dominated Poincar\ee duality
space of formal dimension $n$ and $\xi$ is the Spivak normal
fibration of $X$, we also construct a characteristic element
\[  \sigma(X)\in F(X,\xi,n):= \Omega^{\infty+n}\bV\bLA\ubul(X,\xi,n)\,. \]
This refines the Mishchenko--Ranicki (visible) symmetric signature
of $X$ which, in the nonlinear setting, is an element of
$\Omega^{\infty+n}\bV\bL\ubul(X,\xi)$ or of
$\Omega^{\infty+n}\bL\ubul(X,\xi)$. The construction has certain
naturality properties. As in \cite[\S1]{DwyerWeissWilliams}, these
properties imply that every \emph{family} of finitely dominated
formally $n$--dimensional Poincar\ee duality spaces $X_b$~,
depending on a parameter $b\in B$, determines a characteristic
section $\sigma(p)$ of a fibration on $B$ whose fibers are,
essentially, the spaces $F(X_b,\xi_b,n)$.

\medskip
Let $\sR$ be any Waldhausen category with stable SW-product satisfying the axioms of \cite[\S2]{WWduality}.
We also need the shifted cousins of this SW-product and write
$\sR^{(n)}$
for $\sR$ with the shifted SW-product, as in section~\ref{sec-controlsusp}.
The SW-product determines a duality involution on $\bK(\sR^{(n)})$ as explained
in \cite{WWduality}. (See remark~\ref{rem-corrduality} below for a
correction.) Using that, we have the map
\begin{equation} \Xi\co \bL\ubul(\sR^{(n)}) \lra \bK(sR^{(n)})\thf
\end{equation}
of \cite{WWduality}. We also write this in the form
\begin{equation} \label{eqn-happyXi}
\Xi\co \Omega^n\bL\ubul(\sR) \lra \bK(\sR^{(n)})\thf
\end{equation}
noting that $\bL\ubul(\sR^{(n)})\simeq \Omega^n\bL\ubul(\sR)$, in the spirit of remark~\ref{rem-apologies}. \newline
Now suppose that $\sR$ is $s\sR(X)$ with the $\xi$--twisted stable SW--product
introduced in section~\ref{sec-twisted}. Then we have good reasons to write
\begin{equation}
\Xi\co \Omega^n\bL\ubul(X,\xi) \lra \bA(X,\xi,n)\thf
\end{equation}
for the map~(\ref{eqn-happyXi}). We can also restrict to
$\Omega^n\bVL\ubul(X,\xi)$ to get
\begin{equation} \label{eqn-geoXi}
\Xi\co \Omega^n\bVL\ubul(X,\xi) \lra \bA(X,\xi,n)\thf.
\end{equation}
This brings us to the definition of $\bVLA\ubul$ anticipated in section~\ref{sec-intro}:

\begin{defn} \label{defn-LK} 
{\rm Let $\Omega^n\bVLA\ubul(X,\xi,n)$ be the homotopy pullback of
\[
\xymatrix@C=50pt{
\Omega^n\bVL\ubul(X,\xi) \ar[r]^{\Xi}_{(\ref{eqn-geoXi})} & \bA(X,\xi,n)\thf & \ar[l]_-{\textup{incl.}} \bA(X,\xi,n)\hf~.
}
\]
}
\end{defn}

\bigskip
Now let $X$ be a finitely dominated Poincar\ee duality space of
formal dimension $n$, with Spivak normal bundle $\xi\co E\to X$ and
a choice of a stable map $\eta\co S^{n+d}\to E/X$ such that $[\eta]$
maps to a fundamental class for $X$ under the Thom isomorphism. (We
assume as usual that $\xi$ has fibers $\simeq S^d$ and comes with a
distinguished section. Together, $\xi$ and $\eta$ are determined by
$X$ up to contractible choice and $\xi$ can therefore be
regarded as the Spivak normal fibration of $X$.) \newline Let
$s\sR=s\sR(X)$ with the $\xi$--twisted $SW$ product. Let
\[ F(X,\xi,n) = \Omega^{\infty+n}\bVLA\ubul(X,\xi,n) \]
which can also be described as the homotopy pullback of the diagram of infinite
loop spaces
\begin{equation} \label{eqn-F(X,xi,n)}
\xymatrix@M=4pt@C=25pt{ \Omega^{\infty}\bVL\ubul(s\sR^{(n)}) \ar[r]^-{\Xi} &
\Omega^\infty(\bK(s\sR^{(n)})\thf) & \ar[l]_-{\textup{incl.}}
\Omega^\infty(\bK(s\sR^{(n)})\hf)\,.
}
\end{equation}
To construct $\sigma(X)\in F(X,\xi,n)$, we begin as usual with
$S^0\times X$, viewed as a retractive space over $X$ alias object of
$s\sR$. The composite stable map
\[
\CD S^{n+d} @>>> E/X @>>> X_+\wedge (E/X)
\endCD
\]
defines an element $\varphi$ in $(S^0\times X)\odot_0(S^0\times X)$,
for the $SW$ product in $s\sR^{(n)}$. This is fixed under the
symmetry involution and nondegenerate, and so $S^0\times X$ and
$\varphi$ together determine a vertex, say $v_L$~, in the standard
simplicial set model for $\Omega^\infty\bVL\ubul(s\sR^{(n)})$.
They also determine a (homotopy) fixed point for the duality
involution on $\Omega^\infty\bK(s\sR^{(n)})$, which we may regard
as a vertex $v_K$ in $\Omega^\infty(\bK(s\sR^{(n)})\hf)$.
The images of $v_L$ and $v_K$ respectively under the two arrows in diagram~(\ref{eqn-F(X,xi,n)})
agree by the definition of $\Xi$, for which we refer again to \cite[\S9]{WWduality}.

\begin{defn}
\label{defn-homotopycharacteristic} {\rm Let $\sigma(X)\in
F(X,\xi,n)$ be the element determined by $v_L$, $v_K$ and the
constant path connecting their images under the two arrows in diagram~(\ref{eqn-F(X,xi,n)}).
}
\end{defn}

\medskip
It remains to be said how the assignment $\sigma$ is a
\emph{characteristic} \cite[1.1]{DwyerWeissWilliams} on a suitable
category of formally $n$--dimensional Poincar\ee duality spaces and
homotopy equivalences. This forces on us another revision of
$F(X,\xi,n)$, in fact an enlargement, which will as usual leave the
homotopy type unchanged. \newline We have some freedom in
interpreting the term $\Omega^{\infty}\bV\bL\ubul(s\sR^{(n)})$ from
the definition of $F(X,\xi,n)$. For a start we can regard it as the
0-th infinite loop space in the $\Omega$--spectrum
$\bV\bL\ubul(s\sR^{(n)})$. This was defined, following
\cite[\S9]{WWduality}, as the geometric realization of a
$\Delta$-set
\[  [m] \mapsto \textrm{vsp}_0(s\sR^{(n)}(m)). \]
Here $s\sR^{(n)}(m)$ is a category of certain functors from the
poset of faces of $\Delta^m$ to $s\sR^{(n)}$. The notation
$\textrm{vsp}_0(s\sR^{(n)}(m))$ means: the \emph{set} of visible
symmetric Poincar\ee objects of formal dimension zero in
$s\sR^{(n)}(m)$, that is, objects of $s\sR^{(n)}(m)$ with an
appropriate nondegenerate visible symmetric structure. Let
\[\textrm{vsp}'_0(s\sR^{(n)}(m)) \]
be the classifying space of the \emph{category} of visible symmetric
Poincar\ee objects of formal dimension zero in $s\sR^{(n)}(m)$. A
morphism between visible symmetric Poincar\ee objects of formal
dimension zero in $s\sR^{(n)}(m)$ is a morphism in $s\sR^{(n)}(m)$
respecting the visible symmetric structures; such a morphism is
automatically a weak equivalence in $s\sR^{(n)}(m)$. Then
\[  [m] \mapsto \textrm{vsp}'_0(s\sR^{(n)}(m))
\]
is a $\Delta$-space. Its geometric realization is an enlarged
version of the above construction of the 0-th infinite loop space in
the $\Omega$--spectrum $\bV\bL\ubul(s\sR^{(n)})$. It is homotopy
equivalent to the original, by a standard argument which exploits
the fact that several ways of realizing a bisimplicial set give the
same result. The map $\Xi$ extends easily to this enlarged version.
Consequently we end up with an enlarged version of $F(X,\xi,n)$.
Using this, it is easy to promote $\sigma(X)\in F(X,\xi,n)$ of
definition~\ref{defn-homotopycharacteristic} to a characteristic.
\newline Namely, let $X_i$ for $i=0,1,\dots,k$ be Poincar\ee duality
spaces of formal dimension $n$, with Spivak normal fibrations
$\xi_i\co E_i\to  X$ with fibers $\simeq S^d_i$ where
\[  d_k\le d_{k-1}\le \dots \le d_0 \]
and preferred zero sections which are fiberwise cofibrations. Let
$\eta_i\co S^{n+d_i}\to E_i/X_i$ be stable maps representing
fundamental classes for the $X_i$. Let homotopy equivalences $u_i\co
X_i \to X_{i-1}$ be given for $i=1,\dots,k$, covered by maps
\[  v_i \co \Sigma^{d_{i-1}-d_i}_{X_i}E_i \lra E_{i-1} \]
respecting the zero sections, and such that $v_i\eta_i=\eta_{i-1}$.
We can describe these data by a diagram
\[ (X_0,\xi_0,\eta_0)\leftarrow (X_1,\xi_1,\eta_1)\leftarrow \cdots\leftarrow
(X_{k-1},\xi_{k-1},\eta_{k-1})\leftarrow(X_k,\xi_k,\eta_k)~.\] With
the new definition of $F$, the diagram determines a map from the standard $k$-simplex
$\Delta^k$ to $F(X_0,\xi_0,n)$, which we could call the \emph{characteristic of
the diagram}. This assignment extends
definition~\ref{defn-homotopycharacteristic} and it has the
naturality properties which make it into a characteristic (on a
certain category $\sP_n$) as defined in
\cite[1.1]{DwyerWeissWilliams}. The objects of $\sP_n$ are triples
$(X,\xi,\eta)$ as above. It is very fortunate that we can allow
continous variation of $\eta$~; that is to say, the characteristic
$\sigma$ depends continuously on the ``reductions'' $\eta$.

\medskip
We spell out what \cite[1.6]{DwyerWeissWilliams} means here. Let
$p\co Y\to B$ be a fibration whose fibers $Y_b$ are finitely
dominated Poincar\ee duality spaces of formal dimension $n$. We
assume for simplicity that $B$ is a simplicial complex. We also need
a fiberwise Spivak normal fibration. Suppose that this comes in the
shape of a fibration $\xi\co E\to Y$ with preferred section, with
fibers $\simeq S^d$, and a map
\[ \eta\co B_+\wedge S^{n+d} \lra E/\sim \]
over $B$, where $E/\sim$ is the pushout of $E\leftarrow Y\to B$.
Every simplex $K\subset B$ determines a Poincar\ee duality space
$Y_K=p^{-1}(K)$ with Spivak normal fibration $\xi_K\co E_K\to Y_K$~,
where $E_K$ is the portion of $E$ above $K$, and a family of
reductions $\eta_b\co S^{n+d}\to E_K/Y_K$ where $b\in K$. These data
determine a map
\[  \sigma(Y_K)\co K \to F(Y_K,\xi_K,n) \]
using the continuity property of $\sigma$. As $K$ varies, we have
maps
\[  \hocolim_K \, K \lra \hocolim_K F(Y_K,\xi_K,n)\lra \hocolim_K \pt~. \]
The space on the right is the barycentric subdivision of $B$~; the
map on the right is a quasifibration which we can also describe as
$F_B(Y,\xi,n)\to B$~; the composite map from left to right is a
homotopy equivalence. Hence the map on the left can be viewed as a
``homotopy section'' $\sigma(p)$ of $F_B(Y,\xi,n)\to B$, and this is
what we wanted.

\bigskip
Next we need a generalization or adaptation of
definitions~\ref{defn-LK} and~\ref{defn-homotopycharacteristic} to
the controlled setting. Let $(\bar X,X)$ be a control space with
compact $\bar X$, let $n$ be an integer $\ge 0$ and let $\xi$ be a
spherical fibration on $X$, with fibers $\simeq S^d$ and with a
distinguished section which is a fiberwise cofibration. We have the
$\xi$--twisted stable $SW$-product on $s\sR\ld(\bar X,X)$ from
section~\ref{sec-twisted}.

\begin{defn} \label{defn-LKcontrol} 
{\rm Let $\Omega^n\bVLA\ubul((\bar X,X),\xi,n)$ be the homotopy pullback of
\[
\xymatrix@C=35pt{
\Omega^n\bVL\ubul((\bar X,X),\xi) \ar[r]^{\Xi} & \bA((\bar X,X),\xi,n)\thf &
\ar[l]_{\textup{incl.}} \bA((\bar X,X),\xi,n)\hf~.
}
\]
}
\end{defn}

Let $F((\bar X,X),\xi,n)=\Omega^{\infty+n}\bVLA\ubul((\bar X,X),\xi,n)$,
to be thought of as $\Omega^\infty$ of the spectrum $\Omega^n\bVLA\ubul((\bar X,X),\xi,n)$
just defined. Now, in order to
get a signature element $\sigma(\bar X,X)\in F((\bar X,X),\xi,n)$, we
have to throw in a finite domination assumption and a Poincar\ee
duality assumption. Suppose therefore that $S^0\times X$~, as a
retractive space over $X$, is finitely dominated in the controlled
sense. Suppose that in addition to the data $(\bar X,X)$ and $\xi$,
we are given a stable map
\[   \eta\co S^{n+d}\lra E/\!\!/X\,. \]
There
is a diagonal inclusion of $E/\!\!/X$ in $(S^0\times X)\odot_0(S^0\times X)$
for the $SW$ product in $s\sR^{(n)}$ which we are considering. The
composition of $\eta$ and the diagonal is then a 0-dimensional
visible symmetric structure on $S^0\times X$ as an object of
$s\sR^{(n)}$.

\begin{defn}
\label{defn-controlledhomotopycharacteristic}
 {\rm If this is nondegenerate, we call $(\bar X,X)$
together with $\xi$ and $\eta$ a controlled Poincar\ee duality
space. In that case let $\sigma(\bar X,X)\in F((\bar X,X),\xi,n)$ be
the element determined (as in
definition~\ref{defn-homotopycharacteristic}) by $S^0\times X$ with
the above nondegenerate visible symmetric structure.}
\end{defn}

The naturality properties of $\sigma(\bar X,X)$ are similar to those
of $\sigma(X)$ in definition~\ref{defn-homotopycharacteristic}. The
morphisms we are most interested in have the form
\[   ((\bar X_0,X_0),\xi_0,\eta_0) \leftarrow ((\bar X_1,X_1),\xi_1,\eta_1) \]
with an underlying map $f\co \bar X_1\to \bar X_0$ which takes $X_1$
to $X_0$ and restricts to a homeomorphism from $\bar X_1\smin X_1$
to $\bar X_0\smin X_0$. Existence of an inverse up to homotopy $\bar
X_0\to \bar X_1$ is required; the inverse and the maps in the
homotopy are subject to the same conditions as $f$.

\medskip\nin
\begin{rem} {\rm Definition~\ref{defn-controlledhomotopycharacteristic}
is good enough for our purposes, but it is not exactly a
generalization of definition~\ref{defn-homotopycharacteristic}
because of the compactness condition on $\bar X$. In the case where $\bar X=X$,
this amounts to compactness of $X$, which we did not assume in~\ref{defn-homotopycharacteristic}.
}
\end{rem}

\section{Excisive characteristics and signatures} \label{sec-excicharsig}
The goal here is very easy to
formulate:
\begin{itemize}
\item[(i)] For a closed $n$-manifold $X$ with normal bundle $\xi$,
we need to specify a preferred lift
\[\sigma\upr(X)\in F\upr(X,\xi,n)\]
of $\sigma(X)\in F(X,\xi,n)$, across a suitable assembly map.
\item[(ii)] For any compact control space $(\bar X,X)$ in which
$X$ is an $n$-manifold with normal bundle $\xi$ (without boundary
but not necessarily compact), we need to specify a preferred lift
\[ \sigma\upr(\bar X,X)\in F\upr(X,\xi,n) \]
of $\sigma(\bar X,X)\in F((\bar X,X),\xi,n)$, across a suitable
assembly map.
\end{itemize}
To clarify, $\sigma(X)\in F(X,\xi,n)$ and $\sigma(\bar X,X)\in
F((\bar X,X),\xi,n)$ are the characteristic elements of
definition~\ref{defn-homotopycharacteristic} and
definition~\ref{defn-controlledhomotopycharacteristic}. The space
$F\upr(X,\xi,n)$ is designed in such a way that its homotopy groups
are the homology groups of $X$ with \emph{locally finite} coefficients
$\Omega^n\bVLA\ubul(x,\xi_x,n)$, depending on $x\in X$.  

\medskip
The details are slightly more complicated, but we follow
\cite[\S7]{DwyerWeissWilliams} closely. In the case (i) we define
$F\upr(X,\xi,n)$ and the assembly map to $F(X,\xi,n)$ by means of a
homotopy fiber sequence
\[  F\upr(X,\xi,n) \lra F(X,\xi,n) \lra F(\JJ X,\xi,n) \]
where $\JJ X$ is the control space $(X\times[0,1],X\times[0,1[\,)$.
The second map in the sequence is induced by the inclusion of
$X\cong X\times\{0\}$ in $X\times[0,1[$~; the spherical fibration
$\xi$ is extended to $X\times[0,1[$ in the obvious way. We have the
definition of $F(\JJ X,\xi,n)$ from the end of the previous section.
In the case (ii) we introduce
\[ \JJ(\bar X,X)=\left(\frac{\bar
  X\times[0,1]}{\sim},X\times[0,1[\,\right) \]
where $\sim$ identifies points in $(\bar X\smin X)\times[0,1]$ with
the same coordinate in $\bar X\smin X$. Then again $F\upr((\bar
X,X),\xi,n)$ and the assembly map to $F((\bar X,X),\xi,n)$ are
defined by means of a homotopy fiber sequence
\[  F\upr((\bar X,X),\xi,n) \lra F((\bar X,X),\xi,n) \lra F(\JJ(\bar
X,X),\xi,n). \] The second map in the sequence is induced by the
inclusion of $(\bar X\times\{0\},X\times\{0\})$ in $\JJ(\bar X,X)$.

\medskip
In case (i), this leaves us with the task of trivializing the image
of $\sigma(X)$ under $F(X,\xi,n) \to F(\JJ X,\xi,n)$. We also need
to show that $F\upr(X,\xi,n)$ as defined above is excisive etc., so
that the forgetful map from $F\upr(X,\xi,n)$ to $F(X,\xi,n)$ can be
called an assembly map. In the more general case (ii), our task is
to trivialize the image of $\sigma(\bar X,X)$ under $F((\bar
X,X),\xi,n)\to F(\JJ(\bar X,X),\xi,n)$, and to establish excision
properties etc. for $F\upr((\bar X,X),\xi,n)$.

\bigskip\nin
\emph{Trivializing $\sigma(X)$ in $F(\JJ X,\xi,n)$.} We use
$F(\JJ\pt,\zeta,0)$ where $\pt$ means ``a point'' and $\zeta$ is the
trivial spherical fibration (with fiber $S^0$). We have $\sigma(\pt)
\in F(\JJ\pt,\zeta,0)$. Multiplication with $\sigma(X)$ can be
regarded as a based map
\[ F(\JJ\pt,\zeta,0) \lra F(\JJ X,\xi,n) \]
which takes $\sigma(\pt)$ to $\sigma(X)$. Hence it is enough to show
that $F(\JJ\pt,\zeta,0)$ is contractible. This can be proved by a
standard Eilenberg swindle argument. Think of $\JJ\pt$ as
$([0,\infty],[0,\infty[\,)$. Translation by +1 is an endofunctor $f$
of $s\sR(\JJ\pt)$. Let
\[  g= \bigvee_{i=0}^{\infty} f^i  \]
which is again an endofunctor of $s\sR(\JJ\pt)$. For the induced
self-maps $f_*$ and $g_*$ of $F(\JJ\pt,\zeta,0)$, we clearly have
$f_*\simeq \id$, hence $g_*f_*\simeq g_*$. But since $g\cong \id\vee
gf$, it is also true that $g_*\sim \id+g_*f_*$. Hence the identity
of $F(\JJ\pt,\zeta,0)$ is nullhomotopic. \qed

\medskip
\emph{Remark.} This line of reasoning simplifies
\cite[7.8]{DwyerWeissWilliams}. Observe that it is precisely the
manifold property of $X$ which ensures that there is a map
\emph{multiplication with $\sigma(X)$} from $F(\JJ\pt,\zeta,0)$ to
$F(\JJ X,\xi,n)$. The construction does therefore not generalize to
Poincar\ee duality spaces $X$.

\medskip
\emph{Remark.} What the argument really gives us is a lift of
$\sigma(X)\in F(X,\xi,n)$ to an element $\sigma\upr(X)$ in the
homotopy pullback of the following diagram of pointed spaces:
\[
\CD F(X,\xi,n) @>\textrm{incl.}>> F(\JJ X,\xi,n) @<\times\sigma(X)<<
F(\JJ\pt,\zeta,0).
\endCD
\]
Since the right--hand term is contractible, this homotopy pullback
is an acceptable substitute for the homotopy fiber of $F(X,\xi,n)
\to F(\JJ X,\xi,n)$.

\bigskip\nin
\emph{Trivializing $\sigma(\bar X,X)$ in $F(\JJ(\bar X,X),\xi,n)$.}
Again we use $F(\JJ\pt,\zeta,0)$ and the element $\sigma(\pt) \in
F(\JJ\pt,\zeta,0)$. Multiplication with $\sigma(\bar X,X)$ is a
based map
\[ F(\JJ\pt,\zeta,0) \lra F(\JJ(\bar X,X),\xi,n) \]
which takes $\sigma(\pt)$ to $\sigma(X)$. Since $F(\JJ\pt,\zeta,0)$
is contractible, this achieves the trivialization. \qed

\bigskip
In order to establish excision properties for $F\upr(X,\xi,n)$, we
begin by clarifying how the categories
\[ s\sR(X),\quad s\sR\ld(\JJ X),\quad s\sR\sG\ld(\JJ X) \]
fit together. (Recall that the decoration ``ld'' stands for
\emph{locally finitely dominated} and the $\sG$ in $\sR\sG$ stands for \emph{germs}.)
This will be done from the point
of view of Waldhausen's fibration and approximation theorems
\cite[\S1.6]{WaldRutgers}. We assume for now that $X$ is a compact
space, not necessarily a manifold. In addition to the standard
subcategory $w(s\sR\ld(\JJ X))$ of weak equivalences in $s\sR\ld(\JJ
X)$, we introduce a coarser notion of weak equivalence $\kappa$,
that is, a larger subcategory $\kappa(s\sR\ld(\JJ X))$. A morphism
in $s\sR\ld(\JJ X)$ is a $\kappa$-equivalence if its mapping cone is
equivalent to zero in $s\sR\sG\ld(\JJ X)$. Adopting Waldhausen's
notation, we write
\[ s\sR\ld(\JJ X)^{\kappa} \]
for the full subcategory of $s\sR\ld(\JJ X)$ consisting of all
objects $\kappa$-equivalent to the zero object. This is a Waldhausen
category with the usual $w$-equivalences as the weak equivalences.
We also write
\[ s\sR\ld(\JJ X)_{\kappa} \]
for $s\sR\ld(\JJ X)$ equipped with the coarser notion of weak
equivalence, that is, $\kappa$-equivalence.

\begin{lem}
\label{lem-Waldapprox} The forgetful functor $s\sR\ld(\JJ X)_{\kappa}\to s\sR\sG\ld(\JJ X)$ satisfies the hypotheses App1 and
App2 of Waldhausen's approximation theorem.
\end{lem}

\proof
Property App1 means that morphisms which are taken to weak
equivalences by the functor in question are already weak
equivalences. This is trivially true in our case. To establish App2, we abbreviate
$\sC_1=s\sR\ld(\JJ X)_{\kappa}$ and $\sC_2=s\sR\sG\ld(\JJ X)$. Fix an
object $Y_1$ in $\sC_1$ and a morphism $f\co Y_1\to Y_2$ in $\sC_2$.
We need to find an object $\bar Y_2$ in $\sC_1$ and morphisms
$g\co Y_1\to \bar Y_2$ in $\sC_1$ as well as $h\co \bar Y_2\to Y_2$ in $\sC_2$
such that $f=hg$ and $h$ is a weak equivalence. This will be called an
\emph{App2 factorization}. (Waldhausen's App2 also requires that $g$
be a cofibration. But if $g$ is not a cofibration, it can easily be
converted into one by means of a mapping cylinder construction, so
there is no need to pay any attention to that.) \newline
\emph{Case 1.} Here we assume in addition to all the above that $Y_1$ is the zero
object of $\sC_1$. It will be convenient to use the following terminology: an
object $Y$ of $\sC_2$ is \emph{clean} if there exists a weak equivalence $\bar Y\to Y$
where $\bar Y$ belongs to $\sC_1$. Our task is then to show that all objects of $\sC_2$
are clean. This is not easy. It is however mostly axiomatic business
and we start by summarizing related results of \cite{WeissHam}. \newline
Let $\sC$ be any Waldhausen category satisfying the saturation axiom for
weak equivalences \cite[\S0]{WeissHam}, and equipped with a cylinder functor satisfying the cylinder axiom.
For objects $C$ and $D$ of $\sC$ let $\sM(C,D)$ be the following category. The objects
are diagrams in $\sC$ of the form $C\to D'\leftarrow D$ in $\sC$, where the
second arrow is both a cofibration and a weak equivalence; a morphism from
$C\to D'\leftarrow D$ to $C\to D''\leftarrow D$ is a morphism $D'\to D''$ in $\sC$
making the diagram
\[
\xymatrix@R=7pt{
& D' \ar[dd] &  \\
C \ar[ur] \ar[dr] && D \ar[ul] \ar[dl] \\
& D'' &
}
\]
commutative. (Then $D'\to D''$ is a weak equivalence in $\sC$.) There is a
composition law in the shape of a functor
\[  \sM(D,E)\times\sM(C,D) \lra \sM(C,E)~. \]
It follows that we can make a new category $\sH\sC$ with the same objects as $\sC$,
where $\mor_{\sH\sC}(C,D)=\pi_0 B\sM(C,D)$. (The $\sH$ is meant to be reminiscent of \emph{homotopy}, but it
does not mean that $\sH\sC$ is obtained from $\sC$ by organizing morphism sets in $\sC$ into equivalence classes
for a relation called homotopy.) The functor $(C,D) \mapsto B\sM(C,D)$ of two variables
takes weak equivalences in either the left or right variable to homotopy equivalences of CW-spaces
\cite[\S1]{WeissHam}. Therefore:
\begin{itemize}
\item[] All weak equivalences in $\sC$ become invertible in $\sH\sC$. (Consequently
$\sH\sC$ can also be obtained from $\sC$ by formally adding inverses for all morphisms in $\sC$ which are
both cofibrations and weak equivalences.)
\end{itemize}
The functor $C\mapsto B\sM(C,D)$,
for fixed $D$ in $\sC$, takes cofiber squares in $\sC$ to homotopy pullback squares of CW-spaces
\cite[\S2]{WeissHam}. Here a
\emph{cofiber square} in $\sC$ is a (commutative) pushout square in which either the vertical arrows or
the horizontal arrows are cofibrations. This has the following consequence:
\begin{itemize}
\item[] If the functor \emph{suspension} from $\sH\sC$ to $\sH\sC$ is an equivalence of categories, then it
makes $\sH\sC$ into a triangulated category where the distinguished triangles are obtained (up to isomorphism)
from cofiber sequences in $\sC$:
\[
\xymatrix{  C \ar[r]^-{f} & D \ar[r] & \cone(f) \ar[r] & \Sigma C
}
\]
\end{itemize}
(We use Verdier's axioms for a triangulated category as in \cite{WeibelBook}.
The invertibility of the suspension functor is not used in the formulation of the remaining axioms.)
Now we add two more conditions on $\sC$. For the first, suppose that
\[
\CD
A @>>> B @>>> C \\
@VVV @VVV @VVV \\
A' @>>> B' @>>> C'
\endCD
\]
is a commutative diagram in $\sC$ whose rows are cofiber sequences. The condition is that if two
of the vertical arrows are weak equivalences, then so is the third. The second condition is that
if $A,B$ are objects of $\sC$ such that the unique morphism $\pt \to A\vee B$ is a weak equivalence,
then $\pt \to A$ and $\pt\to B$ are weak equivalences. When these conditions hold, in addition
to the above, then we have from \cite[\S3]{WeissHam}:
\begin{itemize}
\item[] A morphism $f\co C\to D$ in $\sC$ is a weak equivalence if and only if it becomes invertible in $\sH\sC$.
\end{itemize}
Now we specialize to $\sC:=\sC_2$. Our task is (still) to show that all objects of $\sC_2$ are clean.
We proceed by showing first that the class of clean objects in $\sC_2$
determines and is determined by a full triangulated subcategory of $\sH\sC_2$. This breaks down
into two statements:
\begin{itemize}
\item[(i)] If $Y\to Z$ is a weak equivalence in $\sC_2$ and $Z$ is clean, then $Y$ is clean.
\item[(ii)] If $Y\to Z$ is a cofibration where $Y$ and $Z$ are both clean, then the cofiber $Z/Y$ is also
clean.
\end{itemize}
Both of these are easy to verify and we only sketch the first. Suppose that $g\co \bar Z\to Z$ is a weak equivalence,
where $\bar Z$ belongs to $\sC_1$. Since $Y\to Z$ is a weak equivalence in $\sC_2$, there exists a homotopy
inverse $f\co Z\to Y$ in the controlled setting, relative to the reference space $X\times[0,1[~$ but not necessarily
\emph{over} it. Then $fg\co \bar Z\to Y$ is a controlled homotopy equivalence and we can view it as a morphism in $\sC_2$
with source in $\sC_1$ if we are willing to equip $\bar Z$ with a new retraction map to $X\times[0,1[~$. Therefore
$Y$ is clean. \newline
Now we know that the clean objects form a triangulated full subcategory of $\sH\sC_2$. Also, it follows easily
from the definition of $\sC_2$ that every object in $\sH\sC_2$ is a direct summand of a clean object.
Now general triangulated category principles imply that there is a subgroup of
the Euler-Poincar\ee group of $\sH\sC_2$ such that the clean objects are precisely those whose
Euler-Poincar\ee characteristic is in that subgroup.
(The Euler-Poincar\ee group of $\sH\sC_2$ is an abelian group with one generator $[Y]$ for every isomorphism class
of objects $Y$ of $\sH\sC_2$, and relations $[Y_1]=[Y_0]+[Y_2]$ for every distinguished triangle
$Y_0 \to Y_1 \to Y_2 \to \Sigma Y_0$.) To finish the argument, it is therefore enough to show that
the Euler-Poincar\ee group is trivial. But the Euler-Poincar\ee group is clearly a quotient of $K_0$ of
$\sC_2$, which we know is zero. \newline
\emph{Case 2.} Here we assume that $Y_1$ is in $s\sR\lf(\JJ X)_{\kappa}$, so it has a
locally finite and finite dimensional controlled CW structure.
By case 1, we can enlarge $f\co Y_1\to Y_2$ to a diagram
\[
\xymatrix{
& \bar Y_2 \ar[d]^e \\
Y_1 \ar[r]^-f & Y_2
}
\]
in $\sC_2$, where $\bar Y_2$ belongs to $\sC_1$ and $e$ is a weak equivalence in $\sC_2$.
We can find a CW subobject $Y'_1\subset Y_1$ such that the inclusion $Y'_1\to Y_1$ is
a weak equivalence in $s\sR\ld(\JJ X)_{\kappa}$ and $f|Y'_1$ admits a factorization
up to controlled homotopy (relative to zero sections) through $e$, say $f|Y'_1\simeq e\varphi$. The factorization
and the homotopy together allow us to think of the homotopy pushout (relative to zero sections) of
\[  \CD Y_1 @<\textup{incl.}<< Y'_1 @>\varphi>> \bar Y_2 \endCD \]
as an object $P$ in $s\sR\ld(\JJ X)_{\kappa}$. Then $\bar Y_2\subset P$ while $e\co \bar Y_2\to Y_2$
and $f\co Y_1\to Y_2$ extend simultaneously and canonically to a weak equivalence $P\to Y_2$ in $\sC_2$.
This solves our factorization problem. \newline
\emph{Case 3: the general case.} Let $r\co Y_1\to X\times[0,1[$ be the retraction. We are assuming that $Y_1$
can be dominated by a locally finite and finite dimensional controlled CW-object $W$. The domination
can be thought of in the following way. There are a morphism $W\to Y_1$ in $\sC_1$~, a controlled
map $Y_1\to W$ relative to zero sections which need not be a morphism in $\sC_1$~, and a controlled
homotopy $(u_t:Y_1\to Y_1)_{t\in[0,1]}$ relative to zero sections such that $u_0=\id$ and $u_1$ is the composition
$Y_1\to W\to Y_1$. Write $Y_{1,t}$ for $Y_1$ with the retraction $ru_t$.
The above description of the domination shows us that
the morphism $fu_1\co Y_{1,1}\to Y_2$ in $\sC_2$ admits a factorization
\[  Y_{1,1} \lra W \lra Y_2 \]
where the first arrow is a morphism in $\sC_1$ and the second arrow is in $\sC_2$.
Applying case 2 above to the second arrow, $W\to Y_2$, and pre-composing with $Y_{1,1}\to W$ afterwards, we deduce that
$fu_1\co Y_{1,1} \lra Y_2$ has a factorization
\[ \CD Y_{1,1} @>g>> \bar Y_2 @>h>> Y_2 \endCD \]
where $g$ is in $\sC_1$ and $h$ is a weak equivalence in $\sC_2$.
Enlarge $\bar Y_2$ to the mapping cylinder $Z(g)$ of $g$ (relative to zero sections),
with a retraction to $X\times[0,1[$ so that
the slice $Y_1\times t$ of the cylinder has the retraction $ru_t$. Then there is an inclusion
\[ Y_{1,0}\to Z(g) \]
in $\sC_1$, and it only remains to
extend $h\co \bar Y_2\to Y_2$ to a morphism $Z(g) \to Y_2$ in $\sC_2$.
Define the extension so that it agrees with $fu_t$ on the slice $Y_1\times t$ of the cylinder. \qed

\begin{lem}
\label{lem-otherapprox} The maps from $\bK(s\sR(X))$ to
$\bK(s\sR\ld(\JJ X)^{\kappa})$ and from $\bV\bL\ubul(s\sR(X))$ to
$\bV\bL\ubul(s\sR\ld(\JJ X)^{\kappa})$ induced by the inclusion
$s\sR(X)\to s\sR\ld(\JJ X)^{\kappa}$ are homotopy equivalences.
\end{lem}

\proof For the $K$-theory case, apply Waldhausen's approximation
theorem to the inclusion $s\sR(X\times[0,1[\,)\to s\sR\ld(\JJ
X)^{\kappa}$. The inclusion of $s\sR(X)$ in $s\sR(X\times[0,1[\,)$
clearly also induces a homotopy equivalence in $K$-theory. For the
$\bV\bL\ubul$-theory case, argue similarly, using a
$\bV\bL\ubul$-theory version of the approximation theorem. (This has
essentially the same hypotheses App1 and App2. The conclusion, that
the functor in question induces a homotopy equivalence of
$\bV\bL\ubul$-spectra, comes from a direct comparison of homotopy
groups.) \qed

\begin{cor}
\label{cor-twohofise} The inclusions $s\sR(X)\to s\sR\ld(\JJ X)$ and
the ``passage to germ'' functor $s\sR\ld(\JJ X)\to s\sR\sG\ld(\JJ
X)$ lead to homotopy fiber sequences of spectra,
\[
\CD \bA(X) @>>> \bA(\JJ X) @>>> \bA(\JJ X_{\infty})~,
\endCD
\]
\[
\CD \bV\bL\ubul(X,\xi,n) @>>> \bV\bL\ubul(\JJ X,\xi,n) @>>>
\bV\bL\ubul(\JJ X_{\infty},\xi,n).
\endCD
\]
\end{cor}

\proof By Waldhausen's fibration theorem there is a homotopy fiber
sequence
\[
\CD \bK(s\sR\ld(\JJ X)^{\kappa}) @>>> \bK(s\sR\ld(\JJ X) @>>>
\bK(s\sR\ld(\JJ X)_{\kappa}).
\endCD
\]
With the definition $\bA(\JJ X)=\bK(s\sR\ld(\JJ X)$ and the
identications of lemma~\ref{lem-Waldapprox} and
lemma~\ref{lem-otherapprox}, this turns into a homotopy fiber
sequence
\[
\CD \bA(X) @>>> \bA(\JJ X) @>>> \bA(\JJ X_{\infty}).
\endCD
\]
The same reasoning applies in the $VL\ubul$-theory case. Of course
one needs to know that Waldhausen's fibration theorem has an
analogue in $L$-theory. There is such an analogue, as follows.
\newline Let $\sC$ be any Waldhausen category with weak equivalences
$w\sC$ and with an $SW$ product $\odot$ satisfying the axioms of
\cite[\S2]{WWduality}. Suppose that another subcategory $\kappa\sC$
of $\sC$ with $w\sC\subset \kappa\sC\subset \sC$ is specified.
Suppose that $\sC$ with weak equivalences $\kappa\sC$ and the same
$SW$ product $\odot$ also satisfy the axioms of
\cite[\S2]{WWduality}. As before, let $\sC_{\kappa}$ stand for $\sC$
with the coarse notion $\kappa$ of weak equivalences, and
$\sC^{\kappa}$ for the full Waldhausen subcategory of $\sC$
consisting of the objects which are $\kappa$-equivalent to zero.
Then there is a homotopy cartesian square of spectra
\[
\CD
\bL\ubul(\sC^{\kappa}) @>>> \bL\ubul(\sC) \\
@VVV @VVV \\
\bL\ubul((\sC^{\kappa})_{\kappa}) @>>> \bL\ubul(\sC_{\kappa})
\endCD
\]
with contractible lower left-hand term. We shorten this as usual to
a homotopy fiber sequence
\[
\CD \bL\ubul(\sC^{\kappa}) @>>> \bL\ubul(\sC) @>>>
\bL\ubul(\sC_{\kappa})\,.
\endCD
\]
This is the basic fibration theorem for symmetric $L$-theory. There
is a version for quadratic $L$-theory and also one for visible
symmetric $L$--theory, when that is defined. \newline The proof of
the symmetric $L$-theory fibration theorem is as follows, in
outline. Suppose that $(C,D,\varphi)$ is a symmetric Poincar\ee pair
in $\sC$. In more detail, we assume that $D$ and $C$ are related by
a cofibration $D\to C$ and that $\varphi$ is a homotopy fixed point
for the action of $\ZZ/2$ on $\Omega^n(C\odot C/D\odot D)$,
satisfying the appropriate nondegeneracy condition. Then $n$ is the
formal dimension of the Poincar\'e pair. The symmetric structure
$\varphi$ on the pair $(C,D)$ descends to a symmetric structure
$\varphi/\partial \varphi$ on $C/D$, which may be degenerate. It is
a basic fact of Ranicki's algebraic theory of surgery that the
passage from the Poincar\ee pair $(C,D,\varphi)$ to the (single)
symmetric object $(C/D,\varphi/\partial\varphi)$ is reversible. We
have already used it in the proof of theorem~\ref{thm-homolothy};
see \cite{RanickiLMS1}, \cite{WeissLMS} and \cite{VogelNouvelle} for
more details on the inverse construction. These details on the
inverse construction imply that $C/D$ with the symmetric structure
$\varphi/\partial\varphi$ is nondegenerate in $\sC_{\kappa}$ if and
only if $D$ belongs to $\sC^{\kappa}$. Hence the bordism theory of
Poincar\ee pairs $(C,D,\varphi)$ with $C$ in $\sC$ and $D$ in
$\sC^{\kappa}$ is ``equivalent'' to the bordism theory of Poincar\ee
objects in $\sC_{\kappa}$. This amounts to a homotopy equivalence of
spectra, from the homotopy cofiber of
\[  \bL\ubul(\sC^{\kappa})\to \bL\ubul(\sC) \]
to $\bL\ubul(\sC_{\kappa})$. That in turn can be reformulated as a
homotopy fiber sequence of spectra $\bL\ubul(\sC^{\kappa})\to
\bL\ubul(\sC)\to\bL\ubul(\sC_{\kappa})$. \qed

\begin{cor}
\label{cor-lasthofise} The inclusions $s\sR(X)\to s\sR\ld(\JJ X)$
and the \emph{passage to germs} functor $s\sR\ld(\JJ X)\to
s\sR\sG\ld(\JJ X)$ lead to a homotopy fiber sequence of spectra
\[
\CD \bV\bLA\ubul(X,\xi,n) @>>> \bV\bLA\ubul(\JJ X,\xi,n) @>>>
\bV\bLA\ubul(\JJ X_{\infty},\xi,n).
\endCD
\]
\end{cor}

\proof What we are really saying is that $\bV\bLK\ubul$ turns the
square of Waldhausen categories with $SW$-duality
\begin{equation}  \label{eqn-SWstar}
\CD
s\sR\ld(\JJ X)^{\kappa} @>>> s\sR\ld(\JJ X) \\
@VVV @VVV \\
(s\sR\ld(\JJ X)^{\kappa})_{\kappa} @>>> s\sR\ld(\JJ X)_{\kappa}
\endCD 
\end{equation}
into a homotopy cartesian square of spectra, with
contractible lower left-hand term. This follows directly from the
analogous statements for $\bK$--theory and $\bV\bL\ubul$-theory, which
we have from the previous corollary, and the definition of
$\bV\bLK\ubul$ in terms of $\bK$ and $\bV\bL\ubul$. \qed

\bigskip\nin
\emph{Excision properties of $F\upr(X,\xi,n)$.} We have defined
$F\upr(X,\xi,n)$ as the homotopy fiber of a map
\[ F(X,\xi,n) \lra F(\JJ X,\xi,n). \]
Unravelling that and using corollary~\ref{cor-lasthofise}, we obtain
an identification of $F\upr(X,\xi,n)$ with
$\Omega^{\infty+1}\bV\bLK\ubul(s\sR\sG\ld(\JJ X))$. We know from
theorem~\ref{thm-visiblehomolothy} that $\Omega^{\infty}\bV\bL\ubul$
applied to $s\sR\sG\ld(\JJ X)$ gives a homotopy invariant and
excisive functor of $X$. Also, $\Omega^{\infty}\bK$ applied to
$s\sR\sG\ld(\JJ X)$ is an excisive and homotopy invariant functor of
$X$ by \cite[\S6-9]{Weissexci}. It follows that
$\Omega^{\infty}\bV\bLK\ubul$ applied to $s\sR\sG\ld(\JJ X)$ is an
excisive and homotopy invariant functor of $X$. \qed

\bigskip\nin
\emph{Excision properties of $F\upr((\bar X,X),\xi,n)$.} This is
similar to the case of $F\upr(X,\xi,n)$. We assume that $\bar X$ is
compact; we need not assume that $X$ is a manifold. There is a
commutative square of Waldhausen categories
\begin{equation} \label{eqn-SWtwostar}
\CD
s\sR\ld(\bar X,X) @>>> s\sR\ld(\JJ(\bar X,X)) \\
@VVV   @VVV  \\
\pt @>>> s\sR\sG\ld(\JJ X)
\endCD 
\end{equation}
which in the case $X=\bar X$ deviates very little from~(\ref{eqn-SWstar}).
Waldhausen's fibration theorem can be
applied to this, after some minor re-definitions. Hence
$F\upr((\bar X,X),\xi,n)$ is identified (via a chain of natural homotopy
equivalences) with
\[
\Omega^{\infty+1}\bV\bLK\ubul(s\sR\sG\ld(\JJ X)).
\]
Therefore it
is essentially a functor of the locally compact space $X$ alone, and
it is excisive in the locally finite sense of
theorem~\ref{thm-homolothy} and theorem~\ref{thm-visiblehomolothy}.
Indeed this follows from theorem~\ref{thm-visiblehomolothy} and the
analogous theorem for $K$-theory, \cite[\S6-9]{Weissexci}. \qed

\bigskip\nin
\emph{Naturality properties of $\sigma\upr(X)$.} The naturality
properties of $\sigma\upr$ are analogous to those of $\sigma$ in the
previous section. That is, a diagram of closed manifolds,
homeomorphisms and stable normal bundle isomorphisms
\[ (X_0,\xi_0) \leftarrow (X_1,\xi_1) \leftarrow \cdots
\leftarrow (X_{k-1},\xi_{k-1})\leftarrow (X_k,\xi_k) \] determines a
map from $\Delta^k$ to $F\upr(X_0,\xi_0,n)$. This assignment extends
$\sigma\upr$ and commutes with the usual face and degeneracy
operators acting on such diagrams. (We omit the details, except for
pointing out that each $\xi_i$ can be viewed as a spherical
fibration on $X_i$ via fiberwise one-point compactification.)

\bigskip\nin
\emph{Naturality properties of $\sigma\upr(\bar X,X)$.} Let $(\bar
X_i,X_i)$ for $i=0,1,\dots,k$ be compact control spaces where each
$X_i$ is an $n$-manifold. Suppose that they are arranged in a
diagram of homeomorphisms of control spaces and stable normal bundle
isomorphisms
\[ ((\bar X_0,X_0),\xi_0) \leftarrow ((\bar X_1,X_1),\xi_1) \leftarrow \cdots\leftarrow ((\bar
X_k,X_k),\xi_k). \] The diagram then determines a map from $\Delta^k$ to
$F\upr((\bar X_0,X_0),\xi_0)$. This assignment commutes with the
usual face and degeneracy operators acting on such diagrams.

\medskip
\section{Algebraic approximations to structure spaces: Set-up}
\label{sec-approxsetup}

\begin{notn} {\rm We use the symbol $\%$ in the \emph{subscript} position
to describe homotopy fibers of assembly maps. For example:
$\bL\ubul\lpr(X)$ is the homotopy fiber of the assembly map
$\bL\ubul\upr(X)\to \bL\ubul(X)$, assuming $X$ has the homotopy type
of a CW-space. Similarly: $\bL\lbul\lpr(X)$ is the homotopy fiber of
the assembly map in quadratic $L$-theory and
$\bLA\lbul\lpr(X,\xi,n)$ is the homotopy fiber of the assembly map
in quadratic $LA$-theory. \newline
 From now on we will often suppress the $\bar X$ in a control space
$(\bar X,X)$ and make up for that with a semicolon followed by a
$c$ to indicate a controlled context. For example, instead of
writing $F\upr((\bar X,X),\xi,n))$
and $F((\bar X,X),\xi,n))$ as in section~\ref{sec-excicharsig},
we may write
\[ F\upr(X,\xi,n\cc)~,\quad F(X,\xi,n\cc). \]
In that spirit, $F\lpr(X,\xi,n\cc)$ is
the homotopy fiber of the assembly map alias forgetful map
$F\upr(X,\xi,n\cc)\to F(X,\xi,n\cc)$ which we have
from section~\ref{sec-excicharsig}. This assumes that $X$ is
the nonsingular part of a control space $(\bar X,X)$
with compact $\bar X$.
}
\end{notn}

\begin{con}
\label{con-mainmap} {\rm Let $M$ be a compact $m$-manifold with
normal bundle $\nu$. We introduce a map $\varphi$ as in~(\ref{eqn-mainmap}) from $\sS(M)$ to
\[ F\lpr(M,\nu,m)\,=\,\Omega^{\infty+m}\bV\bLA\ubul\lpr(M,\nu,m)\,\simeq\,
\Omega^{\infty+m}\bLA\lbul\lpr(M,\nu,m). \] }
\end{con}

\medskip\nin
\emph{Remark.} The setting is \emph{away from} $\partial M$.
Points of $\sS(M)$ correspond to homotopy equivalences of pairs
$f\co (N,\partial N)\to (M,\partial M)$ where the induced map
$\partial N\to\partial M$ is a homeomorphism.
We tend to think of $F\lpr(M,\nu,m)$ as the homotopy fiber of the assembly map
\[ \Omega^{\infty+m}\bV\bLA\ubul\upr(M,\nu,m) \lra \Omega^{\infty+m}\bV\bLA\ubul(M,\nu,m) \]
over the visible symmetric signature $\sigma(M)$. Because this homotopy fiber has a canonical base point, determined
by $\sigma\upr(M)$ etc., it can easily be identified by a translation with the homotopy fiber of assembly over
the base point, whenever that is convenient.

\begin{con}
\label{con-degree}
{\rm Notation being as in the previous construction,
there is a ``local degree'' homomorphism from
$\pi_m\bL\lbul\upr(M,\nu)$ to the group $L_0(\ZZ)^{\pi_0M}$. It is onto. }
\end{con}

This is the composition of Poincar\ee duality, $\pi_m\bL\lbul\upr(M,\nu)\cong
H^0(M;\bL\lbul(\pt))$, with an evaluation map from $H^0(M;\bL\lbul(\pt))$ to
$L_0(\ZZ)^{\pi_0M}$. We admit that $L_0(\ZZ)$ could also be described as $L_0(\pt)$.

\medskip
\begin{thm}
\label{thm-main} Assume $\dim(M)\ge 5$. The diagram
\[
\CD \sS(M) @>\varphi>>
\Omega^{\infty+m}\bLA\lbul\lpr(M,\nu,m) @>\textup{ local degree }>>
L_0(\ZZ)^{\pi_0M}
\endCD
\]
is a homotopy fiber sequence in the concordance stable range.
\end{thm}

\medskip\nin
In this theorem, $L_0(\ZZ)\cong 8\ZZ$ is viewed as a discrete space.
The local degree is defined on $\pi_m\bLA\lbul\lpr(M,\nu,m)$ via the
forgetful map to $\pi_m\bL\lbul\upr(M,\nu)$. Details on the meaning
of ``concordance stable range'' are given in the definition which
follows.

\begin{defn}
\label{defn-concordstable}
 {\rm For a compact manifold $N$ let $k_N$ be the minimum of
all positive integers
$k$ such that the stabilization map of concordance spaces
\[ \Omega\sH(N\times D^i) \lra \Omega\sH(N\times D^{i+1}) \]
is $k$-connected for $i=0,1,2,\dots~$. The precise meaning of
theorem~\ref{thm-main} is that
\begin{itemize}
\item the map from $\pi_0\sS(M)$ to
$\ker[\pi_m\bLA\lbul\lpr(M,\nu,m)\to L_0(\ZZ)^{\pi_0M}]$ determined by $\varphi$
is bijective;
\item the restricted map, from a component of $\sS(M)$ represented by a homotopy equivalence
$(N,\partial N)\to(M,\partial M)$ to the corresponding component of
$\Omega^{\infty+m}\bLA\lbul\lpr(M,\nu,m)$, is $(k_N+1)$-connected.
\end{itemize}
}
\end{defn}


\medskip\nin
\emph{Remark.} It follows from \cite{Igusa88} and smoothing theory
that
\[ k_N\ge \min\left\{\frac{m-7}{2},\frac{m-4}{3}\right\} \]
if $N$ admits a smooth structure,
$m=\dim(N)$. Similar estimates for non-smooth manifolds were part of the topology
folklore in the 1970's, but completely convincing proofs of these have
apparently not been found. Unfortunately, even if $M$ admits a
smooth structure, there may be some components of $\sS(M)$,
represented by $(N,\partial N)\to (M,\partial M)$ say, whose source
manifold $N$ does not admit a smooth structure.

\medskip\nin
Our only goal in the rest of the paper is to prove theorem~\ref{thm-main}. An overview
of the proof has already been given in section~\ref{sec-outline}. We are going to break this up into
smaller pieces.

\begin{lem}
\label{lem-maincomponents}
The map in construction~\ref{con-mainmap} identifies $\pi_0\sS(M)$ with
the kernel of the local degree homomorphism
from $\pi_m\bLA\lbul\lpr(M,m)$ to $L_0(\ZZ)^{\pi_0M}$.
\end{lem}

For $i\ge 0$ we introduce the controlled
structure space $\sS(M\times\RR^i\cc)$, using the compactification
$M*S^{i-1}$ of $M\times\RR^i$ to define the control criteria. Let
\[ \sS^\red(M\times\RR^i\cc) \subset \sS(M\times\RR^i\cc) \]
be the union of the connected components which are reducible, i.e., in the image
of $\pi_0\sS(M)\to \pi_0\sS(M\times\RR^i\cc)$.

\begin{con}
\label{con-extendedmainmap} {\rm For $i\ge 0$ we construct a map
$\varphi$ from $\sS^\red(M\times\RR^i\cc)$ to
\[ F\lpr(M\times\RR^i,\nu,m+i \cc)\,\simeq\,\Omega^{\infty+m+i}\bLA\lbul\lpr(M\times\RR^i,\nu,m+i\cc)~.
\]
This map agrees with construction~\ref{con-mainmap}
when $i=0$. The resulting square~(\ref{eqn-mainmapsquare}) commutes for every $i\ge0$. (We are actually
slightly more general than in~(\ref{eqn-mainmapsquare}) because $M$ can have a nonempty boundary.)
}
\end{con}

Let $f\co N\to M$ be a homotopy equivalence of compact $m$-manifolds, relative to a homeomorphism
$\partial N\to \partial M$. Thus $f\co N\to M$ determines a point in $\sS(M)$.

\begin{lem}
\label{lem-oldhofiber} The homotopy fiber of
\[
\CD
\sS^\red(M\times\RR^i\cc) @>{\times\RR}>> \sS^\red(M\times\RR^{i+1}\cc)
\endCD
\]
over the point determined by $f\times\RR^{i+1}\co N\times\RR^{i+1}\to M\times\RR^{i+1}$ is homotopy
equivalent to a union of components of $\sH(N\times\RR^i\cc)$. Indeed it is
homotopy equivalent to the pullback of the following diagram.
\[
\xymatrix@C=50pt@M=10pt{ & \sS^\red(M\times\RR^i\cc) \ar[d]^-{\textup{inclusion}} \\
\sH(N\times\RR^i\cc) \ar[r]^-{\textup{upper bdry}} & \sS(N\times\RR^i\cc)
}
\]
\end{lem}

\begin{lem}
\label{lem-newhofiber} The homotopy fiber of
\[
\xymatrix@M=10pt{
F\lpr(M\times\RR^i,\nu,m+i\cc) \ar[r]^-{\times\RR} &
F\lpr(M\times\RR^{i+1},\nu,m+i+1\cc)
}
\]
is homotopy equivalent to $\Omega^{\infty}\bA\lpr(M\times\RR^i\cc)$.
\end{lem}

In combination with some results from~\cite{DwyerWeissWilliams}, lemmas~\ref{lem-oldhofiber} and~\ref{lem-newhofiber}
lead to the following technical statement for the square~(\ref{eqn-mainmapsquare}). Again we allow nonempty $\partial M$.
\emph{Notation}: Select
$x\in \sS^\red(M\times\RR^{i+1}\cc)$ with image
\[ y\in \Omega^{\infty+m+i+1}\bLA\lbul\lpr(M\times\RR^{i+1},\nu,m+i+1\cc)~. \]
Let $\Phi_x$ and $\Psi_y$ be the vertical homotopy fibers over $x$ and $y$ respectively in~(\ref{eqn-mainmapsquare}).
Select $a\in \pi_0\Phi_x$ with image $b\in \pi_0\Psi_y$. Let $\Phi_{x,a}$ and $\Psi_{y,b}$ be the
connected components determined by $a$ and $b$, respectively. Select an element of $\pi_0\sS(M)$
whose image in $\pi_0\sS^\red(M\times\RR^i\cc)$ agrees with the image of $a$ under the forgetful map.
Represent that element of $\pi_0\sS(M)$ by a homotopy equivalence $N\to M$ (restricting to a
homeomorphism $\partial N\to \partial M$).

\begin{prop} \label{prop-whatvertical}
The map $\Phi_x\to \Psi_y$ determined by~(\ref{eqn-mainmapsquare}) induces an injection on $\pi_0$.
The restricted map $\Phi_{x,a}\to \Psi_{y,b}$ is $(k_N+i+1)$-connected.
\end{prop}

\medskip\nin
One more ingredient that we need for our downward induction
procedure is the compatibility of the
Casson-Sullivan-Wall-Quinn-Ranicki homotopy equivalence~(\ref{eqn-industart}) and
diagram~(\ref{eqn-mainmapladder}). In the following formulation we make no distinction between
(telescopic) homotopy colimits and ordinary colimits.

\begin{lem}
\label{lem-inthelimit} In the limit $i=\infty$, the horizontal maps of
diagram~(\ref{eqn-mainmapladder}) agree with the map~(\ref{eqn-industart}): there is a
homotopy commutative square
\[
\xymatrix@M=5pt{
{\rule{0mm}{5.5mm}\colimsub{i}
\sS(M\times\RR^{i+1}\cc)} \ar[r]^-{(\ref{eqn-industart})}_-{\simeq}
& {\rule{0mm}{5.5mm}\colimsub{i}\Omega^{\infty+m+i}\bL^{cs}\lbul\lpr(M\times\RR^i,\nu\cc)}     \\
{\rule{0mm}{5.5mm}\colimsub{i}\sS^\red(M\times\RR^{i+1}\cc)} \ar[r]^-{(\ref{eqn-mainmapladder})}
\ar[u]^-{\textup{incl.}} & \ar[u]^-{\simeq}
{\rule{0mm}{5.5mm}\colimsub{i}\Omega^{\infty+m+i}\bLA\lbul\lpr(M\times\RR^i,\nu,m+i\cc)}
}
\]
\end{lem}

\medskip
\proof[Proof of theorem~\ref{thm-main}, assuming lemma~\ref{lem-maincomponents}, lemma~\ref{lem-inthelimit} and
proposition~\ref{prop-whatvertical}]
This follows the outline given in section~\ref{sec-outline}. Note the added precision coming from
definition~\ref{defn-concordstable}.  \qed

\section{Algebraic approximations to structure spaces: Constructions} \label{sec-approxdetails}
The constructions are based on one general principle which we recall from \cite{DwyerWeissWilliams}.
Let $\sC$ and $\sD$ be small categories, $\sC\subset \sD$. Let $U_0$ and $U_1$ be functors from $\sD$ to spaces,
with a natural transformation $a\co U_0\to U_1$. Suppose that $U_0$ and $U_1$ take every morphism
in $\sD$ to a homotopy equivalence of spaces. Then every choice of point $x=(x_0,x_1,\gamma)$ in the homotopy pullback of
\[
\xymatrix{
& \holim~U_1 \ar[d]^{\textup{res}} \\
\holim~U_0|\sC \ar[r]^{a_*} & \holim~U_1|\sC
}
\]
and choice of object $\delta$ in $\sD$ determine a map
\begin{equation} \label{eqn-DWWprinciple}
\hofiber_\delta[B\sC\to B\sD] \lra \hofiber_{x_1(\delta)}[U_0(\delta)\to U_1(\delta)]~.
\end{equation}
Here $x_0\in \holim~U_0|\sC$ and $x_1\in \holim~U_1$~, while $\gamma\co [0,1]\to \holim~U_1|\sC$
is a path connecting the images of $x_0$ and $x_1$. We have written $x_1(\delta)\in U_1(\delta)$ for the coordinate
of $x_1$ in $U_1(\delta)$. \newline
\emph{Sketch of construction of~}(\ref{eqn-DWWprinciple}). Because $U_0$ and $U_1$ take every morphism
in $\sD$ to a homotopy equivalence of spaces, the projections
\[ \hocolim U_0 \lra B\sD~, \qquad \hocolim U_1 \lra B\sD \]
are quasi-fibrations. The element $x_1\in \holim~U_1$ determines
a section $s_1$ of the fibration associated to $\hocolim U_1\to B\sD$, and $x_0$
determines a section $s_0$ of the fibration associated to $\hocolim U_0|\sC\to B\sC$. Pulling these data
back to
\[ \sS:=\hofiber_\delta[B\sC\to B\sD] \]
yields two trivialized fibrations on $\sS$ with fibers $U_0(\delta)$ and
$U_1(\delta)$, respectively, and two sections $t_0$ and $t_1$ of these. We may view
$t_0$ and $t_1$ as maps $\sS\to U_0(\delta)$ and $\sS\to U_1(\delta)$, respectively.
The path $\gamma$ determines a homotopy $a_*t_0 \simeq t_1$.
The map $t_1$ is trivialized (equipped with a homotopy to a constant map)
because, as a section, it is the pullback of a section defined on $B\sD$. \qed
\newline
\emph{Variants.} We allow $\sC$ and $\sD$ to be categories \emph{internal to the category of
topological spaces}. This means that $\sD$ for example has a space of objects and a space of morphism,
the maps \emph{source} and \emph{target} from morphism space to object space are continuous, and so on.
In such a case we try our best to ensure that the functors $U_0$ and $U_1$ from $\sD$ to spaces factor through
the discrete category $\pi_0\sD$ with object set $\pi_0\textup{ob}(\sD)$ and morphism set $\pi_0\textup{mor}(\sD)$.

\bigskip
We now specialize~(\ref{eqn-DWWprinciple}) to obtain construction~\ref{con-mainmap}
Suppose to begin that $M$ has empty boundary.

\begin{defn}
\label{defn-manifoldcats}
{\rm Let $\sD$ be the following category. The objects are finitely dominated
Poincar\ee duality spaces $X$ of formal dimension $m$, together with
a spherical fibration $\xi\co E\to X$ with fiber $\simeq S^d$~, a
preferred section for that, and a stable map $\eta\co S^{m+d}\to
E/X$ which carries a fundamental class. Mostly for convenience, we
require $X$ to be compact Hausdorff and homotopy equivalent to $M$. A morphism from
$(X_1,\xi_1,\eta_1)$ to $(X_0,\xi_0,\eta_0)$ is a pair $(u,v)$ where
$u$ is a homotopy equivalence from $X_1$ to $X_0$ and
\[ v \co \Sigma^{d_1-d_0}_{X_1}E_1 \lra E_0 \]
is a homotopy equivalence which covers $u$, respects the zero
sections and satisfies $v\eta_1=\eta_0$. (It is understood that
$d_1\ge d_0$~, where $d_1$ and $d_0$ are the formal fiber dimensions
of $\xi_1\co E_1\to X_1$ and $\xi_0\co E_0\to X_0$.)
\newline
We allow continuous variation of $\eta$ in objects $(X,\xi,\eta)$.
This makes the set of objects of $\sD$ into a space. (Some
conditions on underlying sets should be added to our definition of
objects to ensure that the class of objects of $\sD$ is indeed a set.) We
also allow continuous variation of $v$ in morphisms
$(u,v)$ as above. This makes the set of morphisms of $\sD$ into a space, in such a way
that ``source'' and ``target'' are continuous maps from the morphism space
to the object space. Hence $\sD$ is a topological category. (By allowing
continuous variation of the $\eta$'s and the $v$'s~, we achieve that $B\sD$
is homotopy equivalent to $BG(M)$, where $G(M)$ is the topological monoid of
homotopy automorphisms of $M$. This follows from the uniqueness of Spivak
normal fibrations~; see section~\ref{sec-updates}.)
\newline
The category $\sD$ has a subcategory $\sC$ which is defined like this.
An object $(X,\xi,\eta)$ of $\sD$, as above, belongs to $\sC$
precisely if $X$ is a closed $m$-manifold, $\xi$ is a sphere bundle
and $\eta\co S^{m+d}\to E/X$ restricts to a homeomorphism from
$\eta^{-1}(E\smin X)$ to $E\smin X$. (The condition on $\eta$ means
that $\eta$ is the ``Thom collapse'' associated with an embedding of
$N$ in some euclidean space.) A morphism $(u,v)$ in $\sD$, as above,
belongs to $\sC$ if its source and target are in $\sC$, and both $u$
and $v$ are homeomorphisms. Continuous variation of $\eta$ in objects
$(X,\xi,\eta)$ and of $v$ in morphisms $(u,v)$ are allowed as before.
The result is that $B\sC$ is homotopy equivalent to a disjoint union of
classifying spaces $B\Homeo(N)$, where $N$ runs through a maximal set
of pairwise non-homeomorphic closed $m$-manifolds which are homotopy
equivalent to $M$. Recall from section~\ref{sec-discrete}) that $\Homeo(N)$ is
the homeomorphism group of $N$ with the discrete topology.
}
\end{defn}

The manifold $M$ itself, equipped with a euclidean normal bundle $\nu$
etc., can be viewed as an object of $\sC$. We note that
\[ \Str(M):=\hofiber_M[ B\sC\hookrightarrow B\sD\,] \]
is a good combinatorial model for $\sS(M)$. More precisely,
$\Str(M)$ comes with a forgetful
map to $\sS(M)$ which is a homology equivalence. We like to think
of that map as an inclusion.
\newline The homotopy invariant signature $\sigma$ which we have
constructed in section~\ref{sec-hocharsig}
determines or is a point in $\holim\, F|\sD$. To be more precise,
$F|\sD$ means the functor taking $(X,\xi,\eta)$ in $\sD$ to
$F(X,\xi,m)$ as defined in section~\ref{sec-hocharsig}~; then we have
$\sigma(X)\in F(X,\xi,m)$ with the naturality properties discussed
at the end of that same section. Similarly, the excisive
signature $\sigma\upr$ which we have constructed in
section~\ref{sec-excicharsig} is a point in
$\holim\, F\upr|\sC$. Here $F\upr|\sC$ means the functor taking
$(N,\xi,\eta)$ to $F\upr(N,\xi,m)$ as defined in
section~\ref{sec-excicharsig}. \newline
We apply~(\ref{eqn-DWWprinciple}) with this choice of $\sD$ and $\sC$, with $U_0=F\upr|\sD$
and $U_1=F|\sD$, and with $x_1=\sigma$, $x_0=\sigma\upr$. The result is a map
\begin{equation} \label{eqn-instanceDWWprin}
 \Str(M) \lra \hofiber_{\sigma(M)}[~F\upr(M,\nu,m)\to F(M,\nu,m)~]~.
\end{equation}
We can write~(\ref{eqn-instanceDWWprin}) in the form
\[ \Str(M) \lra F\lpr(M)~. \]
By obstruction theory, for which we can use that $F\lpr(M)$ is an $H$-space,
this map has an extension to the mapping cylinder of $\Str(M)\to \sS(M)$.
The extension is unique up to contractible choice. Restrict to $\sS(M)$ to obtain
\begin{equation} \label{eqn-thatmap} \varphi_M\co \sS(M) \lra F\lpr(M,\nu,m)~.
\end{equation}
To conceal the dependence on a contractible choice, and also to make the
contractibility of the choice more obvious, we could replace the target $F\lpr(M,\nu,m)$
by the homotopy equivalent
\[  \hocolim[~\sS(M) \leftarrow \Str(M) \to F\lpr(M,\nu,n)~]. \]

\medskip
A systematic way to extend the construction to the case where $M$ has a nonempty
boundary is to use algebraic $K$-theory of pairs.
This requires a few definitions which we have banished to appendix~\ref{sec-Kpairs}.
Meanwhile the notation that comes with these definitions is self-explanatory. The construction
of $\varphi_M$ in~(\ref{eqn-thatmap}) generalizes mechanically to produce a map
\[  \varphi_{(M,\partial M)}\co \sS(M,\partial M) \lra F\lpr(\partial M\subset M,\nu,m) \]
where the definition of the target is based on $K$-theory and $L$-theory of pairs, such as $(M,\partial M)$.
The maps $\varphi_{(M,\partial M)}$ and $\varphi_{\partial M}$ fit into a commutative square
\begin{equation} \label{eqn-thatrelativemap}
\begin{split}
\xymatrix@M=10pt@C=50pt{
\sS(M,\partial M) \ar[r]^-{\varphi_{(M,\partial M)}} \ar[d]^-{\textup{forget}} &
F\lpr(\partial M\subset M,\nu,m) \ar[d]^-{\textup{forget}} \\
\sS(\partial M) \ar[r]^-{\varphi_{\partial M}} &  F\lpr(\partial M,\nu,m-1).
}
\end{split}
\end{equation}
We pass to vertical homotopy fibers (over the base points in the lower row) and obtain
a map of the form
\begin{equation} \label{eqn-relphi1}
\varphi_M\co \sS(M) \lra F\lpr(M,\nu,m)~.
\end{equation}
There is a slight generalization which we shall also need. Let $M$ be a compact manifold
such that $\partial M$ is the union of two codimension zero submanifolds $\partial_0M$ and $\partial_1M$,
with $\partial_0M\cap \partial_1M=\partial\partial_0M=\partial\partial_1M$. There is a structure space
\[  \sS\left(\CD \partial(\partial_0M) @>>> \partial_0M \\ @VVV @VVV \\ \partial_1M @>>> M \endCD \right) \]
which fits into a homotopy fiber sequence
\[
\xymatrix{
{\sS(M,\partial_1M)} \ar[r] &
{\sS\left(\CD \partial(\partial_0M) @>>> \partial_0M \\ @VVV @VVV \\ \partial_1M @>>> M \endCD \right)}
\ar[r] &
{\sS(\partial_0M,\partial\partial_0M)}
}
\]
We can therefore make a map of the form
\begin{equation} \label{eqn-relphi2}
\varphi_{(M,\partial_1M)}\co \sS(M,\partial_1M) \lra  F\lpr(\partial M_1\subset M,\nu,m)
\end{equation}
by passing to vertical homotopy fibers in the square
\begin{equation} \label{eqn-thatotherrelativemap}
\begin{split}
\xymatrix{
{\sS\left(\CD \partial\partial_0M @>>> \partial_0M \\ @VVV @VVV \\ \partial_1M @>>> M \endCD \right)}
\ar[r] \ar[d] &
{F\lpr\left(\CD \partial\partial_0M @>>> \partial_0M \\ @VVV @VVV \\ \partial_1M @>>> M \endCD,\nu,m\right)} \ar[d] \\
{\sS(\partial_0M,\partial\partial_0M)} \ar[r] & {F\lpr(\partial\partial_0M\subset \partial_0M,\nu,m-1)}
}
\end{split}
\end{equation}
This completes construction~\ref{con-mainmap}.

\bigskip
Construction~\ref{con-extendedmainmap} uses very similar ideas and recycled notation.
There is a small complication due to the fact that we do not have a good discrete homology approximation theorem
for the topological group of controlled automorphisms of $M\times\RR^i$, assuming that $M$ is compact.
Instead we have such a theorem (in section~\ref{sec-discrete}) for
the topological group of homeomorphisms $M*S^{i-1}\to M*S^{i-1}$
taking $S^{i-1}$ to $S^{i-1}$, and we have to make the best of that.
Suppose to begin with that $M$ is closed.

\begin{defn} {\rm Generalizing definition~\ref{defn-manifoldcats},
we introduce a category $\sD$ whose objects are certain control
spaces $(\bar X,X)$ with compact Hausdorff $\bar X$, together with a
spherical fibration $\xi\co E\to X$ with fiber $\simeq S^d$, a preferred
section for that, and a stable map
$\eta\co S^{m+d}\to E/\!\!/X$. Compare definition~\ref{defn-strangequot}.
These data are required
to satisfy the conditions of definition~\ref{defn-controlledhomotopycharacteristic}
for a controlled Poincar\ee duality space. In addition we require
that $\bar X\smin X=S^{i-1}$ and that $(\bar X,X)$ be homotopy equivalent as a control space
to $(M*S^{i-1},M\times\RR^i)$; and moreover we mean a homotopy equivalence such
that the maps and homotopies involved restrict to the identity on the singular sets $S^{i-1}$.
A morphism from
$(\bar X_1,X_1,\xi_1,\eta_1)$ to $(\bar X_0,X_0,\xi_0,\eta_0)$ is a pair $(u,v)$ where
\begin{itemize}
\item $u$ is a map of control spaces from $(\bar X_1,X_1)$ to $(\bar X_0,X_0)$ which restricts
to the identity $S^{i-1}\to S^{i-1}$ on singular sets;
\item $v\co\Sigma^{d_1-d_0}_{X_1}E_1 \lra E_0$
covers $u|X_1$, respects the zero
sections and satisfies $v\eta_1=\eta_0$.
\end{itemize}
Continuous variation of
$\eta$ in objects $(\bar X,X,\xi,\eta)$ and continuous variation of $v$
in morphisms $v$ is allowed. The result is that
\begin{equation} \label{eqn-ctrlsemiprod}
B\sD\,\,\simeq\,\, BG(M\times\RR^i\cc)\,\simeq\, BG(M)
\end{equation}
where $G(M\times\RR^i\cc)$ is the grouplike topological monoid of controlled
homotopy automorphisms of $M\times\RR^i$, and $G(M)$ is the grouplike topological monoid of
homotopy automorphisms of $M$.
\newline
The category $\sD$ has a subcategory $\sC$
which is defined like this.
An object $(\bar X,X,\xi,\eta)$ of $\sD$ belongs to $\sC$
precisely if $(\bar X,X)$ is homeomorphic
to $(N*S^{i-1},N\times\RR^i)$ for some closed $m$-manifold $N$,
the spherical fibration $\xi$ is a sphere bundle
and $\eta\co S^{m+d}\to E/\!\!/X$ restricts to a homeomorphism from
$\eta^{-1}(E\smin X)$ to $E\smin X$.  A morphism $(u,v)$ in $\sD$, as above,
belongs to $\sC$ if its source and target are in $\sC$, and both $u$
and $v$ are homeomorphisms. Continuous variation of
$\eta$ in objects $(\bar X,X,\xi,\eta)$ and continuous variation of $v$
in morphisms $v$ is allowed. The result is that
\[ B\sC \simeq \coprod_{\beta} B\Homeo(N_\beta\times\RR^i\cc). \]
Here each $N_\beta$ is a closed $m$-manifold homotopy equivalent to $M$.
We select a maximal set of such manifolds $N_{\beta}$ such that the $N_\beta\times\RR^i$
are pairwise controlled non-homeomorphic. Furthermore
$\Homeo(N_\beta\times\RR^i)$ is the discrete group of controlled
homeomorphisms from $N_\beta\times\RR^i$ to itself.
}
\end{defn}

\medskip\nin
Generalizing from $i=0$ to $i\ge 0$, we obtain functors
$F|\sD$ and $F\upr|\sD$ and elements
\[ \sigma\in \holim\, F|\sD~, \qquad \sigma\upr\in \holim\, F\upr|\sC~. \]
They lead as before to a map
\begin{equation} \label{eqn-instanceDWWprin2} \Str(M\times\RR^i\cc) \lra F\lpr(M\times\RR^i\cc) \end{equation}
where
\[ \Str(M\times\RR^i\cc):=\hofiber_{M\times\RR^i}[ B\sC\hookrightarrow B\sD\,]. \]
There is an obvious inclusion map from $\Str(M\times\RR^i\cc)$ to $\sS^\red(M\times\RR^i\cc)$.
It is not clear whether this is a homology equivalence. But there is a fix. Let
\[ \Gamma:= \hofiber[\,\Homeo(S^{i-1})\to \homeo(S^{i-1})\,] \]
where $\homeo(S^{i-1})$ is the homeomorphism group of $S^{i-1}$ and
$\Homeo(S^{i-1})$ is the underlying discrete group. Then $\Gamma$ is a topological group
whose classifying space $B\Gamma$ has the homology of a point. The group $\Gamma$ acts, via projection to
the discrete group $\Homeo(S^{i-1})$, on both $\Str(M\times\RR^i\cc)$ and
$F\lpr(M\times\RR^i\cc)$, preserving base points. The map~(\ref{eqn-instanceDWWprin2}) is a $\Gamma$-map and we
obtain therefore an induced map of reduced Borel constructions
\begin{equation} \label{eqn-instanceDWWprin3} \Str(M\times\RR^i\cc)_{rh\Gamma} \lra F\lpr(M\times\RR^i\cc)_{rh\Gamma}~. \end{equation}
Now $\Str(M\times\RR^i\cc)_{rh\Gamma}$ \emph{is} a good combinatorial model for
$\sS^\red(M\times\RR^i\cc)$ in the sense that the forgetful map
\begin{equation} \label{eqn-approxcontrol}
\Str(M\times\RR^i\cc)_{h\Gamma} \lra \sS^\red(M\times\RR^i\cc)_{rh\Gamma}
\end{equation}
is a homology equivalence. See appendix~\ref{sec-altapprox}. Therefore the homotopy pushout of
\[
\xymatrix{
\sS^\red(M\times\RR^i\cc)_{rh\Gamma}  & \ar[l] \Str(M\times\RR^i\cc)_{rh\Gamma} \ar[r] & F\lpr(M\times\RR^i\cc)_{rh\Gamma}
}
\]
is homotopy equivalent to $F\lpr(M\times\RR^i\cc)$ by an obvious inclusion. Another obvious inclusion,
that of $\sS^\red(M\times\RR^i\cc)$ into the same homotopy pushout, can therefore be described loosely in the form
\begin{equation} \label{eqn-thatmapcontrol}
\sS^\red(M\times\RR^i\cc) \lra F\lpr(M\times\RR^i\cc)~. \end{equation}
As in the case $i=0$, the construction can be extended to the case where
$M$ is compact with nonempty boundary, and some relative cases. The details of this are left to the reader.

\section{Algebraic models for structure spaces: Proofs}
\label{sec-approxproofs}
\proof[Proof of lemma~\ref{lem-oldhofiber}.]
There is the following commutative square
of controlled structure spaces:
\begin{equation} \label{eqn-AndersonHsiang}  
\begin{split}
\xymatrix{
\sS(M\times\RR^i\times\,]-\infty,1],\,M\times\RR^i\times 1\cc) \ar[r] & \sS(M\times\RR^{i+1}\cc) \\
\sS(M\times\RR^i\times[0,1],\,M\times\RR^i\times 1\cc) \ar[r] \ar[u] & \ar[u] \sS(M\times\RR^i\times[0,\infty[\,\cc).
}
\end{split}
\end{equation}
\emph{Details on the spaces involved}: All manifolds which appear in the diagram,
such as $M\times\RR^i\times[0,1]$ etc., are codimension zero submanifolds
of $M\times\RR^{i+1}$. We make them into control spaces by taking
closures inside $M*S^i$. This amounts to adding the equator $S^{i-1}$ of $S^i$
to $M\times\RR^i\times[0,1]$, the closed upper hemisphere of $S^i$ to $M\times\RR^i\times[0,\infty[$
and the closed lower hemisphere of $S^i$ to $M\times\RR^i\times\,]-\infty,1]$.
Beware that structure space notation of the form $\sS(C,D)$ usually assumes that $C$ is a manifold,
$D$ is a codimension zero closed submanifold of $\partial C$ and the
structures considered are trivial over the closure of the
complement of $D$ in $\partial C$. Hence there is a homotopy
fiber sequence $\sS(C)\to \sS(C,D)\to \sS(D)$. \newline
\emph{Details on the maps}: The maps in the square are given by obvious
extension of structures in all cases. For the left-hand vertical map
for example, we extend by gluing with the trivial structure (identity)
on $M\times\RR^i\times\,]-\infty,0]$. For the lower horizontal
map, we extend a structure $f$ by gluing with $\partial f\times\id$,
where $\partial f$ is the boundary structure (on $M\times\RR^i\times 1$)
determined by $f$ and $\id$ means the identity map on $[1,\infty[\,$. \newline
In the square~(\ref{eqn-AndersonHsiang}), the lower left-hand term is clearly
homotopy equivalent to the controlled $h$-cobordism space
$\sH(M\times\RR^i\cc)$. The upper left-hand term
is homotopy equivalent to $\sS(M\times\RR^i\cc)$, via restriction
of structures to the boundary. With this identification, the upper
horizontal arrow is just $\times\RR$ up to homotopy. The lower right-hand
term is contractible.
The square is homotopy cartesian. (Modulo a replacement of controlled
structure spaces by the homotopy equivalent bounded structure
spaces, this goes back to Anderson and Hsiang \cite{AndersonHsiang}.
See also \cite{WW1} for some added details.) Therefore the
square amounts to a homotopy fibration sequence
\[
\CD
\sH(M\times\RR^i\cc) @>>> \sS(M\times\RR^i\cc) @>{\times\RR}>> \sS(M\times\RR^{i+1}\cc).
\endCD
\]
This completes the proof in the case where $f\co N\to M$ is the identity of $M$. The general case
follows by a translation argument: there is a commutative diagram
\begin{equation} \label{eqn-famoustranslate}
\begin{split}
\xymatrix{
\sS(N\times\RR^{i+1}\cc) \ar[r]^-{\simeq} & \sS(M\times\RR^{i+1}\cc)  \\
\sS(N\times\RR^i\cc) \ar[r]^-{\simeq} \ar[u]^-{\times\RR} & \ar[u]^-{\times\RR} \sS(M\times\RR^i\cc)
}
\end{split}
\end{equation}
where the horizontal arrows are given by composition with $f\times\RR^i$ and $f\times\RR^{i+1}$, respectively. \qed

\medskip
\proof[Proof of lemma~\ref{lem-newhofiber}.] This is easy from definitions~\ref{defn-LK}
and~\ref{defn-LKcontrol}, proposition~\ref{prop-timesRandXi} and corollary~\ref{cor-ultimateAstab}.
We allow ourselves to use twisted versions of proposition~\ref{prop-timesRandXi}
and corollary~\ref{cor-ultimateAstab}, in the spirit of section~\ref{sec-twisted}; the twist is determined
by $\nu$, the normal bundle of $M$ (not shown in the notation by our current conventions).
From the first homotopy cartesian square in~\ref{prop-timesRandXi}, twisted version, we get that the square
\[
\xymatrix{
\Omega^{m+i}\bLA\lbul(M\times\RR^i,m+i,\cc) \ar[r]^-{\textup{forget}} \ar[d]^-{\times\RR} &
\ar[d]^-{\times\RR} \bA(M\times\RR^i,m+i\cc)\hf \\
\Omega^{m+i+1}\bLA\lbul(M\times\RR^{i+1},m+i+1\cc)
\ar[r]^-{\textup{forget}} & \bA(M\times\RR^{i+1},m+i+1\cc)\hf
}
\]
is also homotopy cartesian. (Compare the horizontal homotopy fibers.) Hence the homotopy fiber of the
left-hand column, which we are interested in, is identified with the homotopy fiber of the right-hand column.
By corollary~\ref{cor-ultimateAstab}, twisted version, this can be identified with the homotopy fiber of
\[
\CD
\Omega\bA^h(M\times\RR^{i+1},m+i+1\cc)\hf \\
@VV\textup{ev} V \\
\bA(M\times\RR^{i+1},m+i+1\cc)\hf
\endCD
\]
where $\ZZ/2$ acts by flip on the $\Omega$ and the map $\textup{ev}$ is given by evaluation at the center
of a loop. Applying $\Omega^\infty$ to that homotopy fiber, we get
\[  \Omega^{\infty+1}\bA^h(M\times\RR^{i+1}\cc)~\simeq~\Omega^\infty\bA(M\times\RR^i\cc) \]
(using corollary~\ref{cor-ultimateAstab} again). Summarizing, we have established a
homotopy fiber sequence of infinite loop spaces
\[
\xymatrix@R=14pt{
\Omega^\infty\bA(M\times\RR^i\cc) \ar[d] \\
\Omega^{\infty+m+i}\bLA\lbul(M\times\RR^i,m+i,\cc) \ar[d]^{\times\RR} \\
\Omega^{\infty+m+i+1}\bLA\lbul(M\times\RR^{i+1},m+i+1\cc)
}
\]
The entire argument is compatible with assembly; i.e., we could instead have started with a homotopy
cartesian square like the first one in~\ref{prop-timesRandXi}, but with $\bLA\upr$ and $\bA\upr$
or $\bLA\lpr$ and $\bA\lpr$ instead of $\bLA$ and $\bA$. So lemma~\ref{lem-newhofiber} is established.
\qed

\bigskip
We turn to proposition~\ref{prop-whatvertical}.
A simplification which comes from the proof of lemma~\ref{lem-newhofiber}
is that instead of working with the commutative square~(\ref{eqn-mainmapsquare}),
we may work with the simplified version
\begin{equation} \label{eqn-mainmapsquarecut}
\CD
\sS^\red(M\times\RR^{i+1}\cc) @>>> \Omega^{\infty}\bA\lpr(M\times\RR^{i+1},m+i+1\cc)\hf \\
@AAA @AA {\times\RR} A \\
\sS^\red(M\times\RR^i\cc) @>>> \Omega^{\infty}\bA\lpr(M\times\RR^i,m+i\cc)\hf.
\endCD
\end{equation}
The horizontal arrows in diagram~(\ref{eqn-mainmapsquarecut}) are obtained from those
in diagram~(\ref{eqn-mainmapsquare}) by composing with forgetful maps. This can be simplified
some more by concatenating with the homotopy cartesian square of corollary~\ref{cor-moreAstab}.
We also use abbreviations $A$, $A\upr$ and $A\lpr$ for $\Omega^\infty\bA$, $\Omega^\infty\bA\upr$
and $\Omega^\infty\bA\lpr$, respectively, where possible.
The result is, after a rotation,
\begin{equation} \label{eqn-mainmapsquare2cut}
\begin{split}
\xymatrix@C=40pt@M=7pt{
\sS^\red(M\times\RR^i\cc) \ar[r]^-{\times\RR} \ar[d]^-{j\chi\lpr} & \sS^\red(M\times\RR^{i+1}\cc) \ar[d]^-{\chi\lpr} \\
\cone~A\lpr(M\times\RR^i\cc) \ar[r]^-{\times\RR}      & A\lpr(M\times\RR^{i+1}\cc) \\
}
\end{split}
\end{equation}
The labels $\chi\lpr$ refer to constructions which extract
controlled $A$-theory characteristics (playing off the excisive type $\chi\upr$ against
the controlled homotopy invariant type $\chi$, and using discrete approximation technology for that). The left-hand column
is a map of type $\chi\lpr$ followed by the inclusion
\[ j\co A\lpr(M\times\RR^i\cc) \lra \cone~A\lpr(M\times\RR^i\cc). \]
Therefore $\Phi_x$~, $\Psi_y$~, $\Phi_{x,a}$ and $\Psi_{y,b}$ in the statement of proposition~\ref{prop-whatvertical}
can be redefined in terms of diagram~(\ref{eqn-mainmapsquare2cut}), now as horizontal homotopy fibers or components of such.
Another simplification
comes from diagram~(\ref{eqn-famoustranslate}): without loss of generality, $x$ is the base point of
$\sS(M\times\RR^{i+1}\cc)$, so that $y$ is the base point of $A\lpr(M\times\RR^{i+1}\cc)$.
Yet another small simplification: without loss of generality, $M$ is connected and based.

\proof[Proof of proposition~\ref{prop-whatvertical}, first part] We pass to
horizontal homotopy fibers $\Phi_\pt$ and $\Psi_\pt$ in diagram~(\ref{eqn-mainmapsquare2cut}).
Commutativity of the diagram leads to a map
\begin{equation} \label{eqn-thebigmap} \Phi_\pt\to \Psi_\pt~. \end{equation}
We have already identified $\Phi_\pt$ with a union of connected components of the
controlled $h$-cobordism space $\sH(M\times\RR^i\cc)$,
and $\Psi_\pt$ of~(\ref{eqn-thebigmap}) is clearly
homotopy equivalent to $\Omega A\lpr(M\times\RR^{i+1}\cc)$, and also to
$A\lpr(M\times\RR^i\cc)$. After deleting some more connected components of the source, if necessary,
we can therefore write~(\ref{eqn-thebigmap}) in the form
\begin{equation} \label{eqn-therestrictedbigmap}
\xymatrix{ \sH^\red(M\times\RR^i\cc) \ar[r] & A\lpr(M\times\RR^i\cc)~.
}
\end{equation}
Here $\pi_0$ of $\sH^\red(M\times\RR^i\cc)$ can be identified with
the image of $\pi_0$ of $\sH(M)$ in $\pi_0$ of $\sH(M\times\RR^i\cc)$.
The aim is to show that the map~(\ref{eqn-therestrictedbigmap}), \emph{constructed as a
restriction of}~(\ref{eqn-thebigmap}), is identical with a map
\begin{equation} \label{eqn-familiarWald}
\chi^\lpr\co \sH^\red(M\times\RR^i\cc) \lra  A\lpr(M\times\RR^i\cc).
\end{equation}
which merits the label $\chi^\lpr$ because it is obtained by playing off
excisive characteristics against (controlled) homotopy invariant characteristics.
We also want to show that~(\ref{eqn-thebigmap}) itself induces an injection on $\pi_0$. To clarify this
last statement, we note that $\pi_0$ of $\Psi_\pt\simeq A\lpr(M\times\RR^i\cc)$ is
\[
\begin{array}{rl} \Wh(\pi_1M) & \textup{ if }i=0 \\
\tilde K_0(\ZZ\pi_1M) & \textup{ if }i=1 \\
K_{1-i}(\ZZ\pi_1M) & \textup{ if }i>1~.
\end{array}
\]
We aim to show that the map
$\pi_0\Phi_\pt\to \pi_0\Psi_\pt$
induced by~(\ref{eqn-thebigmap}) associates to
a controlled $h$-cobordism its controlled Whitehead torsion. This implies that it is
injective, by the classification of controlled $h$-cobordisms. \newline
Here the proof proper begins.
In order to make a closer connection with diagram~(\ref{eqn-AndersonHsiang}) we enlarge the left-hand column
in diagram~(\ref{eqn-mainmapsquare2cut}), without changing homotopy types, to get
\begin{equation} \label{eqn-mainmapsquare3cut}
\begin{split}
\xymatrix@R=40pt{
{\sS^\red\left(M\times\RR^i\times\,]-\infty,1],M\times\RR^i\times 1\cc\right)}
\ar[r] \ar[d]^-{j\chi\lpr} &  \sS^\red(M\times\RR^{i+1}\cc) \ar[d]^-{\chi\lpr}  \\
{\Omega^\infty\left(\bA\lpr\left(\begin{matrix} M\times\RR^i\times 1 \\ \downarrow \\ M\times\RR^i\times\,]-\infty,1]\end{matrix}
\cc\right)\cup \bC\right)} \ar[r] & {A\lpr(M\times\RR^{i+1}\cc)}
}
\end{split}
\end{equation}
Here $\bC:=\cone~\bA\lpr(M\times\RR^i\cc)$ and the union $\cup\bC$ is taken along the common subspectrum
$\bA\lpr(M\times\RR^i\cc)$. We use the embedding
\[ \bA\lpr(M\times\RR^i\cc) \lra  \bA\lpr\left(\begin{matrix} M\times\RR^i\times 1 \\
\downarrow \\ M\times\RR^i\times\,]-\infty,1]\end{matrix}
\cc\right) \]
induced by exact functors of type \emph{product with the cofibration $\{1\}\to\, ]-\infty,1]$}.
The top row in~(\ref{eqn-mainmapsquare3cut}) is already familiar
as the top row of diagram~(\ref{eqn-AndersonHsiang}) and the lower row
is an $A$-theory counterpart. (A small complication here: the map
\[
\xymatrix@M=4pt{ {\sS^\red\left(M\times\RR^i\times\,]-\infty,1],M\times\RR^i\times 1\cc\right)} \ar[d]^-{\chi\lpr} \\
{A\lpr\left(\begin{matrix} M\times\RR^i\times 1 \\ \downarrow \\ M\times\RR^i\times\,]-\infty,1]\end{matrix}
\cc\right)}
}
\]
implicit in the left-hand column
of~(\ref{eqn-mainmapsquare3cut}) needs to be defined in a roundabout way as in~(\ref{eqn-relphi1}) or~(\ref{eqn-relphi2}),
using $A\lpr$ instead of $F\lpr$.) In order to make the connection with diagram~(\ref{eqn-AndersonHsiang}) closer still,
we think of~(\ref{eqn-mainmapsquare3cut}) as one face of a cubical diagram whose
opposite face is
\begin{equation} \label{eqn-mainmapsquarecomplement}
\begin{split}
\xymatrix@R=40pt{
{\sS^\red\left(M\times\RR^i\times[0,1],M\times\RR^i\times 1 \cc\right)}
\ar[r] \ar[d]^-{j\chi\lpr}
 & \sS^\red(M\times\RR^i\times[0,\infty[\,\cc) \ar[d]^-{\chi\lpr} \\
{\Omega^\infty\left(\bA\lpr\left(\begin{matrix} M\times\RR^i\times 1 \\ \downarrow \\ M\times\RR^i\times[0,1] \end{matrix}
\cc\right)\cup \bC\right)} \ar[r] & {A\lpr(M\times\RR^i\times[0,\infty[\,\cc)}
}
\end{split}
\end{equation}
The upper row in~(\ref{eqn-mainmapsquarecomplement})
is already familiar as the lower row of diagram~(\ref{eqn-AndersonHsiang}) and the lower row
in~(\ref{eqn-mainmapsquarecomplement}) is an $A$-theory counterpart of the upper row.
The promised cube is a map from~(\ref{eqn-mainmapsquarecomplement})
to~(\ref{eqn-mainmapsquare3cut}). The cube therefore has a geometric face
\begin{equation} \label{eqn-AndersonHsiangred}
\begin{split}
\xymatrix{
\sS^\red(M\times\RR^i\times\,]-\infty,1],\,M\times\RR^i\times 1\cc) \ar[r] & \sS^\red(M\times\RR^{i+1}\cc) \\
\sS^\red(M\times\RR^i\times[0,1],\,M\times\RR^i\times 1\cc) \ar[r] \ar[u] & \ar[u] \sS(M\times\RR^i\times[0,\infty[\,\cc)
}
\end{split}
\end{equation}
which is a sub-square of~(\ref{eqn-AndersonHsiang}), obtained by selecting some connected components in each
term of~(\ref{eqn-AndersonHsiang}). Therefore by lemma~\ref{lem-oldhofiber} it is almost a
homotopy pullback square. More precisely, the map
induced from the initial (lower left-hand) term to
the homotopy pullback of the other three terms has all homotopy fibers empty or contractible.
(It might not induce a surjection on $\pi_0$.) The cube also has an $A$-theory face
\begin{equation} \label{eqn-cubeAface}
\begin{split}
\xymatrix@C=15pt{
{\Omega^\infty\left(\bA\lpr\left(\begin{matrix} M\times\RR^i\times 1 \\ \downarrow \\ M\times\RR^i\times\,]-\infty,1]\end{matrix}
\cc\right)\cup \bC\right)}
\ar[r]  &   A\lpr(M\times\RR^{i+1}\,\cc)  \\
{\Omega^\infty\left(\bA\lpr\left(\begin{matrix} M\times\RR^i\times 1 \\ \downarrow \\ M\times\RR^i\times[0,1] \end{matrix}
\cc\right)\cup \bC\right) } \ar[r] \ar[u] &  {A\lpr(M\times\RR^i\times[0,\infty[\,\cc)} \ar[u]
}
\end{split}
\end{equation}
and this is a homotopy pullback square by inspection (using the additivity theorem). The conclusion from
this diagrammatic reasoning is that the map~(\ref{eqn-thebigmap}) agrees with
\begin{equation} \label{eqn-cubeedge}
\begin{split}
\xymatrix{
{\sS^\red(M\times\RR^i\times[0,1],\,M\times\RR^i\times 1\cc)}
\ar[d]^-{j\chi\lpr} \\
{\Omega^\infty\left(\bA\lpr\left(\begin{matrix} M\times\RR^i\times 1 \\ \downarrow \\ M\times\RR^i\times[0,1] \end{matrix}
\cc\right)\cup\bC\right)}
}
\end{split}
\end{equation}
from the cube, after appropriate restriction. It is clear that
\[  \sS^\red(M\times\RR^i\times[0,1],\,M\times\RR^i\times 1\cc)~\simeq~\sH^\red(M\times\RR^i\cc).
\]
The canonical inclusion
\[ \bA\lpr(M\times\RR^i\times[0,1]\cc) \lra \bA\lpr\left(M\times\RR^i\times 1\to M\times\RR^i\times[0,1]\cc\right)\cup\bC \]
is a homotopy equivalence, and $A\lpr(M\times\RR^i\times[0,1]\cc)\simeq A\lpr(M\times\RR^i\cc)$.
With some elementary simplifications we are therefore allowed to write~(\ref{eqn-cubeedge}), alias
restriction of~(\ref{eqn-thebigmap}), in
the form~(\ref{eqn-familiarWald}). This was one of our declared goals.
Our construction or description of~(\ref{eqn-familiarWald}) deviates slightly from the standard because we are in
a controlled setting, and because we have used discrete approximations of homeomorphism groups of manifolds
as in \cite{McDuff80} instead of the technology of simplicial sets and simple maps \cite{WaldRutgers}, \cite{WaldJaRo}. \newline
It remains to prove that the map $\pi_0\Phi_\pt\to \pi_0\Psi_\pt$ induced by~(\ref{eqn-thebigmap}) is the
controlled Whitehead torsion map. This comes from a mild improvement on the diagrammatic reasoning just
developed. Let us admit that by using
\[ \sS^\red(M\times\RR^i\times[0,1],\,M\times\RR^i\times 1\cc) \]
as the initial term in square~(\ref{eqn-AndersonHsiangred}), we did not make the best choice. Instead we could have used
\[  T\cup \sS^\red(\sS^\red(M\times\RR^i\times[0,1],\,M\times\RR^i\times 1\cc) \]
where $T$ is any subset of $\sS(M\times\RR^i\times[0,1],\,M\times\RR^i\times 1\cc)$ which has exactly one
element in each connected component mapping to $\sS^\red(M\times\RR^i\cc)\subset \sS(M\times\RR^i\cc)$ under
the forgetful map (and no other elements). If we enlarge~(\ref{eqn-AndersonHsiangred}) in this manner, we have
a slightly better approximation to a homotopy pullback square: the canonical map from initial term to the
homotopy pullback of the other three terms is still a homotopy equivalence on base point components,
but it also induces a bijection on $\pi_0$.
It is straightforward to extend the map from~(\ref{eqn-AndersonHsiangred}) to~(\ref{eqn-cubeAface}) which we already have to this enlargement of~(\ref{eqn-AndersonHsiangred}), using $\chi\lpr$ constructions as before.
This allows us to see what~(\ref{eqn-thebigmap}) does on $\pi_0$.
\qed

\bigskip
\proof[Proof of proposition~\ref{prop-whatvertical}, second part] Here the goal is to show that~(\ref{eqn-thebigmap})
restricts to a $(k_N+i+1)$-connected map
\[  \Phi_{\pt,a}\to \Psi_{\pt,b} \]
in the notation of proposition~\ref{prop-whatvertical}. Note that $N$ depends on $a$. There is an easy
reduction to the case where $a$ is the base point component, $a=0$. This is again a translation argument.
The component $a$ of the homotopy fiber $\Phi_\pt$
is homotopy equivalent to the base point component of
\[  \hofiber_\pt[\,\sS(N\times\RR^i\cc) \to \sS(N\times\RR^{i+1}\cc)\,]~. \]
Therefore replacing $M$ by $N$ will do the trick. The details are left to the reader. \newline
We are now looking at $\Phi_{\pt,0}\to \Psi_{\pt,0}$ and we may think of that as the restriction
of~(\ref{eqn-familiarWald}) to the base point components. Therefore it only remains to show
that the restriction of~(\ref{eqn-familiarWald}) to base point components is $(k_M+i+1)$-connected.
We use induction on $i$.
For $i=0$ we know that~(\ref{eqn-familiarWald}) is Waldhausen's map and so is $(k_M+1)$-connected after restriction
to base point components; see beginning of section 10 in \cite{DwyerWeissWilliams}. For $i>0$, there is a
homotopy pullback square of controlled $h$-cobordism spaces and inclusion-induced maps,
\begin{equation} \label{eqn-ctrlindu1}
\begin{split}
\xymatrix{
\sH(M\times\RR^{i-1}\times\,]-\infty,1]\cc) \ar[r] & \sH^\red(M\times\RR^i\cc) \\
\sH(M\times\RR^{i-1}\times[-1,1]\cc) \ar[r] \ar[u] & \ar[u] \sH(M\times\RR^{i-1}\times[-1,\infty[\,\cc).
}
\end{split}
\end{equation}
This is closely related to the homotopy pullback square
\begin{equation} \label{eqn-ctrlindu2}
\begin{split}
\xymatrix{
A\lpr(M\times\RR^{i-1}\times\,]-\infty,1]\cc) \ar[r] & A\lpr(M\times\RR^i\cc) \\
A\lpr(M\times\RR^{i-1}\times[-1,1]\cc) \ar[r] \ar[u] & \ar[u] A\lpr(M\times\RR^{i-1}\times[-1,\infty[\,\cc).
}
\end{split}
\end{equation}
which we have from theorem~\ref{thm-changeofdeco}. In both squares, the off-diagonal terms are contractible.
Using discrete approximation of homeomorphism groups as in sections~\ref{sec-discrete} and~\ref{sec-altapprox},
we can construct a map of type $\chi\lpr$ from a mildly
corrupted version of~(\ref{eqn-ctrlindu1}) to~(\ref{eqn-ctrlindu2}).
For our purposes, a sufficiently mild form of
corruption is to preserve all terms in~(\ref{eqn-ctrlindu1}) except the initial term, to
preserve the base component in the initial term as well, and to replace each other component of the initial term
by a selected point in it. The existence of such a map implies that if
\[
\xymatrix@M=5pt{
\sH^\red(M\times[-1,1]\times\RR^{i-1}\cc) \ar[r]^-{\chi\lpr} & A\lpr(M\times[-1,1]\times\RR^{i-1}\cc)
}
\]
is $(k_M+i)$-connected on base point components, then the looping of~(\ref{eqn-familiarWald}) induces a bijection
on $\pi_0$ and is also $(k_M+i)$-connected on base points components. The looping of~(\ref{eqn-familiarWald}) is
therefore unconditionally $(k_M+i)$-connected, and so~(\ref{eqn-familiarWald}) is $(k_M+i+1)$-connected on base
point components. This is our induction step: a reduction from $i$ to $i-1$ for the price of replacing
$M$ by $M\times[-1,1]$.
\qed

\medskip
\begin{rem}
\label{rem-CSWQR}
\label{rem} {\rm Although formula~(\ref{eqn-CSWQRcontrol})
is meant as a quotation, a historical essay around it is in order, with a view to
the proof of lemma~\ref{lem-inthelimit}. Perhaps the best known, most polished and least complicated method for
establishing a homotopy fiber sequence
\[ \widetilde{\sS}^s(M\times\RR^i\cc) \lra
\Omega^{\infty+m+i}\bL^s\lbul\lpr(M\times\RR^i\cc)\lra L_0(\ZZ)^{\pi_0M} \]
is the one developed by Ranicki \cite[\S18]{RanickiTopMan}. We summarize the key
points, referring to \cite{RanickiTopMan} for definitions and clarifications.
\begin{enumerate}
\item The cases $i>0$ can be reduced to the case $i=0$ by a torus trick.
\item A compact polyhedron $X$ and a homotopy equivalence $e\co M\to X$ transverse
to the triangulation of $X$ are chosen. A degree one normal map $f\co N\to M$ between
closed manifolds determines (if $f\circ e$ is transverse to the triangulation of $X$)
chain complexes $C(N)$ and $C(M)$ ``dissected'' over $X$, with dissected
nondegenerate symmetric structures of formal dimension $m$. The map
$f_*\co C(N)\to C(M)$ respects the symmetric structures and so determines a splitting
\[   C(N) \simeq D\oplus C(M) \]
where the ``kernel'' $D$ is again dissected and comes with a nondegenerate symmetric structure.
But $D$ is globally (after assembly) contractible. Hence the dissected symmetric
structure has an automatic refinement to a dissected quadratic structure.
\item $\bL\lbul\upr(M)$ admits a description as the bordism theory of chain complexes
dissected over $X$, with a dissected quadratic Poincar\ee structure.
With that description, the assembly map is induced by the assembly (``universal'' assembly)
of dissected chain complexes.
\item The local degree homomorphism
from $\pi_m\bL\lbul\upr(M)$ to $L_0(\ZZ)^{\pi_0M}$ is zero on elements determined by
degree one normal maps $f\co N\to M$. (The reason is
that an oriented map $f$ between manifolds of the same dimension which has ``global''
degree one, in the sense that it respects fundamental classes, will also have local degree
one, i.e., the generic cardinality of $f^{-1}(x)$ is 1 for any $x$ in the target. It is this
relationship between global and local degree which has been shown to fail spectacularly
\cite{BryantFerryMioWeinberger96} in the world of ANR homology manifolds.)
\end{enumerate}
These ideas, generalized to a parameterized setting (where the parameter
space is $\Delta^k$ for $k=0,1,2,\dots)$ lead to a map
\[ \widetilde{\sS}^s(M\times\RR^i\cc) \lra
\Omega^{\infty+m+i}\bL^s\lbul\lpr(M\times\RR^i\cc).
\]
Proving that the map
is a homotopy equivalence except for a deviation in $\pi_0$ is another matter
and there is no need to go into that here. What we need to do is this:
unravel each of the items (1)--(3) and relate it to the methods (e.g.
the \emph{characteristic element} method and the \emph{control method}) which we have
favored in this paper.

(1a). The torus trick consists in using that $\sS^s(M\times\RR^i\cc)$
and $\bL^s\lbul\lpr(M\times\RR^i\cc)$ are homotopy retracts of
$\sS^s(M\times(S^1)^i)$ and $\bL^s\lbul\lpr(M\times (S^1)^i)$, respectively.
The retracting maps
\[
\begin{array}{c}
\sS^s(M\times(S^1)^i)\lra \sS^s(M\times\RR^i\cc), \\
\bL^s\lbul\lpr(M\times (S^1)^i)\lra \bL^s\lbul\lpr(M\times\RR^i\cc)
\end{array}
\]
are obvious transfer maps in both cases.

(2a). Dissection theory works with retractive spectra over $X$
and visible symmetric structures just as well as with chain complexes and
symmetric structures. It remains true that a dissected visible
symmetric structure on a dissected retractive spectrum $Y$ over $X$ (subject to some
finiteness conditions) automatically lifts to a dissected quadratic structure if
$Y$ is globally weakly equivalent to zero.

(3a). The interpretation of the assembly map in terms of \emph{assembly of
dissected chain complexes} works just as well with dissected retractive
spaces and spectra. Given a degree one normal map
$f\co N\to M$ which is a homotopy equivalence, there are two slightly different
ways of using the dissection argument to extract $L$-theoretic information.
\begin{itemize}
\item One way is to form the \emph{visible symmetric kernel} (mapping cone of the stable
Umkehr map $M\amalg X\to N\amalg X$, with a
nondegenerate visible symmetric structure). It is dissected over $X$ along
with its nondegenerate visible symmetric structure. Since $f$ is a homotopy
equivalence, the kernel is globally contractible and
the dissected nondegenerate visible symmetric structure on it automatically
lifts to a dissected nondegenerate quadratic structure. The global contractibility
of the kernel implies a preferred global nullbordism of the kernel (i.e., a nullbordism
of the assembled kernel with the assembled quadratic structure).
The dissected kernel and the global nullbordism then determine
a point in the homotopy fiber (over the base point, zero) of the assembly map
\[ \Omega^{\infty+m}\bL\lbul\upr(X) \lra \Omega^{\infty+m}\bL\lbul(X)\,. \]
\item Another way is to note that $M\amalg X$ and $N\amalg X$ can themselves
be dissected over $X$, and come with dissected nondegenerate visible symmetric
structures. Before assembly, these may not be equivalent; after assembly
they are certainly equivalent via $f\amalg \id: N\amalg X\to M\amalg X$. Hence these data
determine a point in a homotopy fiber of the assembly map
\[ \Omega^{\infty+m}\bV\bL\ubul\upr(X) \lra \Omega^{\infty+m}\bV\bL\ubul(X). \]
This time we take the homotopy fiber over the point in
$\Omega^{\infty+m}\bV\bL\ubul(X)$ determined by $M\amalg X$ with its (assembled)
nondegenerate visible symmetric structure.
\end{itemize}
There is an obvious compatibility between the two methods. By
theorem~\ref{thm-VhatLexcision}, we lose no $L$-theoretic information by
relying on the second method. \newline
But now another little problem remains. We have two descriptions of the
assembly map in $L$-theory, one in terms of dissections and another
one using control. How are they related? An easy way to make a connection is
to use both approaches simultaneously. Very schematically, we have one description
of $\bV\bL\ubul\upr$-theory as
\[   \bV\bL\ubul(X,\textup{dissected}) \]
and another as the homotopy fiber of an inclusion
\[ \bV\bL\ubul(X) \lra \bV\bL\ubul(X\times[0,1[\,,\textup{controlled}), \]
as in section~\ref{sec-excicharsig}. There is a third
description of $\bV\bL\ubul\upr(X)$ as the homotopy fiber
of an inclusion
\[  \bV\bL\ubul(X,\textup{dissected}) \lra
\bV\bL\ubul(X\times[0,1[\,,\textup{controlled and dissected})\,. \]
(The dissections are always over $X$, even where we are dealing with retractive
spaces over $X\times[0,1[$. \emph{Control} refers to the control space $\JJ X$.)
Validating this third formula amounts to showing that
\[ \bV\bL\ubul(X\times[0,1[\,,\textup{controlled and dissected}) \]
is contractible. (This is left to the reader.) This makes the connection which we were
after. Moreover it does that in such a way that the two standard methods (by dissection
and by control) of lifting manifold signatures in $VL$-spaces
across the assembly map are seen to agree.
}
\end{rem}

\medskip
\proof[Proof of lemma~\ref{lem-inthelimit}.] The homotopy colimit of
\[ \bA\lpr(M)\stackrel{\times\RR\rule{2.5mm}{0mm}}{\lra}\bA\lpr(M\times\RR\cc)
\stackrel{\times\RR\rule{2.5mm}{0mm}}{\lra}
\bA\lpr(M\times\RR^2\cc)\to\cdots \]
is contractible by lemma~\ref{lem-degenerateprod}. It follows that
\[ \hocolim_i \,\,(\bA\lpr(M\times\RR^i,m+i\cc))\ho \]
is also contractible. (The dimension indicator $m+i$ here specifies the involution or
the $SW$-product, which is obtained from the standard one by $(m+i)$-fold looping).
Therefore the inclusion
\[ \hocolim_i (\bA\lpr(M\times\RR^i,m+i\cc))\hf
\lra \hocolim_i (\bA\lpr(M\times\RR^i,m+i\cc))\thf \]
is a homotopy equivalence. Therefore the forgetful map
\[   \hocolim_i\,\, F\lpr(M\times\RR^i\cc)\lra
\hocolim_i\,\,\Omega^{\infty+m+i}\bL\lbul\lpr(M\times\RR^i\cc) \]
is also a homotopy equivalence. Hence in the limit $i\to \infty$, it does not matter
whether we use the map $\varphi$
in construction~\ref{con-extendedmainmap} or a simplified version
\begin{equation}  \label{eqn-notilde1}
\sS^\red(M\times\RR^i\cc) \lra \Omega^{\infty+m+i}\bL\lbul\lpr(M\times\RR^i\cc)
\end{equation}
constructed purely in terms of $VL$-theoretic signatures (as in sections~\ref{sec-hocharsig}
and \ref{sec-excicharsig}, but without any algebraic $K$-theory). When $i>0$, this map
has an automatic lift
\begin{equation}  \label{eqn-notilde2}
\sS^\red(M\times\RR^i\cc) \lra \Omega^{\infty+m+i}\bL^{cs}\lbul\lpr(M\times\RR^i\cc)  
\end{equation}
because all components of $\sS^\red(M\times\RR^i\cc)$ are represented by (controlled) simple structures.
Furthermore, it is easy to extend~(\ref{eqn-notilde2}) to the block structure space
\[ \widetilde\sS^\red(M\times\RR^i\cc)\,. \]
The reason is that we can extend~(\ref{eqn-notilde2})
by adding on another simplicial direction. In other words, there are maps
\begin{equation} \label{eqn-yestilde}
\sS^\red(M\times\Delta^k_!\times\RR^i\cc)
\lra \Omega^{\infty+m+i}\bL^{cs}\lbul\lpr(M\times\Delta^k_!\times\RR^i\cc) 
\end{equation}
for every $k\ge 0$, generalizing~~(\ref{eqn-notilde2}); the notation $\Delta^k_!$ means that we
think of $\Delta^k$ as a diagram of manifolds (the faces of any codimension),
not as a single manifold with boundary.
By taking geometric realizations over $k$ (and noting, as we have done before, that on the target
side all the face operators are homotopy equivalences), we obtain a single map
\[
\CD
\widetilde\sS^\red(M\times\RR^i\cc)
@>>>
\Omega^{\infty+m+i}\bL\lbul\lpr(M\times\Delta^k_!\times\RR^i\cc)\,.
\endCD
\]
Now we let the new simplicial direction take over by restricting~(\ref{eqn-yestilde}) to
appropriate $0$-skeletons for each $k$, and noting that that does not affect the
homotopy type of the geometric realization over $k$. Then what we have is just the standard map
of~(\ref{eqn-CSWQRcontrol}), restricted to a union of connected components, namely $\sS^\red\subset \sS^{cs}$.
This identification relies on remark~\ref{rem-CSWQR}. \qed

\medskip
\proof[Proof of lemma~\ref{lem-maincomponents}.] Without loss of generality, $M$ is connected.
There is a commutative square
\begin{equation} \label{eqn-maincomponents}
\begin{split}
\xymatrix@C=30pt{
& \pi_0(\bA^h\lpr(M,\nu,m)\ho) \ar[d] \\
\pi_0\sS(M) \ar[r] \ar[d] & \pi_m\bLA\lbul\lpr(M,\nu,m) \ar[d] \\
\pi_0\widetilde\sS^h(M) \ar[r]  &  \ar[d] \pi_m\bL^h\lbul\lpr(M,w_\nu) \ar[r]^-{\textup{loc.\,deg.}} & \ZZ \\
& 0  }
\end{split}
\end{equation}
with exact lower row and exact middle column. The middle column is part of the long exact sequence
in homotopy groups of a homotopy fiber sequence
\[
\xymatrix@M=5pt@C=30pt{
\Omega^{m}\bLA\lbul\lpr(M,\nu,m) \ar[r] & \Omega^{m}\bL^h\lbul\lpr(M,w_\nu) \ar[r]^-{(\ref{eqn-defLA})} &
S^1\wedge \bA\lpr^h(M,\nu,m)\ho~.
}
\]
We need to understand the composition
\[
\xymatrix@M=5pt@C=30pt{
\widetilde\sS^h(M) \ar[r]^-{(\ref{eqn-surgcalc})} & \Omega^{\infty+m}\bL^h\lbul\lpr(M,w_\nu) \ar[r]^-{(\ref{eqn-defLA})}
& \Omega^{\infty-1}(\bA\lpr^h(M,\nu,m)\ho)
}
\]
on the $1$-skeleton, as a map taking the subspace $\sS(M)$ to the base point.
For this purpose we may replace the right-hand term by the space $B(\ZZ\otimes_{\ZZ}\Wh)$, where
$\Wh$ is the Whitehead group of $\pi_1M$ with $\ZZ/2$ acting on $\Wh$ by the standard involution,
and on $\ZZ$ by $(-1)^m$. Then we have
\[ (\widetilde\sS^h(M),\sS(M)) \lra (B(\ZZ\otimes_{\ZZ}\Wh),\pt)~ \]
and this is easy to understand. (See remark~\ref{rem-easyunder} below.)
It takes all $0$-simplices in the source
to the base point of $B(\ZZ\otimes_{\ZZ}\Wh)$. It takes a $1$-simplex corresponding to
an $h$-cobordism (over $M$) of torsion $t$ to the loop in
$B(\ZZ\otimes_{\ZZ/2}\Wh)$ determined by $1\otimes t$. Therefore the left-hand square in~(\ref{eqn-maincomponents})
is a pullback square. \qed

\begin{rem} \label{rem-easyunder} {\rm It cannot hurt to say once again that we like to construct $\widetilde\sS^h(M)$
as the geometric realization of a simplicial \emph{space}. When $M$ is closed, the space of $k$-simplices
can be described (for example) as
\[ \hofiber_{M\times\Delta^k}[ B\sC_k \to B\sD_k ]  \]
where
\begin{itemize}
\item $\sC_k$ is the discrete groupoid whose objects are compact manifolds modelled on $\RR^m\times\Delta^k$ (that is,
compact manifolds with an atlas having charts in $\RR^m\times\Delta^k$ and changes of charts respecting
the $(m+k-1)$-dimensional faces), and whose morphisms are homeomorphisms respecting the face structure;
\item $\sD_k$ is a Poincar\ee duality space analogue of $\sC_k$ (still a discrete category, but not strictly
a groupoid).
\end{itemize}
This description is also useful for us because the space of $0$-simplices
is homology equivalent to $\sS(M)$, and so can be used as a good substitute for $\sS(M)$. In the proof
above, where $0$-simplices and $1$-simplices are mentioned, think \emph{space of $0$-simplices} and
\emph{space of $1$-simplices}.
}
\end{rem}

\begin{appendices}
\medskip
\section{Homeomorphism groups of some stratified spaces}
\label{sec-discrete}
When we evaluate characteristic invariants of manifolds on a manifold $M$, we expect and
normally have enough invariance under the discrete group of automorphisms of $M$. It can be very hard
to establish invariance under the topological group of automorphisms of $M$. In \cite{DwyerWeissWilliams}
we relied on the Mather-McDuff-Segal-Thurston theory \cite{McDuff80} which gives homology isomorphisms between the
classifying spaces of the discrete and topological automorphism groups of a compact topological manifold $M$.
We need a similar tool to prove invariance of various controlled characteristics of $M\times\RR^i$ (for compact $M$)
under the topological group of
\emph{controlled} homeomorphisms $M\times\RR^i\to M\times\RR^i$. Recall that $M\times\RR^i$ is identified
with the nonsingular part of the control space $(M*S^{i-1},S^{i-1})$. Three different options come to mind:
\begin{itemize}
\item[(i)] Look for an extension of the Mather-McDuff-Segal-Thurston theory to controlled automorphisms,
hoping that the classifying spaces of the discrete and topological groups of controlled automorphisms
of $M\times\RR^i$ have the same homology.
\item[(ii)] Try to make the Mather-McDuff-Segal-Thurston theory work for automorphisms of the
control space $(M*S^{i-1},S^{i-1})$, and compare with the existing theory for automorphisms
of the singular part $S^{i-1}$.
\item[(iii)] Attempt a reduction to the control-free setting by using the belt buckle trick, which
implies that the classifying space of the topological controlled automorphism group of $M\times\RR^i$ is a homotopy retract of the
classifying space of the topological automorphism
group of $M\times \RR^i/\ZZ^i$.
\end{itemize}
In the end we decided for option (ii). We could not handle (i).
The relationship between (i) and (ii) is nevertheless simple. In (i), the focus is on automorphisms of the control space
$(M*S^{i-1},S^{i-1})$ which extend the identity on the singular part $S^{i-1}$; in (ii), all automorphisms of
$(M*S^{i-1},S^{i-1})$ as a control space are allowed. Option (iii) looks like an attractive shortcut.
We decided against it mainly because it does not respect symmetries of $\RR^i$ such as linear automorphisms. These
symmetries are important to us in view of \cite{WW1}. Another argument against (iii) is that it generates a demand
for belt buckle machinery on the algebraic side, so that it is not really a shortcut.  \newline
The present chapter is therefore a short review of \cite{McDuff80}, the main point being that everything carries over
from a manifold setting to a more general setting where the geometric objects are well-behaved stratified
spaces with manifold strata. We wonder whether a more abstract formulation can be given, perhaps in
sheaf-theoretic language. In any case we assume some familiarity with \cite{McDuff80}.

\medskip
The main result of the section is theorem~\ref{thm-genMcDuff} with $A=\emptyset$.
In section~\ref{sec-altapprox} we spell out what it means for
controlled automorphism spaces and structure spaces.

\begin{defn}
\label{defn-strat} {\rm A \emph{stratified space} is a space $X$
together with a locally finite partition into locally closed subsets
$X_i$~, called the \emph{strata}, such that the closure of each
stratum $X_i$ in $X$ is a union of strata. By an \emph{automorphism} of
a stratified space $X$, we understand a homeomorphism $X\to X$ mapping each
stratum $X_i$ to itself. }
\end{defn}

\begin{defn}
{\rm A stratified space is a $\TOP$ stratified space if it is
paracompact Hausdorff and each stratum $X_i$ is a manifold (of some
dimension $n_i$~, with empty boundary).}
\end{defn}

\medskip
For us, the most important example of a $\TOP$ stratified space is
the join
\[ X=M*S^{i-1} \]
where $M$ is a closed manifold. We partition
this into two strata. One of these is the embedded copy of
$S^{i-1}$. The other stratum, i.e. the complement of $S^{i-1}$ in
$M*S^{i-1}$, can be identified with $M\times\RR^i$. Automorphisms of
$M*S^{i-1}$, with this stratification, must map the sphere $S^{i-1}$
to itself. (This is not a vacuous condition since, for example, $M$
could also be a sphere in which case $M*S^{i-1}$ is homeomorphic to
a sphere.) \newline
Another important type of stratified space with
two strata is furnished by manifolds with boundary. The boundary can
be regarded as one stratum, the complement of the boundary as the
other stratum. \newline
Combining these two examples, one has a canonical
stratification of $M*S^{i-1}$ into three strata when $M$ is a
manifold with nonempty boundary.

\begin{defn}
{\rm The \emph{open cone} $cL$ on a stratified space $L$ is defined
as the quotient of $L\times[0,1[\,$ by $L\times\set{0}$. It is
canonically stratified with strata $cL_i\smin\pt$ and $\pt$, where
$L_i$ denotes a stratum of $L$ and $\pt$ is the base of the cone.}
\end{defn}

The next definition is due to Siebenmann \cite{Siebenmann72}:

\begin{defn}
\label{defn-Siestrat} {\rm A stratified space $X$ is \emph{locally
conelike} if, for each stratum $X_i$ and each $x\in X_i$~, there
exist an open neighborhood $U$ of $x$ in $X_i$~, a compact
stratified space $L$ and a stratification--preserving homeomorphism
(relative to $U$) of $cL\times U$ with an open neighborhood of $x$
in $X$. }
\end{defn}

For a topological manifold $M$ and $i>0$, the stratified space $X=M*S^{i-1}$
with two strata (as explained above) is locally conelike.

\begin{notn} {\rm A locally conelike $\TOP$ stratified space
will be called a $CS$ space. (This is slightly different from
Siebenmann's definition of a $CS$ space, in which there can be only
one stratum of dimension $n$ for each $n\ge 0$, but the partition
into strata is not required to be locally finite.)
\newline
Generalizing some of the conventions of \cite{McDuff80}, we denote
by $\homeo(X)$ and $\Homeo(X)$ the topological group of
automorphisms of a $CS$ space $X$ (with the compact--open topology),
and the underlying discrete group, respectively. More generally, for
a closed subset $A$ of $X$ let $\homeo(X, \rel A)$ and $\Homeo(X,
\rel A)$ be the topological group and the underlying discrete group
of automorphisms of $X$ which agree with the identity in some
\emph{neighborhood} of $A$. Following \cite{McDuff80}, we topologize
$\homeo(X, \rel A)$ as a direct limit
\[ \colim_U \,\set{h\in\homeo(X)\big| h(x)=x\textup{ for }x\in U} \]
where $\set{h\in\homeo(X)\big| h(x)=x\textup{ for }x\in U}$ has the
subspace topology inherited from $\homeo(X)$ and $U$ runs over the
set of all open neighborhoods of $A$ in $X$. (Hence the inclusion map
from $\homeo(X, \rel A)$ to $\homeo(X)$ is continuous, but it need not be
an embedding.) Let
\[ \bar B\homeo(X, \rel A) \]
be the homotopy fiber of the inclusion $B\Homeo(X,\rel A)\to
B\homeo(X,\rel A)$. }
\end{notn}

\begin{lem}
\label{lem-cpctgen} For any compact $K\subset \homeo(X, \rel A)$,
there exists a neighborhood $U$ of $A$ in $X$ such that $f|U=\id$
for all $f\in K$.
\end{lem}

\proof Suppose that for some compact $K\subset \homeo(X, \rel A)$ there is no such $U$.
Then there exists $x\in A$ and a sequence $(y_n)_{n\in\NN}$ in $X$ converging to $x$,
and a sequence $(f_n)_{n\in\NN}$ in $K$ such that $f_n(y_n)\ne y_n$ for all $n$.
Choose a metric $d$ on $X$ inducing the given topology. Let $V_n\subset\homeo(X,\rel A)$ consist
of all $g$ such that $d(g(y_k),y_k)<d(f_k(y_k),y_k)$
for all $k\ge n$. Then $V_n$ is open in $\homeo(X,\rel A)$ by definition of the topology
on $\homeo(X, \rel A)$. Clearly $V_n\subset V_{n+1}$ and
\[  \bigcup_n V_n = \homeo(X,\rel A) \]
so that the $V_n$ for $n\in\NN$ constitute an open covering of $\homeo(X,\rel A)$.
But $f_n\notin V_n$~, so that none of the $V_n$ contains $K$. Therefore $K$ is not compact.
\qed

\begin{defn}
\label{defn-cleansub} {\rm A \emph{clean subspace} of a stratified
space $X$ is a closed subspace $Y$ of $X$ whose frontier $\Fr(Y)$ in
$X$ admits a stratification and a bicollar neighborhood $V\cong
\Fr(Y)\times\RR$ in $X$ (where the homeomorphism $V\cong
\Fr(Y)\times\RR$ respects the stratifications.) }
\end{defn}

\medskip
One example of a clean subspace of a stratified space which we
need is as follows. Let $X$ be a $CS$
space. Let $z$ be a point in a stratum $X_i$, let $U$ be a
neighborhood of $x$ in $X_i$ and let $e\co cL\times U\to X$ be an
open embedding as in definition~\ref{defn-Siestrat}. Let $U'\subset
U$ be an open ball containing $z$ whose closure in $U$ is a
(compact) disk and let $c'L\subset cL$ be an open subcone $L\times
[0,r[\,$ for some $r$ with $0<r<\infty$. We call the image of
$c'L\times U'$ under the embedding $e\co cL\times U\to X$ a
\emph{quasi-ball} (about $z$) in $X$.

\begin{lem}
\label{lem-quasiballclean}
The complement of a quasi--ball in a $CS$ space $X$ is a clean subset of $X$.
\end{lem}

\proof The quasi--ball is contained in an open subset of $X$ which
is identified with $cL\times U$. (We may as well assume $X=cL\times
U$.) Without loss of generality, $U$ can be identified with a
euclidean open ball of radius $1$. Now the spaces $cL$ and $U$ come
with evident actions of the topological monoid $[0,1[$~, a submonoid
of $(\RR,\cdot)$~, and so we get a diagonal action of $[0,1[\,$ on
$cL\times U$. Without loss of generality the quasi--ball is
$\frac{1}{2}\cdot (cL\times U)$. A bicollar for its frontier is then
defined by the embedding
\[ (x,t)\mapsto \frac{e^t}{e^t+1}\cdot x\in cL\times U \]
for $t\in\RR$ and $x$ in the frontier. \qed

\begin{lem} Any point in a $CS$ space
has an open neighborhood which is a quasi--ball. \qed
\end{lem}

The following theorem generalizes a result of Mather in \cite{Mather71} which may have
been the starting point of the theory developed in \cite{McDuff80}. Mather's result serves
in \cite{McDuff80} as an induction beginning in a handle induction, and our generalization has a
similar purpose here.

\begin{thm}
\label{thm-quasiMather} Let $A$ be closed in a $CS$ space $X$.
Assume that the complement of $A$ is a quasi--ball. Then {\rm
$B\Homeo(X,\rel A)$} is acyclic, $B\homeo(X,\rel A)$ is contractible
and hence $\bar B\homeo(X,\rel A)$ is acyclic.
\end{thm}

\proof The contractibility of $B\homeo(X,\rel A)$ is a consequence
of an Alexander trick showing that in fact $\homeo(X,\rel A)$ is
contractible. \newline Our proof that $B\Homeo(X,\rel A)$ is acyclic
is a straightforward generalization of Mather's proof in the
unstratified case \cite{Mather71}. As in the proof of
lemma~\ref{lem-quasiballclean} we can assume that $X= cL\times U$
and that the quasi--ball is $E=\frac{1}{2}\cdot (cL\times U)$. We
are looking at automorphisms $h\co cL\times U\to cL\times U$ which
are the identity outside $r_h\cdot(cL\times U)$ for some positive
$r_h<1/2$ depending on $h$. Copying Mather's argument in the
unstratified case, we begin by choosing a sequence of disjointly
embedded codimension zero disks $D_i\subset \frac{1}{2}U$, for
$i=1,2,3,\dots$. Let $U_i$ be the relative interior of $D_i$ in $U$.
Then $E_i=(2^{-i}\cdot cL)\times U_i$ is a quasi--ball in $E$ for
each $i$, and $E_1$, $E_2$, $E_3$, \dots have disjoint compact
closures in $E$. These choices can easily be made in such a way that
there exists $h_0\in \Homeo(X,\rel X\smin E)=\Homeo(X,\rel A)$ which
maps $E_i$ onto $E_{i+1}$ for $i=1,2,3,\dots$. \newline Let
$G_i=\Homeo(X,\rel X\smin E_i)$. The inclusion
\[  BG_i\to B\Homeo(X,\rel X\smin E) \]
induces a surjection in integer homology. The reasoning is that any
finitely generated subgroup of $\Homeo(X,\rel X\smin E)$ is
conjugate in $\Homeo(X,\rel X\smin E)$ to a subgroup of $G_i$.
(Namely, for any $g_1,g_2,\dots,g_k$ in $\Homeo(X,\rel X\smin E)$
there is $r<1/2$ such that $g_1,\dots,g_k$ agree with the identity
outside $r\cdot(cL\times U)$. Then there exists $g\in \Homeo(X,\rel
X\smin E)$ mapping the closure of $r\cdot(cL\times U)$ to $E_i$~, so
that conjugation with $g$ takes $g_1,\dots g_k$ to the subgroup
$\Homeo(X,\rel X\smin E_i)$.)
\newline
The next step is to introduce the subgroup $G$ of $\Homeo(X,\rel
X\smin E)$ consisting of all automorphisms $h$ which have $h(x)=x$
for all $x\notin \bigcup_iE_i$. Then the restriction homomorphisms
$h\mapsto h|E_i$ lead to isomorphisms
\[
\begin{array}{ccccc}
G & \cong & \prod_{i=1}^{\infty} G_i & \cong & \prod_{i=1}^{\infty}
G_1\,.
\end{array}
\]
(The first isomorphism uses the fact that a bijective continuous map
between compact Hausdorff spaces is a homeomorphism. The second
isomorphism uses conjugation with appropriate powers of $h_0$.) Let
$\sigma\co G\to \Homeo(X,\rel X\smin E)$ be the inclusion. Define
homomorphisms
\[
u,v,w\co G_1\lra G\cong \prod_{i=1}^{\infty} G_1
\]
by $u(g)=(g,1,1,\dots)$, $v(g)=(1,g,g,\dots)$ and
$w(g)=(g,g,g,\dots)$. Using conjugation with the element $h_0$, we
see that $\sigma v$ and $\sigma w$ induce the same homomorphism in
integral homology,
\[ \sigma_*v_*=\sigma_*w_*\co H_*(BG_1) \lra H_*(B\Homeo(X,\rel X\smin E)). \]
Now assume inductively that our vanishing statement has been
established in degrees $s$ for $0<s<k$. In particular
$H_s(B\Homeo(X,\rel X\smin E))=0$ and $H_s(BG_1)=0$ for $0<s<k$.
Writing $w$ as a composition
\[
\CD G_1 @>\textup{diagonal}>> G_1\times G_1 @>\id\times v'>>
G_1\times G'
\endCD
\]
where $G'=\prod_{i=2}^{\infty}$~, and using the K\"{u}nneth formula
for the homology of $BG_1\times BG_1$, we get for any $z\in
H_k(BG_1)$ that
\[ w_*(z)= (z\times 1) + (1\times v'_*(z))\in
H_k(B(G_1\times G')= H_k(BG)\,. \] Hence
$\sigma_*w_*(z)=\sigma_*u_*(z)+\sigma_*v_*(z)$ in
$H_k(B\Homeo(X,\rel X\smin E))$. Since $\sigma_*v_*=\sigma_*w_*$,
this means $\sigma_*u_*(z)=0$ and since $\sigma_*u_*$ is surjective
this implies
\[ H_k(B\Homeo(X,\rel X\smin E))=0.  \qed \]

\medskip
\begin{notn} {\rm Let $Y$ be a clean subspace of a $CS$ space $X$ and let
$A\subset X$ be closed. We assume that $\Fr(Y)\smin A$ has compact
closure in $X$. Generalizing some of the notation in \cite{McDuff80}
again, we write
\[ \homeo(X,Y,\rel A) \]
for the topological submonoid of $\homeo(X,\rel A)$ consisting of
the automorphisms $X\to X$ in $\homeo(X,\rel A)$ which embed $Y$ in
itself. Let $\homeo_0(X,Y,\rel A)$
be the identity component of $\homeo(X,Y,\rel A)$.
The underlying discrete monoid is $\Homeo_0(X,Y,\rel A)$.
We write
\[ \bar B\homeo_0(X,Y,\rel A):= \hofiber[~B\Homeo_0(X,Y,\rel A) \lra B\homeo_0(X,Y,\rel A)~]~. \]
The set--theoretic quotient of $\Homeo_0(X,Y,\rel A)$ by the equivalence relation
\[
\begin{array}{ccc}
h\sim h' & \Leftrightarrow & h=h'\textup{ on some neighborhood of }Y
\end{array}
\]
will be written as $\Emb^X_0(Y,\rel A)$. There is no preferred topology
on this monoid. But let
\[ \emb(Y,\rel A) \]
be the space of all
embeddings $Y\to Y$ which are the identity on some neighborhood of
$Y\cap A$, with the direct limit topology as before; and let
$\emb^X_0(Y,\rel A)$ be the image of the restriction map
\[ \homeo_0(X,Y,\rel A) \lra \emb(Y,\rel A) \]
with the subspace topology induced from $\emb(Y,\rel A)$.
As noted in \cite[\S2]{McDuff80}, this topology agrees with the
quotient topology obtained by viewing $\emb^X_0(Y,\rel A)$ as a
quotient of $\homeo_0(X,Y,\rel A)$. For this McDuff uses an isotopy extension
theorem as in \cite{EdwardsKirby}; we need a variant for $CS$ spaces
\cite[6.5]{Siebenmann72}.
Finally let
\[  \bar B\emb^X_0(Y, \rel A):= \hofiber[~B\Emb^X_0(Y,\rel A) \lra
B\emb^X_0(Y,\rel A)~]. \]
}
\end{notn}


The following propositions~\ref{prop-McDuffs2-1}
and~\ref{prop-McDuffs2-2} as well as
corollary~\ref{cor-unionintersection} generalize propositions 2.1,
2.2 and corollary 2.4 in \cite{McDuff80}, respectively. For the
proofs, see \cite{McDuff80} and the remark following
corollary~\ref{cor-unionintersection} below.

\begin{prop}
\label{prop-McDuffs2-1}
 Let $X$ be a $CS$ space, $Y$ a clean subset of $X$ and
$A$ a closed subset of $X$ such that $\Fr(Y)\smin A$ has compact
closure. Then the inclusion
\[ \bar B\homeo_0(X,Y,\rel A) \lra \bar B\homeo(X,\rel A)  \]
is a weak homotopy equivalence.
\end{prop}

\begin{prop}
\label{prop-McDuffs2-2} If $X$, $Y$ and $A$ are as in
proposition~\ref{prop-McDuffs2-1}, then the sequence
\[
\CD \bar B\homeo_0(X,\rel Y\cup A) @>>> \bar B\homeo_0(X,Y,\rel A)
@>\bar\rho>> \bar B\emb_0^X(Y,\rel A)
\endCD
\]
is an integer homology fibration sequence. In other words, the
inclusion of the space $\bar B\homeo_0(X,\rel Y\cup A)$ into the
homotopy fiber of the restriction map $\bar\rho$ induces an
isomorphism of (untwisted) integer homology groups.
\end{prop}

\begin{cor}
\label{cor-unionintersection} Let $X$, $Y$ and $A$ be as in
proposition~\ref{prop-McDuffs2-1}. If $\bar B\homeo(X,\rel Z)$ is
acyclic for $Z=A$, $Y$ and $Y\cup A$, then it is acyclic for $Z=Y\cap A$
also.
\end{cor}

\medskip
\nin\emph{Remark.} Corollary~\ref{cor-unionintersection} is a formal
consequence of propositions~\ref{prop-McDuffs2-1}
and~\ref{prop-McDuffs2-2}. See \S2 of \cite{McDuff80} for the
deduction.
\newline The proof of McDuff's proposition 2.1 in
\cite{McDuff80}, which is the unstratified case of
proposition~\ref{prop-McDuffs2-1} above, occupies \S4 of
\cite{McDuff80}. Its backbone is a lemma about topological monoids,
lemma 4.1 in \cite{McDuff80}. The deduction of McDuff's proposition
2.1 from that lemma occupies only half a page (right after the
statement of the lemma) and carries over to the stratified case with
only trivial changes. \newline Similarly, the proof of McDuff's
proposition 2.2 in \cite{McDuff80}, which is the unstratified case
of proposition~\ref{prop-McDuffs2-2} above, occupies \S3 of
\cite{McDuff80}. It relies on a string of lemmas about topological
groups and monoids. The deduction of McDuff's proposition 2.2 from
these lemmas occupies only half a page (at the end of \S3 in
\cite{McDuff80}) and carries over to the stratified case with only
trivial changes. \newline

\medskip
\begin{thm} \label{thm-genMcDuff} The space $\bar B\homeo(X,\rel A)$ is acyclic if
$A$ is closed in $X$ and $X\smin A$ has compact closure in $X$.
\end{thm}

\proof Imitating the strategy of \cite{McDuff80}, we deduce this directly
from theorem~\ref{thm-quasiMather} and corollary~\ref{cor-unionintersection}.
Let $Z$ be a stratum of $X$ which is not contained in $A$ and
which is minimal among the strata of $X$ with that property (i.e.,
all strata of $X$ which belong to the closure of $Z$ in $X$, except
$Z$ itself, are contained in $A$). By induction, we may assume:
\begin{itemize}
\item[(i)] $\bar B\homeo(X,\rel N)$ is acyclic for each closed $N\subset X$
which is a neighborhood of $A\cup Z$.
\end{itemize}
In addition, we will assume to begin with that
\begin{itemize}
\item[(ii)] there exists $z\in Z\smin A$ and a quasi--ball
neighborhood $V$ about $z$ which contains all of $Z\smin A$.
\end{itemize}
(This assumption will be removed at a later stage.) Because of (ii),
the open subset $Z\smin A$ of $Z$ is identified with an open subset
of a standard euclidean space, and can be triangulated. We fix a
triangulation.
\newline
For every closed $A'$ which is a neighborhood of $A$ in $X$, we can
choose a triangulation of $Z\smin A$ such that $Z\smin A'$ is
covered by finitely many open stars $\st(z_1),\dots,\st(z_r)$ where
$z_1,\dots,z_r$ are vertices of the triangulation. Making $V$
sufficiently slim, we can arrange that the portion $V_i$ of $V$
lying over $\st(z_i)$ has empty intersection with $A$. Let
$A_i=X\smin V_i$. For $S\subset \set{1,\dots,r}$ let $A_S$ be the
union of the $A_i$ with $i\in S$. Then
\begin{itemize}
\item[(iii)] each $A_S$ is a clean subspace of $X$ and its complement
is a quasi--ball, or the empty set. Hence $\bar B\homeo(X,\rel A_S)$
is acyclic by theorem~\ref{thm-quasiMather}.
\end{itemize}
Now $\bar B\homeo(X,\rel A_1\cap A_2\cap\cdots A_i)$ is acyclic for
$i=1,\dots,r$. This can easily be shown by induction on $i$, using
corollary~\ref{cor-unionintersection} and the last part of (iii)
just above. \newline Finally choose a closed neighborhood $N$ of
$A\cup Z$ such that
\[ N\cap A_1\cap A_2\cap\cdots\cap A_r\quad \subset \quad A'\,. \]
By assumption (i), we also know that $\bar B\homeo(X,\rel N)$ is
acyclic. Hence, by corollary~\ref{cor-unionintersection} again,
\[ \bar B\homeo(X,\rel N\cap A_1\cap A_2\cap\cdots A_r) \]
is acyclic. We have now shown that, for any closed neighborhood $A'$
of $A$, there exists another closed neighborhood $A''$ of $A$ with
$A''\subset A'$ such that $\bar B\homeo(X,\rel A'')$ is acyclic. By
lemma~\ref{lem-cpctgen}, this implies that $\bar B\homeo(X,\rel A)$
is acyclic.
\newline
Now we repeat the argument, but without hypothesis (ii). As before,
we fix a closed neighborhood $A'$ of $A$ in $X$ and we choose
finitely many quasi--balls $V_1,\dots V_r$ about points in $Z\smin
A$ such that the union of the $V_i$ contains $Z\smin A'$ and is
contained in $X\smin A$. We define $A_i$ and $A_S$ as before, for
$S\subset \set{1,\dots,r}$. We have less information about the $A_S$
this time, but at least we know that condition (ii) is satisfied
with $A_S$ in place of $A$. Therefore
\begin{itemize}
\item[(iv)] $\bar B\homeo(X,\rel A_S)$ is acyclic for each
$S\subset \set{1,\dots,r}$,
\end{itemize}
by the first part of this proof, which relied on (ii). We can now
finish the argument as in the first part of the proof, using (ivh
instead of (iii). \qed

\section{Controlled homeomorphism groups}  \label{sec-altapprox}
For a closed manifold $M$, view $M\times\RR^i$ as part of control space $(M*S^{i-1},M\times\RR^i)$. We do not
know whether the inclusion
\[  B\Homeo(M\times\RR^i\cc) \to B\homeo(M\times\RR^i\cc) \]
is a homology equivalence and we have given up on that. Let $\Gamma$ be the homotopy
fiber of $\Homeo(S^{i-1})\to \homeo(S^{i-1})$, a topological group acting on $B\Homeo(M\times\RR^i\cc)$
and on $B\homeo(M\times\RR^i\cc)$ via the forgetful homomorphism to $\Homeo(S^{i-1})$.
The Borel construction $B\Homeo(M\times\RR^i\cc)_{h\Gamma}$ comes with a projection
\[  B\Homeo(M\times\RR^i\cc)_{h\Gamma} \lra B\Gamma \]
which has a zero section $B\Gamma\to B\Homeo(M\times\RR^i\cc)_{h\Gamma}$ since the action of $\Gamma$ respects
the base point of $B\Homeo(M\times\RR^i\cc)$. By the \emph{reduced} Borel construction
\[  B\Homeo(M\times\RR^i\cc)_{rh\Gamma} \]
we mean the quotient of $B\Homeo(M\times\RR^i\cc)_{h\Gamma}$ by the image of the zero section. The quotient map
\[  B\Homeo(M\times\RR^i\cc)_{h\Gamma}\lra B\Homeo(M\times\RR^i\cc)_{rh\Gamma} \]
induces an isomorphism in homology (for any local coefficient system on the target) because
$B\Gamma$ has the homology of a point. Similarly, there is a reduced Borel construction
$B\homeo(M\times\RR^i\cc)_{rh\Gamma}$.

\begin{lem} \label{lem-otherniceapprox}
The inclusion $B\homeo(M\times\RR^i\cc)\to B\homeo(M\times\RR^i\cc)_{rh\Gamma}$
is a homotopy equivalence.
\end{lem}

\proof The action of $\Gamma$ on $B\homeo(M\times\RR^i\cc)$ extends to an action of the
contractible topological group $\hofiber[\homeo(S^{i-1})\to \homeo(S^{i-1})$. Therefore
\[ B\Homeo(M\times\RR^i\cc)_{h\Gamma}~\simeq~ B\Homeo(M\times\RR^i\cc)\times B\Gamma \]
and so the reduced Borel construction $B\Homeo(M\times\RR^i\cc)_{rh\Gamma}$ is homotopy equivalent to the
pushout or homotopy pushout of
\[  B\Homeo(M\times\RR^i\cc)\times B\Gamma \longleftarrow B\Gamma \lra \pt~. \]
Hence the inclusion $B\homeo(M\times\RR^i\cc)\to (B\homeo(M\times\RR^i\cc))_{rh\Gamma}$ induces
an isomorphism in homology for any local coefficient system on the target. By the Seifert-van Kampen theorem, it
also induces an isomorphism on fundamental groups. \qed

\begin{lem} \label{lem-niceapprox} The inclusion of reduced Borel constructions
\[
B\Homeo(M\times\RR^i\cc)_{rh\Gamma} \lra B\homeo(M\times\RR^i\cc)_{rh\Gamma}
\]
is a homology equivalence.
\end{lem}

\proof It is enough to show that the inclusion of ordinary Borel constructions
\[
B\Homeo(M\times\RR^i\cc)_{h\Gamma} \lra B\homeo(M\times\RR^i\cc)_{h\Gamma}
\]
is a homology equivalence. And it is also enough to show that the inclusion
\[
B\Homeo(M\times\RR^i\cc)_{h\Gamma} \lra B\homeo(M\times\RR^i\cc)_{h\Gamma'}
\]
is a homology equivalence, where $\Gamma'=\hofiber[\homeo(S^{i-1})\to \homeo(S^{i-1})]$.
For that we identify
$B\Homeo(M\times\RR^i\cc)_{h\Gamma}$ with the homotopy fiber of the composition
\[
\xymatrix@M=5pt{
B\Homeo((M*S^{i-1},M\times\RR^i)) \ar[r]^-{\textup{forget}} & B\Homeo(S^{i-1}) \ar[r]^-{\textup{incl.}} & B\homeo(S^{i-1})~.
}
\]
That composition fits into a commutative square
\[
\xymatrix{
B\Homeo((M*S^{i-1},M\times\RR^i)) \ar[d] \ar[r] &  B\homeo(S^{i-1}) \ar[d]^= \\
B\homeo((M*S^{i-1},M\times\RR^i)) \ar[r] & B\homeo(S^{i-1})
}
\]
whose left-hand column is a homology isomorphism by section~\ref{sec-discrete}. Given that
\[ \pi_1B\homeo(S^{i-1})=\pi_0\homeo(S^{i-1})\cong \{\pm1\} \]
is as uncomplicated as it is, it is not hard to deduce
a homology isomorphism between the homotopy fibers of the rows. This is exactly what we need
since the homotopy fiber in the lower row is $B\homeo(M\times\RR^i\cc)_{h\Gamma'}$. \qed

\begin{expl} {\rm Let $M$ be a compact manifold such that $\partial M$ is the union of two codimension
zero submanifolds $\partial_0M$ and $\partial_1M$ with intersection
$\partial_0M\cap \partial_1M=\partial\partial_0M=\partial\partial_1M$. What is the recommended
approximation (by the classifying space of a discrete group etc.) to
$B\homeo(M\times\RR^i,\partial_1M\times\RR^i\cc)$? We also write
$B\homeo((M,\partial_1M)\times\RR^i\cc)$ for $B\homeo(M\times\RR^i,\partial_1M\times\RR^i\cc)$.
\begin{itemize}
\item[(i)] In the case where $\partial M=\emptyset$, we can use $B\Homeo(M\times\RR^i\cc)_{rh\Gamma}$.
\item[(ii)] In the case where $i=0$, we can use $B\Homeo(M \rel \partial_0M)$. This is covered
by theorem~\ref{thm-genMcDuff}, or the unstratified version in \cite{McDuff80}. Note however that
$\Homeo(M\rel \partial_0M)$ consists of homeomorphisms $M\to M$ which agree with the
identity \emph{on some neighborhood} of $\partial_0M$.
\item[(iii)] In the general case, we resort to approximating the forgetful map
\[ \quad B\homeo((M,\partial_0M,\partial_1M)\times\RR^i\cc) \lra
B\homeo((\partial_0M,\partial\partial_0M)\times\RR^i\cc)~ \]
by the forgetful map
\[\qquad\quad\,\, B\Homeo((M,\partial_0M,\partial_1M)\times\RR^i\cc)_{rh\Gamma} \lra
B\Homeo((\partial_0M,\partial\partial_0M)\times\RR^i\cc)_{rh\Gamma}. \]
\end{itemize}
}
\end{expl}

\section{$K$-theory of pairs and diagrams} \label{sec-Kpairs}
Let $f\co \sW_s\to \sW_t$ be an exact functor between Waldhausen categories.
Make a new Waldhausen category $\sP=\sP(\sW_s\to \sW_t)$ whose objects are triples $(a_s,a_t,g)$ where
$a_s$ and $a_t$ are objects of $\sW_s$ and $\sW_t$, respectively, and $g\co f(a_s)\to a_t$
is a cofibration in $\sW_t$. A morphism in $\sP$ from $(a_s,a_t,g)$ to $(b_s,b_t,h)$
is a pair of morphisms $(q_s\co a_s\to b_s,q_t\co a_t\to b_t)$ such that the square
\[
\CD
f(a_s) @>g>> a_t \\
@VV q_s V   @VV q_t V \\
f(b_s) @>h>> b_t
\endCD
\]
in $\sW_1$ commutes. Call $(q_s,q_t)$ a weak equivalence if $q_s$ and $q_t$ are weak equivalences; call it a
cofibration if $q_s$ and $q_t$ and the induced map from the pushout of
\[
\CD
f(a_s) @>g>> a_t \\
@VV q_0 V   @. \\
f(b_s) @.
\endCD
\]
to $b_t$ are cofibrations. The additivity theorem implies immediately that
\[  \bK(\sP)~\simeq~\bK(\sW_s)\times \bK(\sW_t)~\simeq~\bK(\sW_s)\vee \bK(\sW_t) \]
by means of the coordinate functors $(a_s,a_t,g)\mapsto a_s$ and $(a_s,a_t,g)\mapsto a_t$. \newline
If both $\sW_s$ and $\sW_t$ are equipped with SW products, denoted $\odot$, and the functor $f$ respects these,
then $\sW_f$ has a preferred SW product given by
\[   (a_s,a_t,g) \odot (b_s,b_t,h):= \hofiber[ g\odot h\co a_s\odot b_s \lra a_t\odot b_t~]~. \]
The splitting $\bK(\sP)\simeq \bK(\sW_s)\times \bK(\sW_t)$ may not respect the
involution on $\bK(\sP)$ determined by the SW product on $\sP$. But there is still a
homotopy fiber sequence of spectra with involution,
\begin{equation} \label{eqn-pairKsequence}
\xymatrix{
\bK(\sW^{(1)}_t) \ar[r] & \bK(\sP) \ar[r]^-{\textup{forget}} &  \bK(\sW_s)
}
\end{equation}
where $\sW_t^{(1)}$ is $\sW_t$ with a new SW product, defined in terms of the
original one by $(a_t,b_t)\mapsto \Omega(a_t\odot b_t)$.

\medskip
In the case where $(X,Y)$ is a
pair of spaces, $\sW_s$ and $\sW_t$ are certain categories of retractive spaces over $Y$ and $X$
which we use to define the spectra $\bA(Y)$ and $\bA(X)$, respectively, and $f\co \sW_s\to \sW_t$ is the functor
induced by the inclusion $Y\to X$, we may write $\bA(Y\subset X)$ or $\bA(Y\to X)$
instead of $\bK(\sP)$.
Furthermore, if $X$ comes with a spherical fibration $\nu$, then we have
$n$-duality SW products on $\sW_s$ and $\sW_t$ determined by $\nu$ and $n$, and we
may use notation such as
\begin{equation}
\xymatrix@C=35pt{
\bA(X,\nu,n+1) \ar[r] & \bA(Y\subset X,\nu,n+1) \ar[r]^-{\textup{forget}} &  \bA(Y,\nu,n)
}
\end{equation}
as an alternative to the notation in~(\ref{eqn-pairKsequence}). \newline
These definitions generalize mechanically to the situation where we have a commutative square of spaces
and inclusion maps,
\[
\xymatrix@R=10pt@C=10pt{ X_{01} \ar[r] \ar[d] & \ar[d] X_0 \\
X_1 \ar[r] & X
}
\]
giving rise to a square of exact functors between categories of retractive spaces,
\[
\xymatrix@R=10pt@C=10pt{ \sW_{01} \ar[r] \ar[d] & \ar[d] \sW_0 \\
\sW_1 \ar[r] & \sW~.
}
\]
Namely, we let
\[  \bA\left(\begin{aligned}\xymatrix@R=8pt@C=8pt{ X_{01} \ar[r] \ar[d] & \ar[d] X_0 \\
X_1 \ar[r] & X
}\end{aligned}\right) := \bK(\sP(\sP_u\to \sP_\ell))
\]
where $\sP_u$ and $\sP_\ell$ are the Waldhausen category of pairs determined by $\sW_{01}\to \sW_0$
and $\sW_1\to \sW$, respectively. If $X$ comes equipped with a spherical fibration $\nu$, then there is a homotopy
fiber sequence of spectra with involution
\begin{equation}
\begin{split}
\xymatrix@R=15pt@M=10pt{
{\bA(X_1\to X,\nu,n+1)} \ar[d] \\
{\bA\left(\CD X_{01} @>>> X_0 \\ @VVV @VVV \\
X_1 @>>> X \endCD,\nu,n+1\right)} \ar[d] \\
{\bA(X_{01}\to X_0,\nu,n).}
}
\end{split}
\end{equation}

\section{Corrections and Elaborations}
\label{sec-updates}
\begin{rem}
\label{rem-corrduality} {\rm There is an unfortunate oversight in
\cite{WWduality}, as follows. In \cite[7.1]{WWduality} we have an
enlarged model $xK(\sC)$ of the $K$--theory space $K(\sC)$ of a
Waldhausen category $\sC$ with Spanier--Whitehead product,
\[ xK(\sC) = \Omega|xw\sS\lbul\sC|\,. \]
Here $xw\sS\lbul\sC$ is an enlarged model of Waldhausen's
$w\sS\lbul\sC$. It is a simplicial category with a degreewise
involution which \emph{anticommutes} with the simplicial operators.
The involution at the category level induces an involution on
$|xw\sS\lbul\sC|$. It should have been pointed out just before
\cite[7.2]{WWduality} that this must be combined with the ``reverse
loops'' operation to give the preferred involution on
$\Omega|xw\sS\lbul\sC|=xK(\sC)$. With that convention, the standard
inclusion of $|xw\sC|$ in $\Omega|xw\sS\lbul\sC|$ respects the
preferred involutions. This is used in \cite[\S9]{WWduality}. }
\end{rem}

\begin{rem}
\label{rem-Vogelstuff} {\rm We have frequently encountered the following
constellation in this paper: a Waldhausen category $\sD$ with $SW$-product
satisfying all the usual axioms and a Waldhausen subcategory $\sC\subset \sD$
closed under weak equivalences and ``duals''. Then it is often useful to
have something like a ``quotient'' of $\sD$ by $\sC$, designed in such a way
that the algebraic $K$-theory spectrum of the quotient is homotopy
equivalent to the mapping cone of $\bK(\sC)\to \bK(\sD)$, and similarly
for the various $L$-theory spectra. From the point of
view of algebraic $K$-theory the easiest and best approach is to continue using
$\sD$, but with a new notion of weak equivalence where all morphisms in $\sD$ whose
mapping cones are in $\sC$ qualify as weak equivalences. From an $L$-theory
point of view, this seems (at first) less fortunate because the old $SW$
product in $\sD$, with the new notion of weak equivalence, will normally
violate one of the basic conditions for an $SW$ product \cite[1.1]{WWduality}.
It is probably possible to repair this by introducing a new $SW$ product to go
along with the new notion of weak equivalence in $\sD$. But this may not always be
worthwhile. What we have tended to do instead is to continue working with the
old $SW$ product in $\sD$, and to look for pairings in
$\sD$ which are nondegenerate for the new notion of weak equivalence (i.e.,
nondegenerate modulo $\sC$). This is unproblematic, but it is slightly embarrassing
that the situation has not been looked at in detail in \cite{WWduality}. It is all
the more reassuring that Ranicki in \cite{RanickiTopMan} has got it right~;
but Ranicki has it only in the setting of chain complex categories.
}
\end{rem}

\begin{rem}
\label{rem-corr} {\rm The proof of
proposition 4.3 in \cite{Weissexci} is slightly
wrong. (The ``functor'' $\tau$ as defined there is not a functor;
its ``induced morphisms'' do not satisfy the required control
condition.) Here is a corrected proof. For each integer $i\ge 0$,
define $\psi_i\co [0,1]\to \,]0,1]$ so that
$\psi_i(t)=(1-\log_2(t))/(1+i)$ if $t\ge 2^{-i}$ and $\psi_i(t)=1$
otherwise. There is an endomorphism of the control space
$(Z\times[0,1],Z\times[0,1[\,)$ given by the formula
\[ ((x,s),t) \mapsto ((x,s\cdot \psi_i(t)),t) \]
for $(x,s)\in Z=X\times\,]0,1]$. Denote the induced endofunctor of
$\sA^{\rho}(\JJ Z)^{\wedge}_{\infty}$ by $\sigma_i$. Note that
$\sigma_0=\id$~, and all the $\sigma_i$ are
related by invertible natural transformations (so that
$\sigma_i$ is isomorphic to $\sigma_j$ for any $i,j\ge 0$).
We re-define $\tau$ by the formula
\[ \tau(A)=\bigoplus_{i\ge 0} \sigma_i\]
and this time it is an endofunctor of $\sA^{\rho}(\JJ
Z)^{\wedge}_{\infty}$. There is an Eilenberg swindle in the shape of
a natural isomorphism of functors $\tau\cong \id\oplus\tau$. This
uses the natural isomorphisms of functors $\sigma_i\to
\sigma_{i+1}$ mentioned earlier. Hence, for the self-map $\tau_*$
of the infinite loop space $K(\sA^{\rho}(\JJ Z)^{\wedge}_{\infty})$
we have $\tau_*+\id\simeq\tau_*$. }
\end{rem}

\medskip
\begin{rem}
\label{rem-Spivak} {\rm The existence and uniqueness of Spivak normal fibrations for a
Poincar\ee duality space $X$ of formal dimension $n$ has been repeatedly used in this paper.
What does it mean~? Let $G_k$ be the space of based homotopy automorphisms
of $S^{k-1}$. Let $\xi_k$ be the canonical quasi-fibration on $BG_k$ with fibers $\simeq S^{k-1}$.
Form the space $U_k(X)$ of pairs $(g,\eta)$ where $g\co X\to BF_k$
and $\eta$ is a based map from $S^{n+k}$ to the Thom space of $g^*\xi_k$.
(Here, \emph{Thom space} means the mapping cone of the projection.)
The existence and uniqueness claim is that
\[ U_{\infty}(X) = \colim_k U_k(X) \]
is contractible. Wall \cite{WallPoincare} shows that $U_{\infty}(X)$ is connected,
and it seems likely that similar (Spanier-Whitehead duality) arguments could be employed
to show that $U_{\infty}(X)$ is contractible. A different argument is as follows. Suppose
first that $X$ has the homotopy type of a compact $CW$-complex.
For $k\gg 0$ we consider pairs $(N,e)$
where $N$ is a compact smooth codimension zero submanifold of $\RR^{n+k}$,
the inclusion $\partial N\to N$ induces in isomorphism on $\pi_1$~, and
$e\co N\to X$ is a homotopy equivalence. We can call such a thing a regular neighborhood
of $X$ in $\RR^{n+k}$. These regular neighborhoods can be regarded as $0$-simplices
of a suitable simplicial set where the j-simplices are certain regular
neighborhoods of $X\times\Delta^j$ in $\RR^{n+k}\times\Delta^j$. Let $V_k(X)$ be its
geometric realization. It is relatively easy to verify that
\[ V_{\infty}(X) = \colim_k V_k(X) \]
is contractible; this amounts to an existence and uniqueness statement for
regular neighborhoods of $X$ in euclidean space. Hence it is enough
to show that there exist compatible homotopy equivalences
\[ V_k(X) \to U_k(X) \]
for large enough $k$. Indeed,
for any $(N,e)$ in $V_k(X)$, the homotopy fibers of $e_{|\partial N}\co \partial N\to X$
are $(k-1)$-spheres \cite{RanickiLMS1}. Then $N/\partial N$ can be regarded as
the Thom space of a spherical fibration on $X$ and the Pontryagin-Thom collapse
from $\RR^{n+k}\cup\infty$ to $N/\partial N$ is a Spivak reduction. This gives
maps $V_k(X)\to U_k(X)$. To understand why
these maps $V_k(X)\to U_k(X)$ should be highly connected, fix some $(g,\eta)$
in $U_k(X)$. The mapping cone of $g^*\xi_k$ is a quotient of the mapping
cylinder. The mapping cylinder comes with a map $\delta$ to $[0,1]$ measuring the ``distance'' to
the zero section $X$. If $\eta$ is transverse to $\delta^{-1}(1/2)$, then $\eta^{-1}$
of $\delta^{-1}$ of $[0,1/2]$ is a codimension zero smooth compact submanifold $N$ of $S^{n+k}$
avoiding the base point of $S^{n+k}$. This comes with an obvious map $e\co N\to X$, and
more importantly, with a map of pairs $\bar e$ from $(N,\partial N)$ to the
(disk bundle, sphere bundle) pair associated with $g^*\xi_k$.
Now $e$ may not be a homotopy equivalence. But embedded surgery on $(N,\partial N)$,
with the goal of making $\bar e$ embedded bordant to a homotopy equivalence of pairs,
will repair that. Hence the map $V_k(X)\to U_k(X)$ is $0$-connected. A parameterized version
of this argument shows that $V_k(X)\to U_k(X)$ is $j$-connected for any $j$. \newline
If $X$ is finitely dominated with nonzero finiteness obstruction, then $X\times S^1$
has zero finiteness obstruction, thanks to M. Mather; see also \cite{WallPoincare}.
But $R_{\infty}(X)$ is a homotopy retract of $R_{\infty}(X\times S^1)$; hence
$R_{\infty}(X)$ is again contractible.
}
\end{rem}

\end{appendices}


\end{document}